%% file: DualPolynomials.tex
\documentclass[preprint,3p,sort&compress,times]{elsarticle}
\usepackage{amsmath}
\usepackage{amssymb}
\usepackage{amsthm}
\usepackage{mathtools}
\usepackage{placeins}
\usepackage{bm}
\usepackage{caption}
\usepackage{subcaption}
\usepackage[usenames,dvipsnames]{color}

\usepackage{cases}

\usepackage{tikz-cd}

\usepackage{graphicx}			
\usepackage[all]{xy}

\renewcommand{\d}{\mathrm{d}}

\newcommand \bdof {\VJ{\mathbb{N}_1 \widetilde{\mathcal{B}}^0\bb{\hat{\phi}^h}}}

\newcommand{\tabref}[1]{Table~\ref{#1}} 
\newcommand{\figref}[1]{Figure~\ref{#1}} 
\newcommand{\secref}[1]{Section~\ref{#1}} 
\newcommand{\exampleref}[1]{Example~\ref{#1}} 

\newcommand{\dualgrad}{\widetilde{\mathrm{grad}}}
\newcommand{\dualcurl}{\widetilde{\mathrm{curl}}}
\newcommand{\dualdiv}{\widetilde{\mathrm{div}}}
\newcommand{\dualSb}{\widetilde{S}\bb{\Omega} \times \widetilde{D}\bb{\partial \Omega} }
\newcommand{\dualDb}{\widetilde{D}\bb{\Omega} \times \widetilde{C}\bb{\partial \Omega}}

\newcommand{\bb}[1]{\left( #1 \right)}
\newcommand{\bs}[1]{\left[ #1 \right]}
\newcommand{\nl}{\\[1.1ex]}

\newcommand{\VJ}[1]{\textcolor{black}{#1}}
\newcommand{\MIGnew}[1]{\textcolor{black}{#1}}
\newcommand{\varun}[1]{\textcolor{black}{#1}}

\input{math_commands.tex}

\theoremstyle{plain}

\newtheorem{corollary}{Corollary}

\newtheorem{definition}{Definition}
\newtheorem{example}{Example}

\newtheorem{lemma}{Lemma}

\newtheorem{proposition}{Proposition}
\newtheorem{remark}{Remark}

\journal{}
\begin{document}
\begin{frontmatter}

\title{{Construction and application of algebraic dual polynomial representations for finite element methods on quadrilateral and hexahedral meshes}}
\author[TUD]{V.~Jain\corref{cor}}
\ead{V.Jain@tudelft.nl}
\author[TUD]{Y.~Zhang}
\author[TUD]{A.~Palha}
\author[TUD]{M.~Gerritsma}

\address[TUD]{Delft University of Technology, Faculty of Aerospace Engineering, P.O. Box 5058, 2600 GB Delft, The Netherlands}

\cortext[cor]{Corresponding author. Tel. +31 15 2789670.}

\begin{abstract}
Given a sequence of finite element spaces which form a de Rham sequence, we will construct a dual representation of these spaces with associated differential operators which connect these spaces such that they also form a de Rham sequence.
The matrix which converts primal representations to dual representations -- the Hodge matrix -- is the mass or Gram matrix. It will be shown that a bilinear form of a primal and a dual representation is equal to the vector inner product of the \varun{expansion} coefficients (degrees of freedom) of both representations. This leads to very sparse system matrices, even for high order methods. The derivative of dual representations will be defined. Vector operations, grad, curl and div, for primal and dual representations are both topological and do not depend on the metric, i.e. the size and shape of the mesh or the order of the numerical method.
Derivatives are evaluated by applying sparse incidence and \varun{inclusion} matrices to the expansion coefficients of the representations. As illustration of the use of dual representations, the method will be applied to i) a mixed formulation for the Poisson problem in 3D, ii) it will be shown that this approach allows one to preserve the equivalence between Dirichlet and Neumann problems in the finite dimensional setting and, iii) the method will be applied to the approximation of grad-div eigenvalue problem on \VJ{affine and} non-affine meshes.
\end{abstract}


\end{frontmatter}

\section{Introduction}
The differential operators grad, curl and div play a fundamental role in the representation of \varun{p}hysical field laws. In terms of differential geometry these operators are either represented by a topological, metric-free operation called the exterior derivative or by the metric-dependent operation called the codifferential or coderivative denoted by $\delta$ or $\mathrm{d}^{\star}$, see \cite{frankel}. The codifferential is the Hilbert adjoint of the exterior derivative and the metric dependence of the codifferential is a direct consequence of the definition of inner product, which involves the metric.
Although the distinction between exterior derivative and codifferential is most clearly seen in terms of differential geometry, it can also be introduced in terms of vector calculus as weak differential operators, see for instance \cite[Appendix A]{BG09}. 
		
The construction of the exterior derivative in a discrete setting is relatively straightforward since it  follows from the continuous definition. 
		The challenge lies in the construction of the sequence of discrete function spaces such that, together with the exterior derivative, \VJ{they constitute} an exact sequence (a de Rham sequence), see for example \cite{AFW10}. 
		
		The subject of constructing discrete representations for the codifferential operator is a more challenging one. This operator contains all metric and material properties and its discretization typically establishes the distinguishing properties of each numerical method. Different aspects of the codifferential operator can be taken into account when constructing its discrete counterpart. Apart from its metric-dependence, the discrete codifferential is not just the restriction of the continuous operator to an appropriate finite dimensional subspace. In weak formulations this generally corresponds to using integration by parts to convert the codifferential into an exterior derivative. This approach stems directly from the fact that finite element formulations typically use one de Rham sequence (a primal sequence of spaces).
		By employing integration by parts, the discrete codifferential can be constructed using only one single sequence of spaces, see for example \cite{AFW10}.
		
\varun{In the finite difference/finite volume community the definition of the codifferential operator follows a different route, relying on the construction of a dual grid, see~\cite{LMS14}.}
\cite{LMS14}. 
For the evaluation of the \varun{codiferential}, the scalar or vector field is converted to the dual grid, there the exterior derivative is computed after which the result is mapped back to the primal grid. The mappings which convert the variables on the primal grid to the variables on the dual grid are discrete representations of the Hodge-$\star$ operator and therefore are called Hodge matrices. This approach, instead of defining the discrete codifferential as the dual of the exterior derivative, directly approximates the explicit representation of the codifferential
		\begin{equation*}
			\mathrm{d}^{\star} = (-1)^{d(k+1)+1}\star\mathrm{d}\star
		\end{equation*}
		where $d$ denotes the dimension of the ambient space, $k$ refers to the $k$-form for which the derivative is taken, and $\star$ denotes the Hodge-$\star$ operator see \cite{frankel, abraham, framework2011}.
		
		In \cite{Palha2014} both approaches are compared for a spectral element method.
		On a single grid, integration by parts is used to convert the codifferential to an exterior derivative, while in the dual grid approach dual spectral element basis functions are constructed on a dual grid.
		
		As mentioned before, the construction of discrete codifferential operators is challenging, especially for high order spectral element methods. The underlying reason has to do with the fact that the codifferential is non-local. The higher the approximating degree of the method, the denser the matrix representation of the codifferential  and the higher the condition number of the resulting system of equations become\MIGnew{s}. These two properties are particularly impactful to any numerical method. This is the focus of this work: to propose a high order discrete codifferential operator with improved sparsity and condition number.

In this paper we use integration by parts and interpret the mass or Gram matrices as Hodge matrices.
These Hodge matrices correspond to discrete Hodge-$\star$ operators that convert the primal representation to a dual representation.
By doing so, the integration by parts formula becomes metric-free and the codifferential is interpreted as the derivative of dual variables.
This leads to much sparser system matrices -- especially for higher order methods.
We set up a de Rham sequence for primary spaces and for each space in this complex we construct a dual space.
{The sequence of these dual spaces also forms a de Rham complex.}
The construction of a dual basis used in this work is similar to the \emph{inverse Gram} constructions described in~\cite{2001Oswald, 2013Wozny}.

Earlier dual bases have been used in isogeometric methods for projection of B-splines~\cite{1973deBoor,1976deBoor}.
In finite elements they have been used for mortar methods~\cite{2000Wohlmuth, 2001Oswald}.
In \cite{2017Dornisch} these ideas have been combined to form isogeometric mortar methods.
For other implementations in isogeometric methods see~\cite{2014Veiga,2019Pinto}.
Different {approaches} for construction of dual basis have also been discussed in~\cite{2013Wozny}, and for construction of dual splines in~\cite{2014Wozny}.

As applications, we will demonstrate the use of dual basis on a constrained minimization problem of the Poisson equation in 3D.
We also show the discrete well posedness of this problem which is expressed in terms of degrees of freedom only.
It will be shown that the use of an algebraic dual basis results in a very sparse matrix where two of the sub-matrices consist of $1$, $-1$ and $0$ only and do not change with the shape and size of \VJ{the elements of the mesh as long as the mesh topology remains the same}.
This observation is relevant for the incompressible flow equation where we encounter a similar div-grad pair.
These techniques may also be valuable in electromagnetism to represent the involution constraint $\mbox{div}\,B=0$ in a way that is very sparse and independent of the shape and size of the mesh.

We will also solve for the pair of Dirichlet-Neumann problems discussed in \cite{CarstensenDemkowiczGopalakrishnan} and prove their equivalence in discrete sense.
This is in general not trivial.
It is shown that the duality relation continues to hold point-wise in these finite-dimensional approximations, on arbitrary grids, through the use of algebraic dual polynomials.

\VJ{Finally, we will solve for a grad-div eigenvalue problem on affine and non-affine meshes using dual degrees of freedom and will show optimal convergence rates for the eigenvalues.}

\VJ{The remainder of this paper is structured as follows:}
The construction of dual polynomial spaces in the one dimensional case is presented in Section~\ref{sec:Construction_of_dual}. In this section it is also shown how nodal sampling and edge sampling from polynomial spaces extend to Sobolev spaces.
The treatment of multiple elements case and derivative of a dual representation will be given in this section.
In Section\MIGnew{s}~\ref{sec:2D_dual_spaces} and \ref{sec:3D_dual_spaces} these constructions are extended to \VJ{two dimensions and three dimensions}\MIGnew{,} respectively.
In Section~\ref{sec:mixed_formulation} a dual polynomial representation is used for the mixed formulation of the Poisson equation in the three-dimensional case with multi-elements.
In Section~\ref{sec:Dirichlet_Neumann} equivalence of the Dirichlet-Neumann problems~\cite{CarstensenDemkowiczGopalakrishnan} is proved and demonstrated by a particular example.
In \secref{sec:eigenvalues2} we address \varun{a} grad-div eigenvalue problem and show optimal convergence rates on affine and non-affine grids.
Finally, in Section~\ref{sec:Conlusions} conclusions are drawn and future work is discussed.
\section{Construction of  {1D} dual finite elements}\label{sec:Construction_of_dual}
We will use the definition of finite element spaces in terms of the triplet $(K,\mathcal{P},\mathcal{N})$ by Ciarlet, \cite{Ciarlet78}, see also Ern and Guermond, \cite[\S1.2]{ErnGuermond04} and Brenner and Scott, \cite[\S3.1]{BrennerScott08}.

\begin{definition}\label{def:finite_element}
A finite element consists of the triplet $(K,\mathcal{P},\mathcal{N})$ where
\begin{enumerate}[i]
\item $K$ is a compact, connected, Lipschitz subset of $\mathbb{R}^d$ with non-empty interior;
\item $\mathcal{P}$ is a (finite dimensional) linear vector space with domain $K$. Usually, $\mathcal{P}$ is a polynomial vector space \varun{of dimension $d_{\mathcal{P}}$};
\item $\mathcal{N}$ is a set of linear functionals $\{ \mathcal{N}_i\}$, $i=1,\ldots,d_{\mathcal{P}}$, acting on elements of $\mathcal{P}$, such that the linear map,
    \begin{equation*}
    p^h \in \mathcal{P} \; \mapsto \; (\mathcal{N}_1(p^h),\ldots,\mathcal{N}_{d_{\mathcal{P}}}(p^h)) \in \VJ{\mathbb{R}^{d_{\mathcal{P}}}} \;,
    \label{eq:bijectivity_linear_functionals}
    \end{equation*}
    is bijective.
\end{enumerate}
\end{definition}
The linear functionals $\{\mathcal{N}_i\}$ are called the {\em local degrees of freedom}.
The following Proposition taken from \cite{ErnGuermond04} defines the basis functions:
\begin{proposition}\label{prop:basis_functions}
There exists a basis $\{ \Psi _1,\ldots,\Psi _{d_{\mathcal{P}}} \}$ in $\mathcal{P}$ such that
\[ \mathcal{N}_i(\Psi _j) =\delta_{ij}\;,\quad 1\leq i,j\leq d_{\mathcal{P}} \;.\]
\end{proposition}
\begin{example}\label{ex:Lagrange_1D}
Consider the interval $K = [-1,1] \subset \mathbb{R}$. Let $\xi_i \in K$, $i=0,\ldots,N$, be the roots of the polynomial $(1-\xi^2)L_N'(\xi)$, where $L_N(\xi)$ is the Legendre polynomial of degree $N$ and $L_N'(\xi)$ its derivative. These nodes are referred to as the Gauss-Lobatto-Legendre (GLL) points, \cite{SpectralBook}. Let $\mathcal{P}$ be the space of polynomials of degree $N$ defined on the interval $K$. For any $p^h \in \mathcal{P}$ define the degrees of freedom by
\begin{equation}
\mathcal{N}^0_i(p^h) := p^h(\xi_i) \;,\quad i=0,\ldots,N \;.
\label{eq:nodal_sampling_1D}
\end{equation}
Because polynomials are continuous, \eqref{eq:nodal_sampling_1D} is well-defined.
The superscript `$0$' in $\mathcal{N} _i ^0$ indicates that we sample the polynomial $p^h$ in points.
The basis which satisfies the Kronecker-delta property from Proposition~\ref{prop:basis_functions} is given by the set of Lagrange polynomials through the GLL-points
\[ h_i(\xi) = \frac{(\xi^2-1)L_N'(\xi)}{N(N+1)L_N(\xi_i)(\xi-\xi_i)} \;, \quad i=0, 1, \hdots, N \;.  \]
\end{example}
\noindent
This example also corresponds to \cite[Prop.1.34]{ErnGuermond04} for $d=1$.

\begin{remark}\label{rem:nodal_sampling_outside_P}
Note that the degrees of freedom are linear functionals on $\mathcal{P}$. The nodal sampling of functions in $\mathcal{P}$ is essentially the Dirac delta distribution which is well defined when the space $\mathcal{P}$ consists of continuous functions, see \cite[Example 2.10.2]{DemkowiczOden}.
Extension of this functional to Sobolev spaces in this way is in general not possible.
The extension to Sobolev spaces will be given in Definition \ref{def:nodal_sammpling_L2}.
\end{remark}

\begin{example}\label{ex:edge_1D}
Let $K$ and $\xi_i$ be defined as in Example~\ref{ex:Lagrange_1D}. Let $\mathcal{Q}$ be the space of polynomials of degree $(N-1)$.
{For this example we choose to define the degrees of freedom as}
\begin{equation}
\mathcal{N}^1_i(p^h)  {:=} \int_{\xi_{i-1}}^{\xi_i} p^h(\xi)\,  {\mathrm{d}\xi} \;,\quad i=1,\ldots,N\;.
\label{eq:edge_sampling_1D}
\end{equation}
For polynomials, the integral in \eqref{eq:edge_sampling_1D} is well-defined. The superscript `$1$' in $\mathcal{N}^1_i$ expresses the fact that the degrees of freedom are associated to line segments $[\xi_{i-1},\xi_i]$ {(geometrical dimension $d =1$)}.
The basis functions, $e_j(\xi)$, which must satisfy the Kronecker-delta property from Proposition~\ref{prop:basis_functions} need to satisfy
\[ \mathcal{N}^1_i(e_j) = \int_{\xi_{i-1}}^{\xi_i} e_j(\xi)\,\mathrm{d}\xi = \delta_{ij} \;.\]
\end{example}
\begin{lemma}
The basis functions $e_j(\xi)$ on the GLL-grid defined in Example~\ref{ex:edge_1D} are given by
\begin{equation}
e_j(\xi) = - \sum_{k=0}^{j-1} \frac{\mathrm{d}h_k}{\mathrm{d}\xi}(\xi) \;, \quad j=1,\ldots,N \;,
\label{def:edge_functions}
\end{equation}
where $h_k(\xi)$ are the Lagrange polynomials through GLL points defined in Example~\ref{ex:Lagrange_1D}.
\begin{proof}
\begin{equation*}
\int_{\xi_{i-1}}^{\xi_i} e_j(\xi)\,\mathrm{d}\xi = - \sum_{k=0}^{j-1} \int_{\xi_{i-1}}^{\xi_i}  {\frac{\mathrm{d}h_k}{\mathrm{d}\xi} (\xi) \, \mathrm{d}\xi} = - \sum_{k=0}^{j-1} \left [ h_k(\xi_i) - h_k(\xi_{i-1}) \right ] =\delta_{ij} \;,
\label{eq:Kronecker_delta_edge}
\end{equation*}
where we repeatedly use the Kronecker-delta property of the Lagrange polynomials. If the Lagrange polynomials $h_k(\xi)$ are polynomials of degree $N$, then $\mathrm{d}h_k(\xi)/\mathrm{d}\xi$ is a polynomial of degree $(N-1)$.
It is easy to show that \VJ{$\left\lbrace e_i \right\rbrace _{i=1}^N $} forms a basis for $\mathcal{Q}$.
\end{proof}
\end{lemma}

\begin{corollary}\label{cor:derivative_incidence}
From \eqref{def:edge_functions} it follows that
\[ \frac{\mathrm{d}h_j}{\mathrm{d}\xi} = e_j(\xi) - e_{j+1}(\xi) \;.\]
So if $p^h \in \mathcal{P}$ is expanded in terms of Lagrange polynomials as
\[ p^h(\xi) = \sum_{i=0}^N \mathcal{N}^{0}_i(p^h) h_i(\xi) \;,\]
then its derivative is given by
\begin{equation}
\frac{\mathrm{d}p^h}{\mathrm{d}\xi}(\xi) = \sum_{i=0}^N \mathcal{N}^0_i(p^h) \frac{\mathrm{d}h_i}{\mathrm{d}\xi} = \sum_{i=0}^N \mathcal{N}^0_i(p^h) \left [ e_i(\xi) - e_{i+1}(\xi) \right ] = \sum_{i=1}^N \left ( \mathcal{N}^0_i(p^h) - \mathcal{N}^0_{i-1}(p^h) \right ) e_i(\xi) \;,
\label{eq:derivative_nodal_expansion}
\end{equation}
where we have used the fact that $e_0(\xi)=e_{N+1}(\xi)=0$.

Let $\mathbb{E}^{1,0}$ be the $N\times (N+1)$ matrix
\begin{equation} \label{eq:incidence_matrix_1D}
\mathbb{E}^{1,0} = \left ( \begin{array}{cccccc}
-1 & 1 & & & & \\
 & -1 & 1 & & & \\
 &  &  \ddots & \ddots & & \\
 &  &  &  -1 & 1 & \\
 &  &  &  &  -1 & 1
 \end{array} \right ) \;,
\end{equation}
then we can write \eqref{eq:derivative_nodal_expansion} as
\[ \frac{\mathrm{d}p^h}{\mathrm{d}\xi}(\xi) = \sum_{i=1}^N \sum_{j=0}^N \mathbb{E}^{1,0}_{i,j} \mathcal{N}^0_j(p^h)\ e_i(\xi) \;.\]
\end{corollary}
Taking the derivative of a nodal expansion changes the nodal degrees of freedom discussed in Example~\ref{ex:Lagrange_1D} to the integral degrees of freedom discussed in Example~\ref{ex:edge_1D}.
The matrix $\mathbb{E}^{1,0}$ is called the {\em incidence matrix}, which converts {the} nodal degrees of freedom of a function $p$ to {the} integral degrees of freedom of its derivative $\frac{\mathrm{d}p}{\mathrm{d}x}$.
\subsection{Construction of dual basis} \label{sec:construction_dual_basis}
Consider the finite element constructed in Example~\ref{ex:Lagrange_1D}. Any element $p^h \in \mathcal{P}$ can be represented as
\[ p^h(\xi) = \sum_{i=0}^N \mathcal{N}^0_i(p^h) h_i(\xi) \;,\]
where $\mathcal{N}^0_i(p^h)$ are the nodal degrees of freedom and $h_i(\xi)$ are the associated basis functions. To simplify the notation, we will write this as
\begin{equation}
p^h(\xi) = \Psi^0(\xi) \mathcal{N}^0(p^h) \;,
\label{eq:representation_nodal_function}
\end{equation}
where
\begin{equation*}
\Psi^0(\xi) = \bb{\begin{array}{ccccc}
h_0(\xi) & h_1(\xi) & \ldots & h_{N-1}(\xi) & h_N(\xi)
\end{array}}
\quad \mbox{and} \quad
\mathcal{N}^0(p^h) = \left ( \begin{array}{c}
\mathcal{N}^0_{0}(p^h) \\ [1.1ex]
\mathcal{N}^0_{1}(p^h) \\ [1.1ex]
\vdots \\ [1.1ex]
\mathcal{N}^0_{N-1}(p^h) \\ [1.1ex]
\mathcal{N}^0_N(p^h)
\end{array} \right ) \;.
\end{equation*}
Here, and in the remainder of the paper, basis functions will always be expressed as a row vector and the degrees of freedom as a column vector.
Let $p^h,q^h \in \mathcal{P}$ be both represented as in \eqref{eq:representation_nodal_function}, then the $L^2$-inner product is given by
\begin{equation*}
\left (p^h, q^h \right )_{L^2(K)} := \int _ K p^h{(\xi)}\, q^h{(\xi)}\, \mathrm{d} K = {\mathcal{N}^0(p^h)}^{\MIGnew{\intercal}} \mathbb{M}^{(0)} \mathcal{N}^0(q^h) \;,
\end{equation*}
{where} $\mathbb{M}^{(0)}$ denotes the mass matrix (or \emph{Gram} matrix) associated with the nodal basis functions
\begin{equation} \label{eq:M0}
\mathbb{M}^{(0)} = \int_K \Psi^0 \bb{\xi} ^\intercal \Psi^0 \bb{\xi} \, \mathrm{d} K \;.
\end{equation}

\begin{definition}\label{def:dual_DOF}
Let $\mathcal{N}^0(p^h)$ be the degrees of freedom for $p^h\in \mathcal{P}$, then the dual degrees of freedom, $\widetilde{\mathcal{N}}^1(p^h)$ for $p^h \in \mathcal{P}$ are defined by
\begin{equation*}
    \widetilde{\mathcal{N}}^1 \bb{p^h} := \mathbb{M}^{(0)} \mathcal{N}^0 \bb{p^h} \;.
\end{equation*}
\end{definition}
\begin{remark} \label{rem:dual}
In Definition \ref{def:dual_DOF} the superscript `$\mathrm{1}$' on $\widetilde{\mathcal{N}}^1$ corresponds to lines (geometric dimension 1) that are geometric duals of nodes (geometric dimension 0) for the one-dimensional domain $K = \bs{-1,1}$.
In \VJ{general, for} $\mathbb{R}^d$ the dual degrees of freedom of $\mathcal{N}^0$ are denoted by $\widetilde{\mathcal{N}}^d$.
\end{remark}
\begin{corollary}\label{cor:dual_nodal_basis_functions}
The dual basis functions are given by
\begin{equation} \label{eq:1d_dual_nodal_basis}
\widetilde{\Psi}^{1}(\xi) {:=} \Psi^{0}(\xi)\, \bb{{\mathbb{M}^{(0)}}}^{-1} \;.
\end{equation}
\begin{proof}
\begin{equation*}
p^h(\xi) = \Psi^0({\xi}) \mathcal{N}^{{0}}(p^h) = \Psi^0({\xi}) {\mathbb{M}^{(0)}}^{-1} \mathbb{M}^{(0)} \mathcal{N}^{{0}}(p^h) = \widetilde{\Psi}^{{1}}(\xi) \widetilde{\mathcal{N}}^{{1}}(p^h) \;.
\label{eq:primal_dual_expansion}
\end{equation*}
\end{proof}
\end{corollary}

\begin{corollary}
The mass matrix $\widetilde{\mathbb{M}}^{({1})}$ is the inverse of the mass matrix $\mathbb{M}^{(0)}$.
\begin{proof}
\[ \widetilde{\mathbb{M}}^{({1})} : = \int_K \widetilde{\Psi}^{{1}}(\xi) ^\intercal \widetilde{\Psi}^{{1}}(\xi) \,\mathrm{d}K \stackrel{\eqref{eq:1d_dual_nodal_basis}}{=} \bb{\mathbb{M}^{(0)}}^{-1} \int_K  \Psi^0({\xi}) ^\intercal \Psi^0({\xi})\, \mathrm{d}K \cdot\, \bb{\mathbb{M}^{(0)}}^{-1} \stackrel{\eqref{eq:M0}}{=} \bb{\mathbb{M}^{(0)}}^{-1} \;,\]
where $\bb{{\mathbb{M}^{(0)}}^{\intercal}}^{-1} = \bb{{\mathbb{M}^{(0)}}}^{-1}$ since $\mathbb{M}^{(0)}$ is symmetric.
\end{proof}
\end{corollary}
\begin{figure}[!htbp]
        \centering
        \begin{subfigure}[b]{0.45\textwidth}
\includegraphics[width=\textwidth]{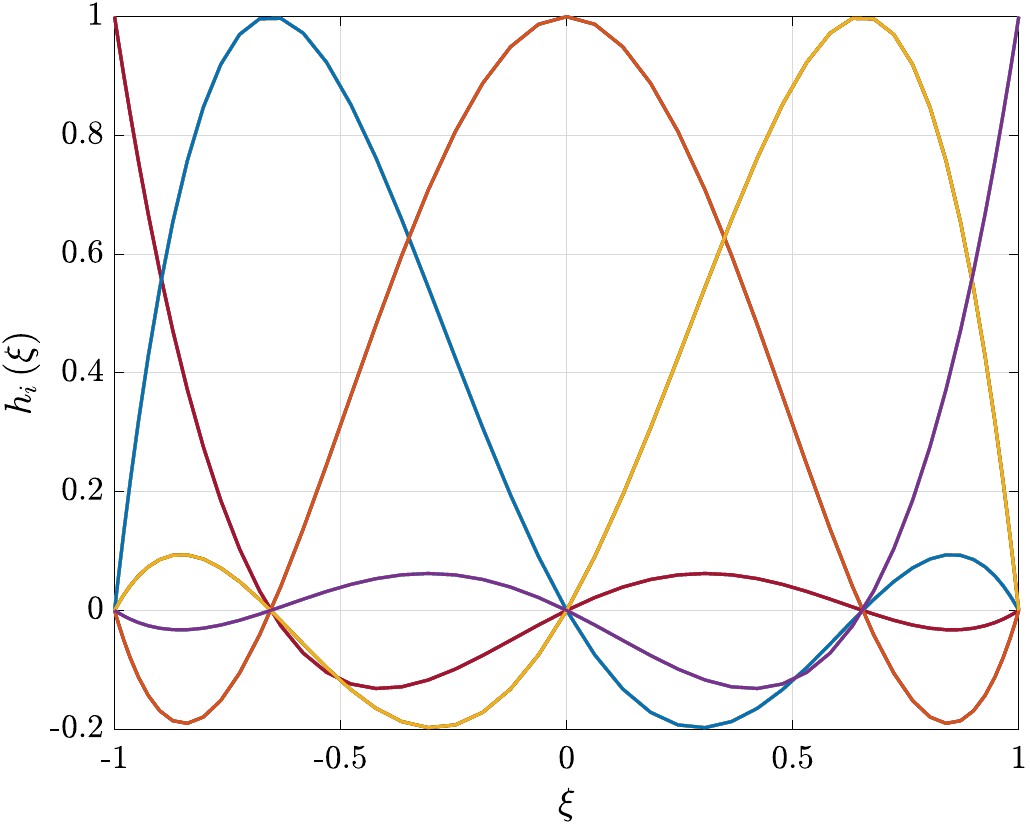}
                \caption{Lagrange polynomials}
        \end{subfigure}%
        \quad
        ~ 
        \begin{subfigure}[b]{0.45\textwidth}
                \includegraphics[width=\textwidth]{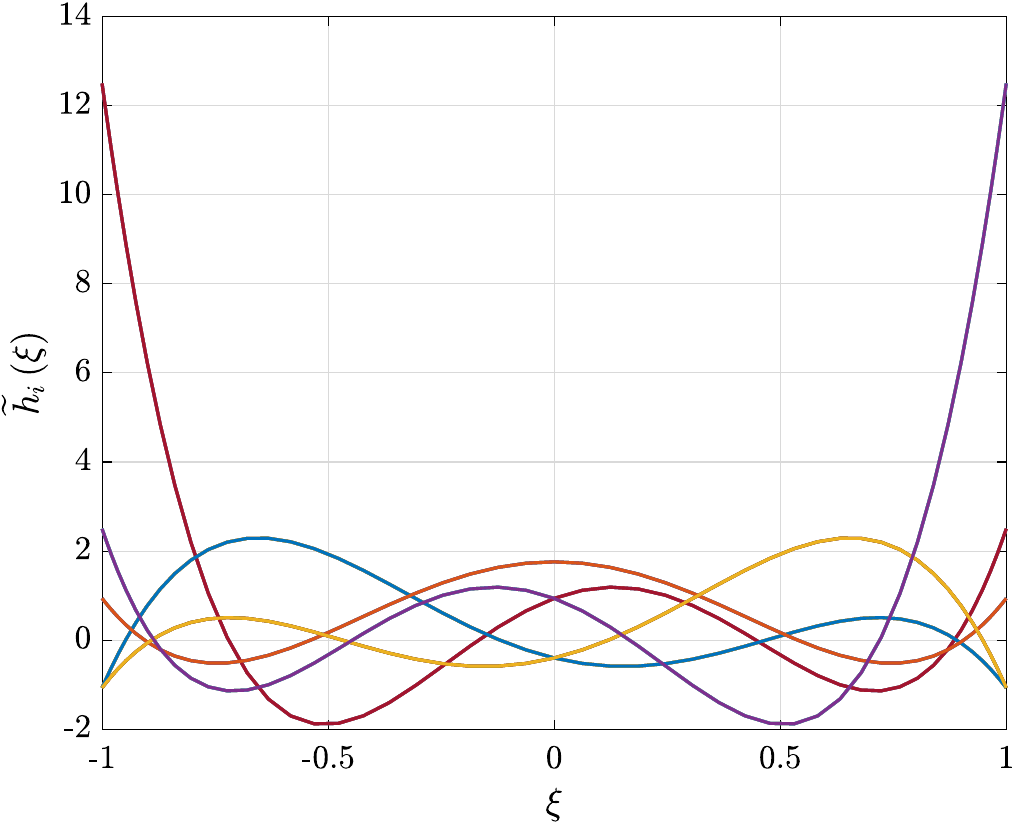}
                \caption{Dual Lagrange polynomials}
        \end{subfigure}
        \caption{The nodal Lagrange polynomial basis functions and the associated dual polynomials {for $N=4$}.} \label{fig:Lagrange_primal_dual}
\end{figure}
In Figure~\ref{fig:Lagrange_primal_dual} the Lagrange polynomials through the GLL-points and the associated dual polynomials are presented for $N=4$.

Analogous to the construction of the dual nodal polynomials, we can also construct the dual polynomials to the edge functions. Let an element $p^h \in \mathcal{Q}$ be represented as
\[ p^h(\xi) = \sum_{i=1}^N \mathcal{N}_i^1(p^h)e_i(\xi)\;.\]
In the simplified notation this can be written as
\[ p^h(\xi) = \Psi^1(\xi) \mathcal{N}^1(p^h) \;,\]
with
\begin{equation} 
\Psi^1(\xi) = \bb{\begin{array}{ccccc}
e_1(\xi) & e_2(\xi) & \ldots & e_{N-1}(\xi) & e_N(\xi)
\end{array}} \quad \mbox{and} \quad \mathcal{N}^1(p^h) = \left ( \begin{array}{c}
\mathcal{N}^1_1(p^h) \nl
\mathcal{N}^1_2(p^h) \nl
\vdots \nl
\mathcal{N}^1_{N-1}(p^h) \nl
\mathcal{N}^1_N(p^h)
\end{array} \right ) \;.
\label{eq:notation_1D_basisvector_DOFvector}
\end{equation}
We can write the $L^2$-inner product for two functions $p^h,q^h \in \mathcal{Q}$ expanded in this way as
\[ (p^h,q^h)_{L^2(K)} =\mathcal{N}^1(p^h)^\intercal \bb{\int_K  \Psi^1(\xi) ^\intercal \Psi^1(\xi) \, \mathrm{d}K} \mathcal{N}^1(q^h) = \mathcal{N}^1(p^h)^\intercal \mathbb{M}^{(1)} \mathcal{N}^1(q^h) \;,\]
with $\mathbb{M}^{(1)}$ the mass matrix (or \emph{Gram} matrix) associated with the edge polynomials
\[ \mathbb{M}^{(1)} = \int_K  \Psi^1(\xi) ^\intercal \Psi^1(\xi) \, \mathrm{d}K \;.\]

\begin{definition}\label{def:dual_edge}
Let $\mathcal{N}^1(p^h)$ be the degrees of freedom for $p^h\in {\mathcal{Q}}$, then the associated dual degrees of freedom $\widetilde{\mathcal{N}}^0(p^h)$ are defined as
\begin{equation*}
    \widetilde{\mathcal{N}}^0 \bb{p^h} := \mathbb{M}^{(1)} \mathcal{N}^1 \bb{p^h} \;.
\end{equation*}
\end{definition}
Here again we follow Remark~\ref{rem:dual} to denote the superscript `$\mathrm{0}$' on $\widetilde{\mathcal{N}}^0$ corresponding to \VJ{nodes (geometric dimension '0') that are} geometric dual of `$\mathrm{1}$' {for} {$K = \bs{-1,1}$.
In {general, for} the $d$-dimensional case the dual degrees of freedom of $\mathcal{N}^1$ will be denoted by $\widetilde{\mathcal{N}}^{d-1}$, see also Section~\ref{sec:2D_dual_spaces}.}

Following Corollary~\ref{cor:dual_nodal_basis_functions}, the dual edge functions are then given by
\begin{equation} \label{eq:1d_dual_edge_basis}
\widetilde{\Psi}^0(\xi) := \Psi^1(\xi) \bb{\mathbb{M}^{(1)}}^{-1} \;.
\end{equation}
\begin{figure}[!htbp]
        \centering
        \begin{subfigure}[b]{0.45\textwidth}
                \includegraphics[width=\textwidth]{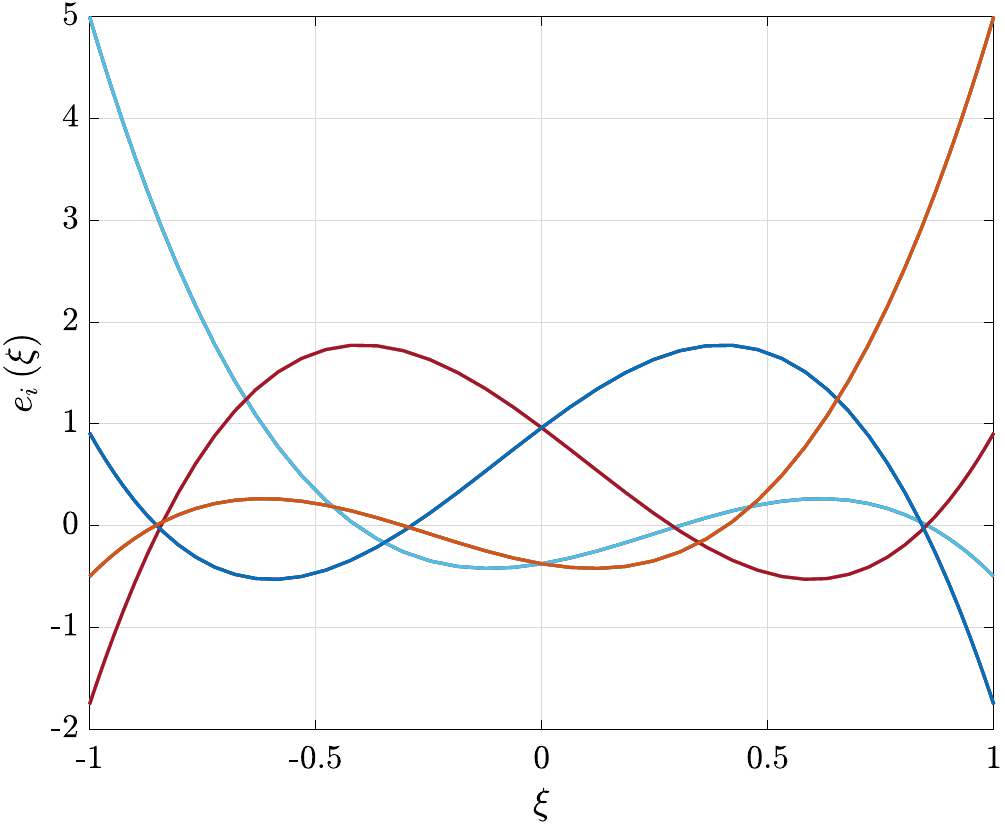}
                \caption{Edge polynomials}
        \end{subfigure}%
        \quad
        \begin{subfigure}[b]{0.45\textwidth}
                \includegraphics[width=\textwidth]{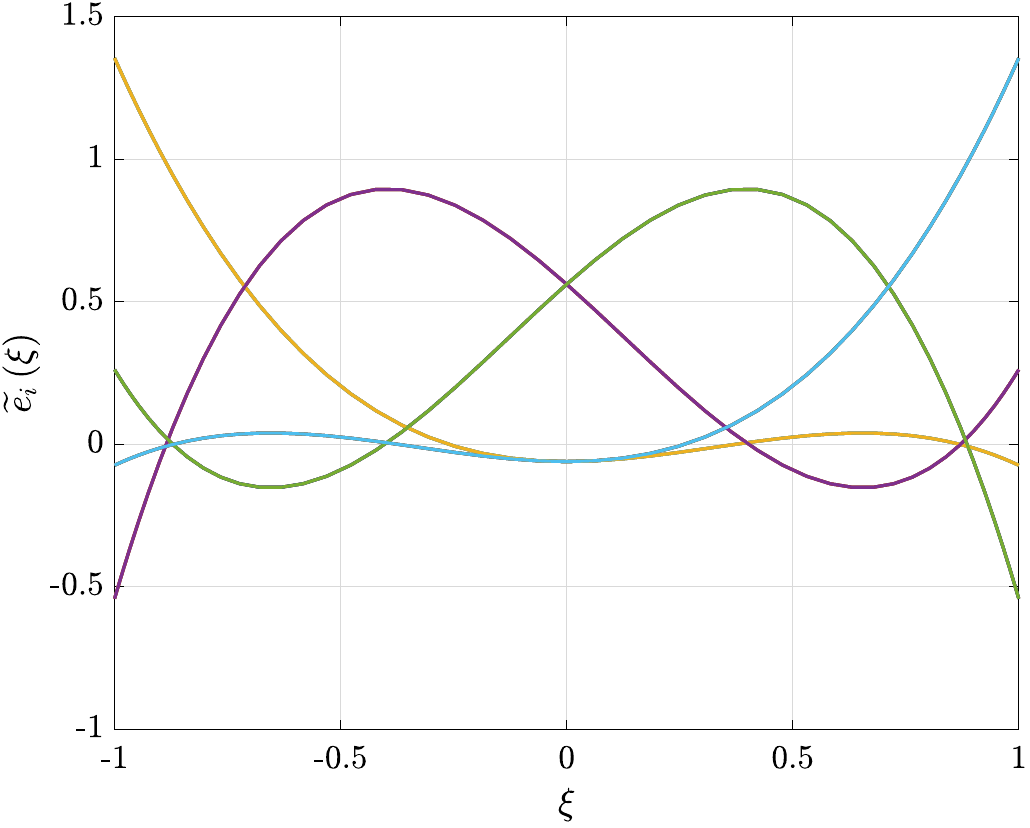}
                \caption{Dual Edge polynomials}
        \end{subfigure}
        \caption{The edge polynomial basis functions and the associated dual polynomials for {$N=4$}.} \label{fig:Edge_primal_dual}
\end{figure}
In Figure~\ref{fig:Edge_primal_dual} the edge polynomials $e_i(\xi)$ and their dual polynomials {$\widetilde{e}_i(\xi)$} are shown for {$N=4$}.
\begin{lemma}\label{lem:canonical_dual_basis}
Let $\Psi^k(\xi)$ and $\widetilde{\Psi}^{1-k}(\xi)$ for $k=0,1$ be the primal and the dual basis as defined above, then these bases are bi-orthogonal with respect to each other
\[ \int_K \widetilde{\Psi}^{1-k}(\xi)^\intercal \Psi^k(\xi)\, \mathrm{d}K = \mathbb{I}\;,\;\;\;k=0,1\;,
\]
where $\mathbb{I}$ is the $(N+1)\times(N+1)$ identity matrix for $k=0$ and the $N\times N$ identity matrix for $k=1$.
\begin{proof}
Using the definition of the dual basis, $\widetilde{\Psi}^{1-k}(\xi) = \Psi^k(\xi) \VJ{\bb{{\mathbb{M}^{(k)}}}^{-1}}$, gives
\[ \int_K \widetilde{\Psi}^{1-k}(\xi)^\intercal \Psi^k(\xi)\, \mathrm{d}K = \bb{\mathbb{M}^{(k)}}^{-1} \int_K \Psi^{k}(\xi)^\intercal \Psi^k(\xi)\, \mathrm{d}K = \bb{\mathbb{M}^{(k)}}^{-1} \mathbb{M}^{(k)} = \mathbb{I} \;,\]
{where in the second term we have used the fact that mass matrix $\mathbb{M}^{(k)}$ is symmetric.}
\end{proof}
\end{lemma}

In Remark~\ref{rem:nodal_sampling_outside_P} it was stated that nodal sampling of a function is only possible in the space $\mathcal{P}$ of continuous functions.
In a Sobolev space the elements consist of equivalence classes of functions that satisfy an integral equation and in this case nodal sampling {may not be} defined.

\begin{lemma}\label{lem:nodal_DOF_via_dual_basis}
Let $p^h \in \mathcal{P}$, then the nodal degrees of freedom are given by
\[ \mathcal{N}^0(p^h) = \int_K \widetilde{\Psi}^1(\xi)^\intercal p^h(\xi) \, \mathrm{d}K \;.\]
\begin{proof}
Every $p^h \in \mathcal{P}$ can be written as $p^h(\xi) = \Psi^{{0}}(\xi) \mathcal{N}^0(p^h)$, therefore
\[ \int_K \widetilde{\Psi}^1(\xi)^\intercal p^h(\xi) \, \mathrm{d}K = \int_K  \widetilde{\Psi}^1(\xi)^\intercal \Psi^{0}(\xi) \mathcal{N}^0(p^h) \, \mathrm{d}K = \mathcal{N}^0(p^h)\;,\]
where in the last equality we used Lemma~\ref{lem:canonical_dual_basis}.
\end{proof}
\end{lemma}
Lemma~\ref{lem:nodal_DOF_via_dual_basis} allows us to extend nodal sampling to square integrable functions.

\begin{definition}\label{def:nodal_sammpling_L2}
For $\VJ{f} \in L^2(K)$ we define the nodal degrees of freedom by
\[ \mathcal{N}^0(\VJ{f}) := \int_K \widetilde{\Psi}^1(\xi)^\intercal \VJ{f}(\xi) \, \mathrm{d}K\;.\]
Using now the fact that $\widetilde{\Psi}^{1}(\xi) = {\Psi^0(\xi)\bb{\mathbb{M}^{(0)}}^{-1}}$ this `nodal sampling' can be written as
\[ \mathcal{N}^0(\VJ{f}) = \bb{\mathbb{M}^{(0)}}^{-1} \int_K  \Psi^0(\xi)^\intercal \VJ{f}(\xi) \, \mathrm{d}K \;,\]
which is just the $L^2$-projection of $\VJ{f}$ on the space $\mathcal{P}$.
Analogously we have
\[ \mathcal{N}^1(\VJ{f}) := \int_K \widetilde{\Psi}^0(\xi)^\intercal \VJ{f}(\xi)\,\mathrm{d}K = {\bb{\mathbb{M}^{\bb{1}}}^{-1} \int _K \Psi ^{{1}} \bb{\xi}^\intercal \VJ{f}\bb{\xi} \mathrm{d}K}\;,\]
\[ \widetilde{\mathcal{N}}^0(\VJ{f}) := \int_K \Psi^1(\xi)^\intercal \VJ{f}(\xi)\,\mathrm{d} K \quad \mbox{and} \quad \widetilde{\mathcal{N}}^1(\VJ{f}) := \int_K \Psi^0(\xi)^\intercal \VJ{f}(\xi)\,\mathrm{d} K \;.\]
\VJ{These definitions are used, for e.g., in the calculation of degrees of freedom of the \varun{right  hand side} term and the boundary conditions in Sections \ref{sec:mixed_formulation} and \ref{sec:Dirichlet_Neumann}.}
\end{definition}
\VJ{
\subsection{Transformation rules for 1D function spaces}
So far we have discussed the construction of polynomial basis on a reference domain $K:= [-1,1]$.
In case of a more general interval $I:= [a,b]$ we have to introduce mappings.
For any $x \in I$, in case of a linear mapping, we have
\[ x = \Phi\bb{\xi} = \frac{a}{2}\bb{1 - \xi} + \frac{b}{2} \bb{1 + \xi} \;, \]
and its Jacobian $\boldsymbol{\mathsf{J}} = \frac{d x}{d \xi} = \frac{\bb{b-a}}{2}$.
In this case we have that the nodal basis and the edge basis transform as
\[ \Psi^{0}(x) = \varun{\overline{\Psi}}^{0}\circ\Phi^{-1}(x)  \;, \qquad \mbox{and} \qquad \Psi^{1}(x) = \frac{1}{\boldsymbol{\mathsf{J}}}\varun{\overline{\Psi}}^{1}\circ\Phi^{-1}(x)   \;, \]
\varun{where $\varun{\overline{\Psi}}^k$, $k=0,1$ are the basis functions in the reference element.}}

\VJ{The construction of mass matrices is then given by
\[ \mathbb{M}^{(0)} = \int _I \Psi ^0 \bb{x} ^\intercal \Psi ^0 \bb{x} dI =  \frac{\bb{b-a}}{2} \int _K \varun{\overline{\Psi}} ^0 \bb{\xi} ^\intercal \varun{\overline{\Psi}} ^0 \bb{\xi} dK \;, \]
and
\[ \mathbb{M}^{(1)} = \int _I \Psi ^1 \bb{x} ^\intercal \Psi ^1 \bb{x} dI = \frac{2}{\bb{b-a}} \int _K \varun{\overline{\Psi}} ^1 \bb{\xi} ^\intercal \varun{\overline{\Psi}} ^1 \bb{\xi} dK \;. \]
\subsubsection{Transformation rules for dual function spaces} \label{sec:mapping_dual_basis_1D}
The construction of dual function spaces follows exactly the same procedure as in \secref{sec:construction_dual_basis}.
\varun{Let $\Psi \bb{x}^k$, for $k=0,1$,  be the transformed basis functions, and $p^h \in \mathcal{P}$ if $k=0$, or $p^h \in \mathcal{Q}$ if $k=1$, then we have that 
\begin{equation}
\widetilde{\mathcal{N}}^{1-k}\bb{p^h} = \mathbb{M}^{(k)} \mathcal{N}^k\bb{p^h} \qquad \mbox{and} \qquad \widetilde{\Psi}^{1-k}\bb{x} = \Psi ^k\bb{x} {\mathbb{M}^{(k)}}^{-1} \qquad \qquad \mbox{for} \quad k =0,1 \;.
\end{equation}}
}
%
\subsection{Multi-element case}\label{sec:multi_elem_1D}
\begin{figure}[!htbp]
        \centering
        \begin{subfigure}[b]{0.45\textwidth}
\includegraphics[width=\textwidth]{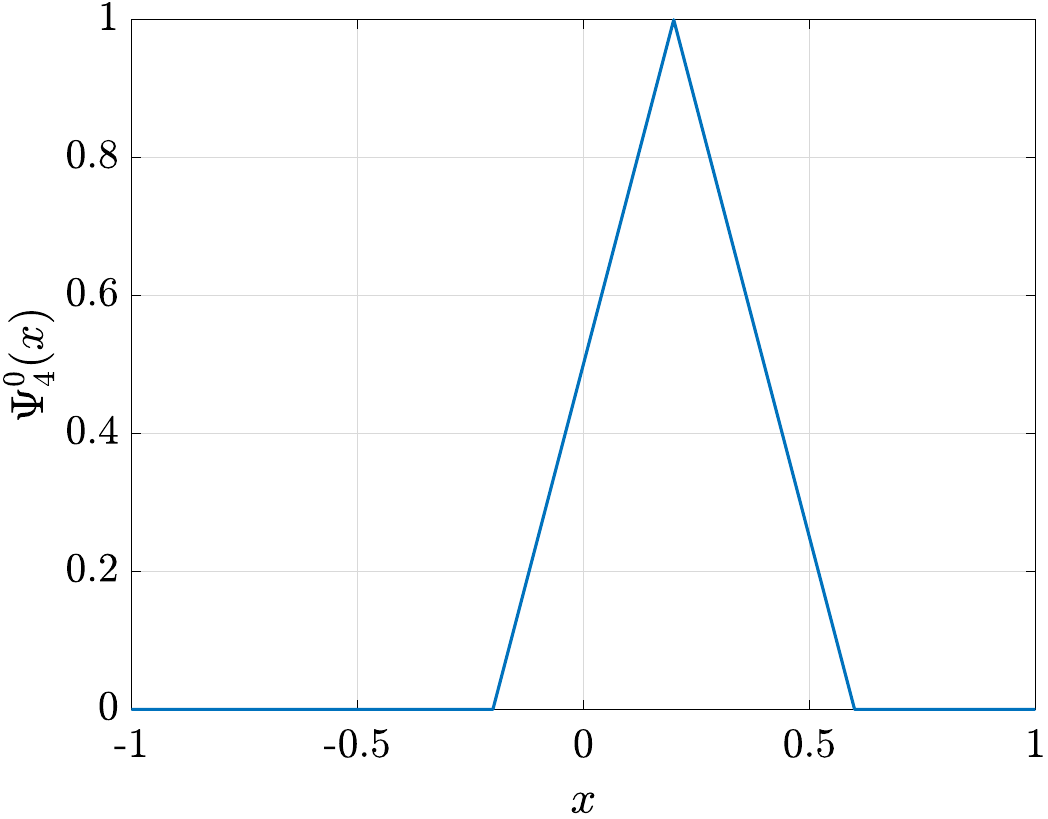}
        \end{subfigure}%
        \quad
        \begin{subfigure}[b]{0.45\textwidth}
                \includegraphics[width=\textwidth]{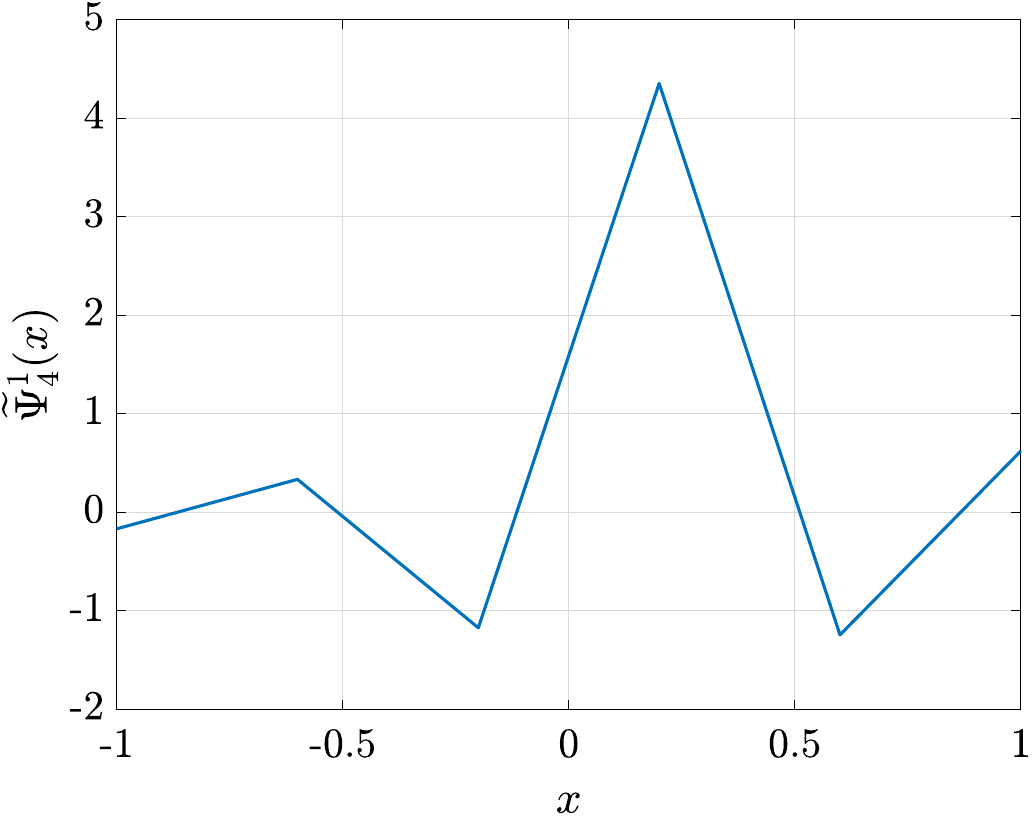}
        \end{subfigure}
        \caption{Global Lagrange basis function for $K=5$ and $N=1$ (left) and the corresponding dual basis function (right).} \label{fig:Global_Lagrange_primal_dual}
\end{figure}

Let $I=[a,b]$ with $a,b \in \mathbb{R}$ and $a<b$. Let the domain $I$ be partitioned in $K_{el}$
elements bounded by the points $a=x_0<x_1 <\ldots < x_{K_{el}-1} < x_{K_{el}}=b$ \MIGnew{and set $I_k=[x_{k-1},x_k]$.}
In the multi-element case we require that polynomials $q^h$, expanded by the Lagrange polynomials in each element, are in the space \MIGnew{$G(I) \subset H^1(I)$}, \MIGnew{with $G(I_k)=\mathcal{P}$}.
\MIGnew{The global basis functions, $\Psi_i^0(x)$ are related to the element basis functions by
\begin{equation}
\Psi_i^0(x)  =
h_j \varun{\circ \Phi ^{-1}\bb{x}}    \quad \mbox{with } i=j+(k-1)N \mbox{ and } x = \varun{\Phi\bb{x}= } \frac{1}{2}\left ( 1 - \xi \right ) x_{k-1} + \frac{1}{2}\left ( 1 + \xi \right ) x_{k} \;, -1 \leq \xi \leq 1 \;,
\label{eq:Global_nodal_basis_1D}
\end{equation}
where $k=1,\dots ,K_{el}$ and $j=0,\ldots ,N$.}
\varun{Here we use a simple linear transformation between the canonical element and the physical elements of the mesh, as presented above}.
\MIGnew{On the standard domain $[-1,1]$ we will use the variable $\xi$, while in the global, multi-element case we will denote the independent variable by $x$.}
Just as in the single element case, the dual degrees of freedom are obtained by premultiplying the global degrees of freedom with the global (assembled) mass matrix. The corresponding dual basis functions are constructed by post-multiplying the row vector of global basis functions by the inverse of the global, assembled mass matrix.
An example with \varun{$I = [-1,1]$}, $K_{\varun{el}}=5$ and $N=1$ is shown in Figure~\ref{fig:Global_Lagrange_primal_dual}.

{The global mass matrix which converts the primal degrees of freedom to the dual degrees of freedom used in the construction of Figure~\ref{fig:Global_Lagrange_primal_dual} is given by
\begin{equation*}
\mathbb{M}^{(0)} = \frac{h}{6} \left (
\begin{array}{cccccc}
2 & 1 &   &   &   &   \\
1 & 4 & 1 &   &   &   \\
  & 1 & 4 & 1 &   &   \\
  &   & 1 & 4 & 1 &   \\
  &   &   & \VJ{1} & 4 & 1  \\
  &   &   &   & 1 & 2
\end{array}
\right ) \;, \quad \quad h=\frac{2}{K} \;,
\end{equation*}
and its (approximate) inverse is given by
\begin{equation*}
\left (\mathbb{M}^{(0)} \right )^{-1} = \left (  \begin{array}{cccccc}
8.6603 &  -2.3206  &  0.6220 &  \underline{-0.1675} &   0.0478 &   -0.0239 \\
-2.3206  &  4.6411 &  -1.2440  &  \underline{0.3349} &  -0.0957 &    0.0478 \\
0.6220  & -1.2440  &  4.3541  &  \underline{-1.1722} &   0.3349 &   -0.1675  \\
-0.1675  &  0.3349 &  -1.1722  &  \underline{4.3541} &  -1.2440 &    0.6220 \\
0.0478  & -0.0957 &   0.3349 &  \underline{-1.2440}  &  4.6411 &  -2.3206 \\
-0.0239 &   0.0478 &  -0.1675  &  \underline{0.6220} &  -2.3206  &    8.6603
\end{array}
\right) \;,
\end{equation*}
where the underlined values are the nodal values plotted for the dual \VJ{basis} in Figure~\ref{fig:Global_Lagrange_primal_dual}.}

\begin{remark}
	Explicit construction of dual basis functions requires the multiplication with the inverse of the global mass matrix. For large, multi-dimensional problems this is unfeasible. 
	But the  explicit construction of dual polynomials is \VJ{in general not required for constructing the algebraic system of equations of the given problem}.
	It is sufficient to use the properties of the dual representation, such as the fact that it forms a bi-orthogonal basis to the \VJ{primal} basis, see Lemma \ref{lem:canonical_dual_basis}.
\VJ{It is in the post processing step when we need to reconstruct dual representations that we solve a linear system of equations to convert dual degrees of freedom to primal degrees of freedom and then use primal basis functions for reconstruction, whenever necessary.}
\end{remark}

We will refer to the space spanned by the dual polynomials as $\widetilde{G}(I)$, which explicitly shows that the functions are expanded in terms of dual polynomials, despite the fact that $G(I) = \widetilde{G}(I)$.

\begin{figure}[!htbp]
	\centering
	\begin{subfigure}[b]{0.45\textwidth}
		\includegraphics[width=\textwidth]{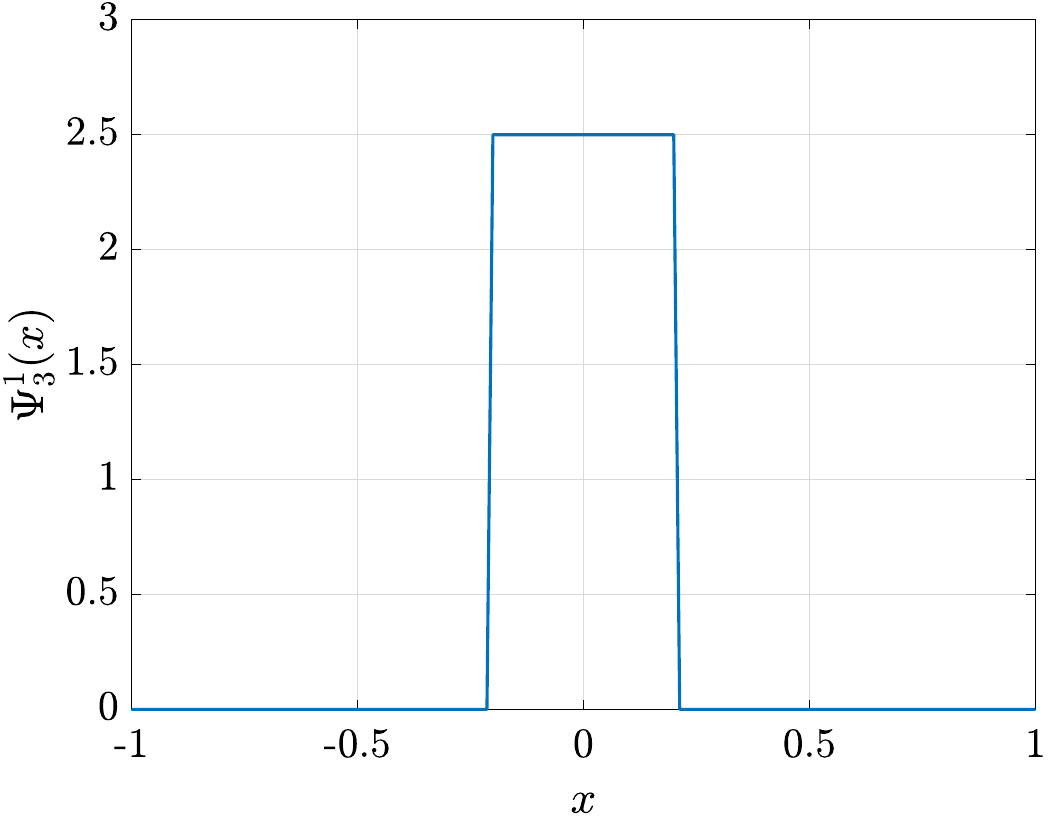}
		\caption{Global primal edge basis function}
	\end{subfigure}%
	\quad
	\begin{subfigure}[b]{0.45\textwidth}
		\includegraphics[width=\textwidth]{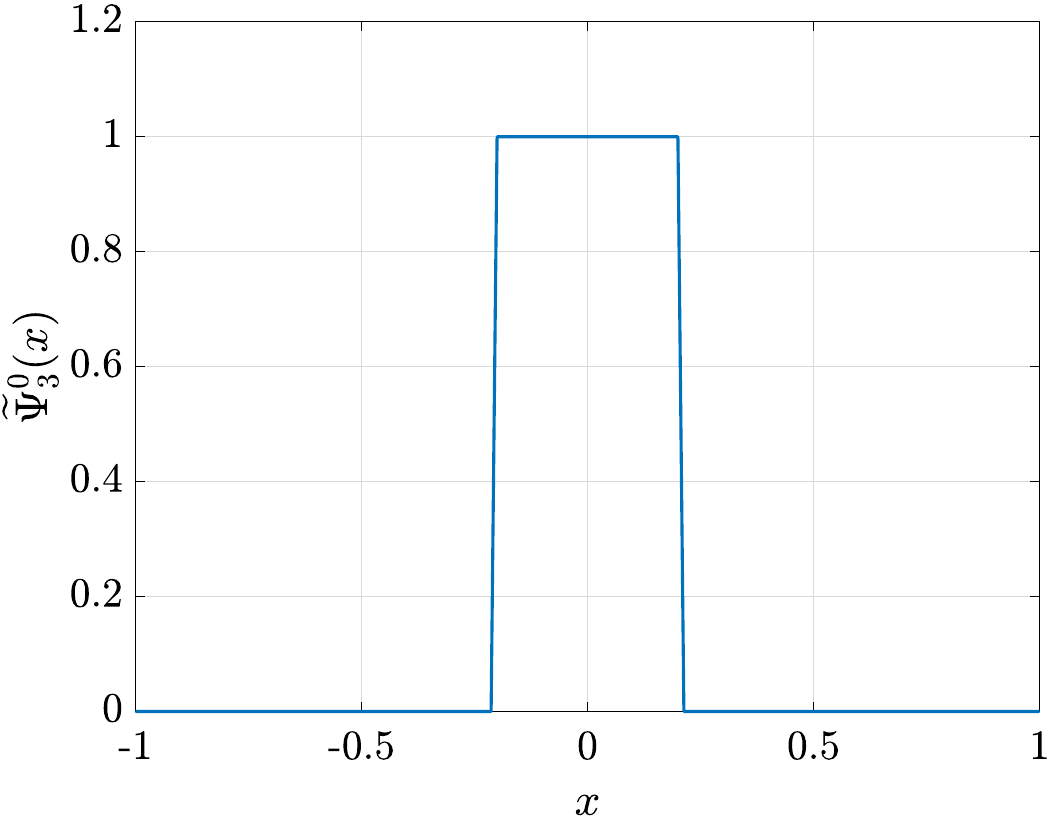}
		\caption{Global Dual of edge basis function}
	\end{subfigure}
	\caption{Global edge basis function for $K_{\varun{el}}=5$ and $N=1$ (left) and the corresponding dual function (right).} \label{fig:Global_edge_primal_dual}
\end{figure}

In the multi-element case the representation in terms of edge functions is discontinuous between elements. Its dual representation, which is a linear combination of primal edge functions, will therefore also be discontinuous between elements. The global representation in terms of edge functions, which consists of discontinuous piecewise polynomials of degree $N-1$, will therefore form elements of the subspace \VJ{$S \bb{I} \subset L^2(I)$} only, \MIGnew{with $S(I_k)=\mathcal{Q}$}.
\MIGnew{The global basis functions, $\Psi_i^1(x)$, are related to the local basis functions through
\begin{equation}
\varun{\Psi_i^1(x) = 
e_j \circ \frac{\Phi^{-1}\bb{x}}{\boldsymbol{\mathsf{J}}}} \;, \quad \quad   
\mbox{with } i=j+(k-1)N \mbox{ and } x = \varun{\Phi\bb{x}= } \frac{1}{2}\left ( 1 - \xi \right ) x_{k-1} + \frac{1}{2}\left ( 1 + \xi \right ) x_{k} \;, -1 \leq \xi \leq 1
\label{eq:global_edge_basis_1D}
\end{equation}
where $i=1,\ldots ,N\cdot K_{el}$, $k=1,\dots ,K_{el}$ and $j=1,\ldots ,N$.}
In Figure~\ref{fig:Global_edge_primal_dual} a global discontinuous representation in terms of edge functions is shown together with its global dual representation for \varun{$I = [-1,1]$,} $K_{\varun{el}}=5$ and $N-1=0$.
\subsection{Differentiation of dual variables}\label{sec:diff_dual_variables_1D}
Since the functions in $G(I)$ are continuous, we can consider the restriction of these functions to the boundary. We will denote this space as $G(\partial I)$.
For the one-dimensional case this space consists of the function's values in the points $x=a$ and $x=b$. The dual representation of a function in a point is in this case equal to the primal representation in a point.

\begin{definition}
Consider the multi-element case and let $q^h \in G\bb{I}$ be expanded in Lagrange polynomials and $\phi^h$ in dual edge polynomials in each element.
Since $\phi^h \notin H^1(I)$, we define its weak derivative\VJ{: For given $\hat{\phi} \in \widetilde{G}\bb{\partial I}$}  we  define $\widetilde{\frac{\mathrm{d}}{\mathrm{d}x}} : \widetilde{S}(I) \times \widetilde{G}(\partial I) \rightarrow \widetilde{G}(I)$ by
\begin{equation}
\int_I \left ( \widetilde{\frac{\mathrm{d}}{\mathrm{d}x}} ( \phi^h, \hat{\phi}^h )\VJ{\bb{x}} \right ) q\VJ{\bb{x}} \ \mathrm{d} x =  \int_I \phi ^h \VJ{\bb{x}} \left ( - \frac{\mathrm{d}}{\mathrm{d}x} q^h \VJ{\bb{x}} \right ) \ \mathrm{d} x + \left . \hat{\phi^h}(x)q^h(x) \right |_a^b \;, \quad \quad \forall q^h \in G(I) \;,
\label{eq:def_weak_derivative}
\end{equation}
where, $(\hat{\phi^h}(a),\hat{\phi}^h(b))^\intercal \in \widetilde{G}(\partial I)$.

\VJ{To clarify the notation we will denote the degrees of freedom on the boundary with $\mathcal{B}^0$ instead of $\mathcal{N}^0$. Therefore, we have
\[ \mathcal{B}^0\bb{q^h}= (q^h(a),q^h(b))^\intercal \in \varun{G}(\partial I) \;. \]}
\end{definition}
Since $\widetilde{\frac{\mathrm{d}}{\mathrm{d}x}} ( \phi^h, \hat{\phi}^h ) \in \widetilde{G}(I)$ we can expand it, \MIGnew{using the global basis functions dual to (\ref{eq:Global_nodal_basis_1D})} as
\MIGnew{\begin{equation*}
\widetilde{\frac{\mathrm{d}}{\mathrm{d}x}} ( \phi^h, \hat{\phi}^h ) (x) = \widetilde{\Psi}^1(x) \ \widetilde{\mathcal{N}}^{1} \left ( \widetilde{\frac{\mathrm{d}}{\mathrm{d}x}} ( \phi^h, \hat{\phi}^h ) \right )  \;,
\end{equation*}}
\MIGnew{where just as in (\ref{eq:notation_1D_basisvector_DOFvector}) $\widetilde{\Psi}^1(x)$ is a row vector of global basis functions and $\widetilde{\mathcal{N}}^{1} \left ( \widetilde{\frac{\mathrm{d}}{\mathrm{d}x}} ( \phi^h, \hat{\phi}^h ) \right )$ is a column vector of global degrees of freedom.}
\MIGnew{\begin{equation*}
q^h(x) = \Psi^0(x) \ \mathcal{N}^{0}(q^h)  \;.
\end{equation*}}
With these global representations the integral on the left hand side in (\ref{eq:def_weak_derivative}) can be written as
\begin{equation*}
\MIGnew{\int_I \left ( \widetilde{\frac{\mathrm{d}}{\mathrm{d}x}} ( \phi^h, \hat{\phi}^h ) \right ) q^h \ \mathrm{d} x = \mathcal{N}^{0}(q^h)^\intercal \widetilde{\mathcal{N}}^{1} \left ( \widetilde{\frac{\mathrm{d}}{\mathrm{d}x}} ( \phi^h, \hat{\phi}^h ) \right ) } \;,
\end{equation*}
because the two dual bases \MIGnew{$\Psi^0(x)$} and \MIGnew{$\widetilde{\Psi}^1(x)$} are bi-orthogonal by construction, see Lemma \ref{lem:canonical_dual_basis}.
Note that we can use this property without actually constructing the dual basis functions. So the integral on the left hand side of (\ref{eq:def_weak_derivative}) reduces to the vector product of the degrees of freedom of the two representations.

\MIGnew{The first integral on the right hand side of (\ref{eq:def_weak_derivative}) evaluates to
\MIGnew{\begin{equation*}
\int_I \phi^h \left ( - \frac{\mathrm{d}}{\mathrm{d}x} q^h \right ) \ \mathrm{d} x = \left ( -\mathbb{E}^{1,0} \mathcal{N}^{0}(q^h) \right )^\intercal \widetilde{\mathcal{N}}^0(\phi^h) \;,
\end{equation*}}
because the derivative of a piecewise Lagrange representation $q$ is expressed in terms of piecewise edge polynomials, Corollary~\ref{cor:derivative_incidence}, and $\phi^h$ is expanded in terms of dual edge polynomials. Both bases are bi-orthogonal therefore the explicit dependence on the basis functions cancels from this integral and the integral can be expressed in terms of degrees of freedom and \VJ{a} topological incidence \VJ{matrix}, only.}

Finally, the boundary terms reduce to
\MIGnew{\begin{equation*}
\left . \hat{\phi}^h(x) q^h(x) \right |_a^b = \mathcal{B}_1^0(q^h) \ \widetilde{\mathcal{B}}_1^0\bb{\hat{\phi}^h}  - \mathcal{B}_{0}^0(q^h) \ \widetilde{\mathcal{B}}_0^0\bb{\hat{\phi}^h}  \;.
\end{equation*}}

Since (\ref{eq:def_weak_derivative}) needs to hold for all $q^h \in G(I)$ -- and therefore for all expansion coefficients \MIGnew{$\mathcal{N}^0(q^h)$} -- this reduces to
\begin{equation}
\widetilde{\mathcal{N}}^{\VJ{1}} \left ( \widetilde{\frac{\mathrm{d}}{\mathrm{d}x}} ( \phi^h, \hat{\phi}^h ) \right ) = - \left ( \mathbb{E}^{1,0} \right )^\intercal \widetilde{\mathcal{N}}^0(\phi^h) + \mathbb{N} \MIGnew{\widetilde{\mathcal{B}}}^0\bb{\hat{\phi}^h} \;,
\label{eq:discrete_deriv_dual_variables}
\end{equation}
\varun{where the inclusion matrix $\mathbb{N}$ maps the degrees of freedom on the boundary $\mathcal{B}^0(\hat{q}^h)$ to the global degrees of freedom.}
For $K_{\varun{el}}=5$ and $N=1$, the incidence matrix is given by
\begin{equation*}
\mathbb{E}^{1,0} = \left ( \begin{array}{cccccc} 
-1  & 1  &  0 & 0 & 0 & 0 \\
0 & -1 & 1 & 0 & 0 & 0 \\
0 & 0 & -1 & 1 & 0 & 0 \\
0 & 0 & 0 & -1 & 1 & 0 \\
0 & 0 & 0 & 0 & -1 & 1
\end{array}
\right ) \;,
\end{equation*}
$\MIGnew{\widetilde{\mathcal{B}}}^0\bb{\hat{\phi}^h} = (\VJ{\hat{\phi}}(a),\VJ{\hat{\phi}}(b))^{\intercal}$ and $\mathbb{N}$ is given by
\begin{equation}
\mathbb{N} = \left ( \begin{array}{cc}
-1 & 0 \\
 0 & 0 \\
 0 & 0 \\
 0 & 0 \\
 0 & 0 \\
 0 & 1
\end{array}
\right ) \;.
\label{eq:inclusion_matrix_1D}
\end{equation}
Note that the matrices $\mathbb{E}^{1,0}$ and $\mathbb{N}$ are independent of the mesh size, $h$, and the shape of the elements and only depend on the topology of the mesh.
\MIGnew{
The flux is considered positive from left to right. On the left side of the interval, the outward unit normal points to the left and therefore the entry in $\mathbb{N}$ is $-1$, while on the right boundary the flux and the outward normal point in the same direction and therefore the entry in the $\mathbb{N}$ matrix is $+1$.}

This makes the derivative of the dual representations also a topological operation in contrast to the codifferential applied to the primal variables which depends explicitly on the metric. Note also that the matrices $\mathbb{E}^{1,0}$ and $\mathbb{N}$ do not depend on the polynomial degree $N$, so even for high order methods the derivative of the dual variables is expressed in very sparse matrices, which only contain the non-zero entries $-1$ and $1$, and the dual degrees of freedom.

\begin{figure}[!htbp]
	\centering
	\begin{subfigure}[b]{0.45\textwidth}
		\includegraphics[width=\textwidth]{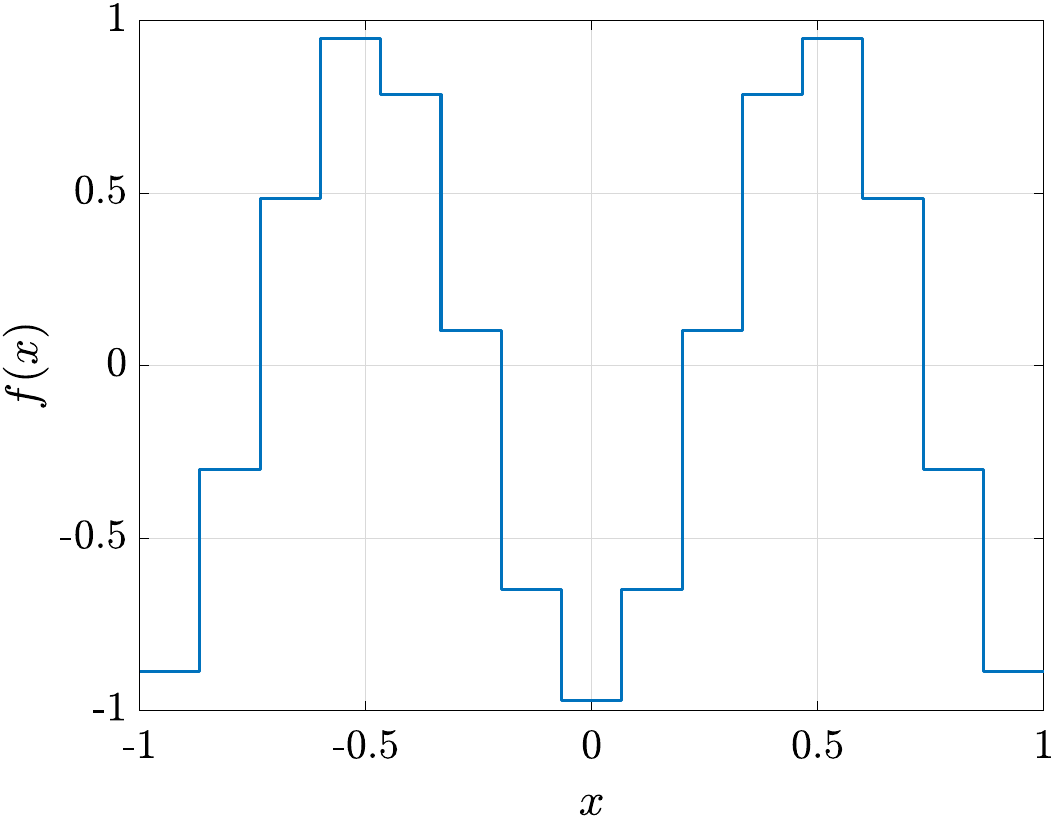}
		\caption{Approximation of $-\cos(2\pi x)$ by piecewise constant functions with $K_{\varun{el}}=15$.}
	\end{subfigure}%
	\quad
	\begin{subfigure}[b]{0.45\textwidth}
		\includegraphics[width=\textwidth]{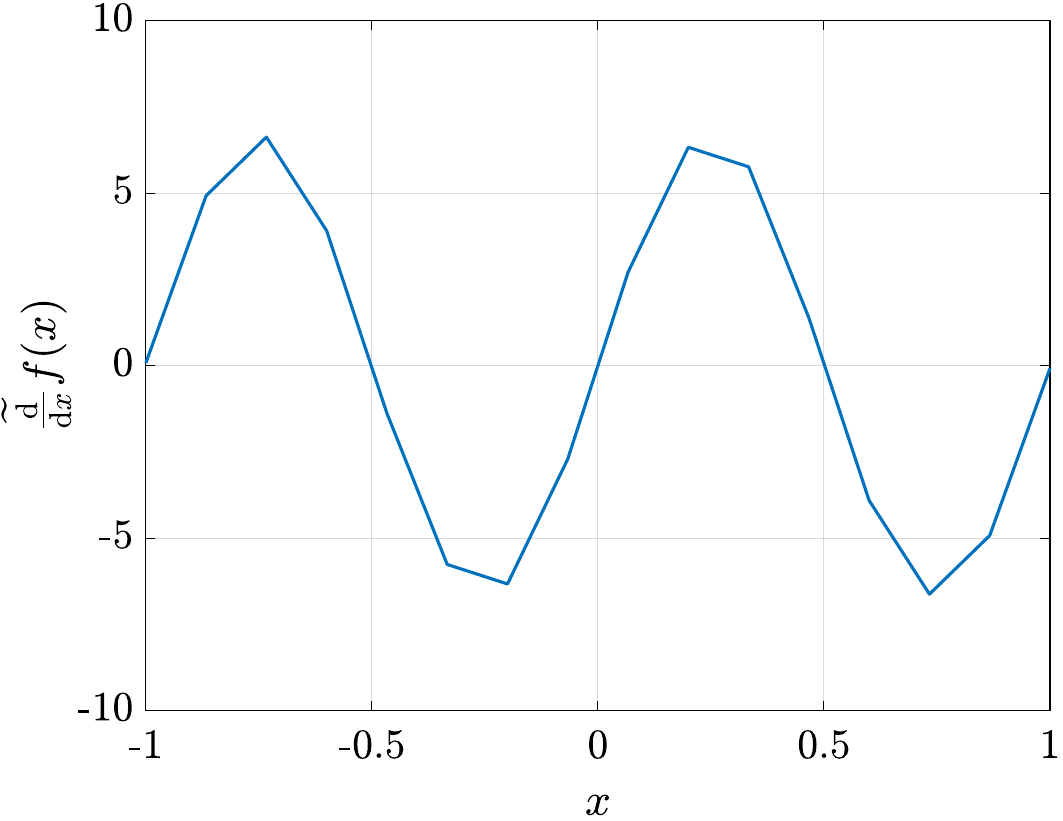}
		\caption{\VJ{Approximation of} differentation of the dual representation of $-\cos(2\pi x)$ for $K_{\varun{el}}=15$. }
	\end{subfigure}
	\caption{Application of the derivative applied to dual variables for $K_{\varun{el}}=15$ and $N=1$} \label{fig:dual_derivative_K15}
\end{figure}

\begin{figure}[!htbp]
	\centering
	\begin{subfigure}[b]{0.45\textwidth}
		\includegraphics[width=\textwidth]{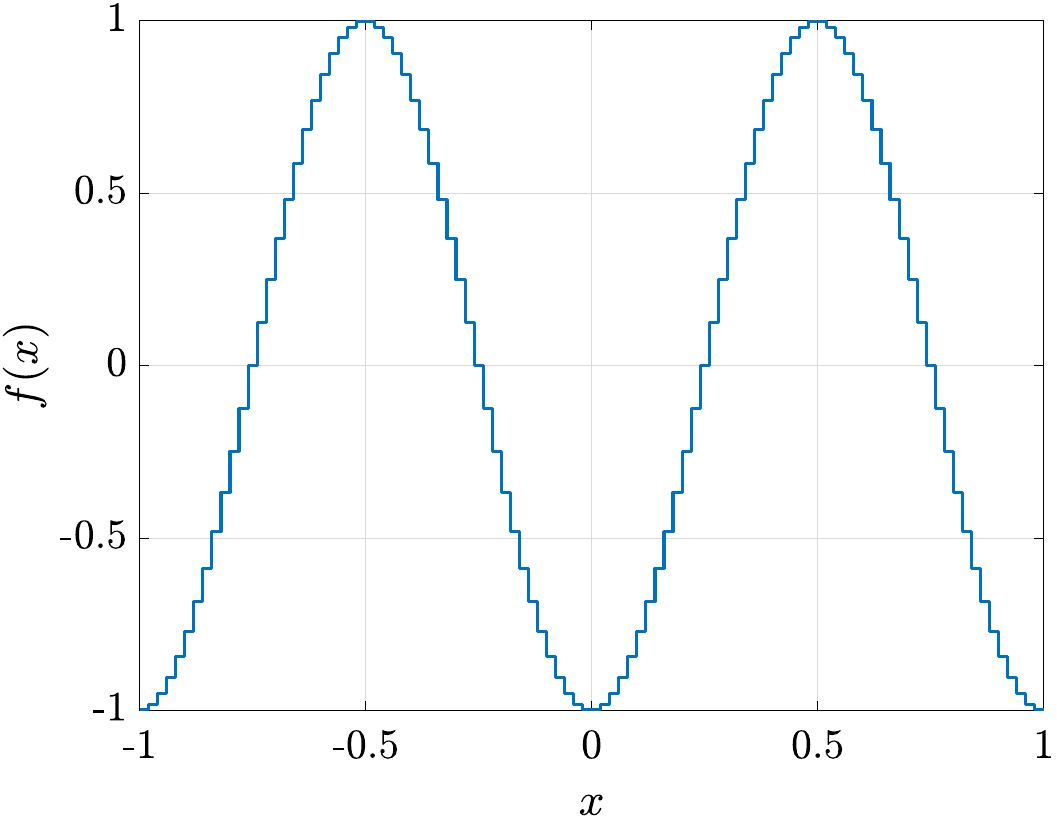}
		\caption{Approximation of $-\cos(2\pi x)$ by piecewise constant functions with $K_{\varun{el}}=100$.}
	\end{subfigure}%
	\quad
	\begin{subfigure}[b]{0.45\textwidth}
		\includegraphics[width=\textwidth]{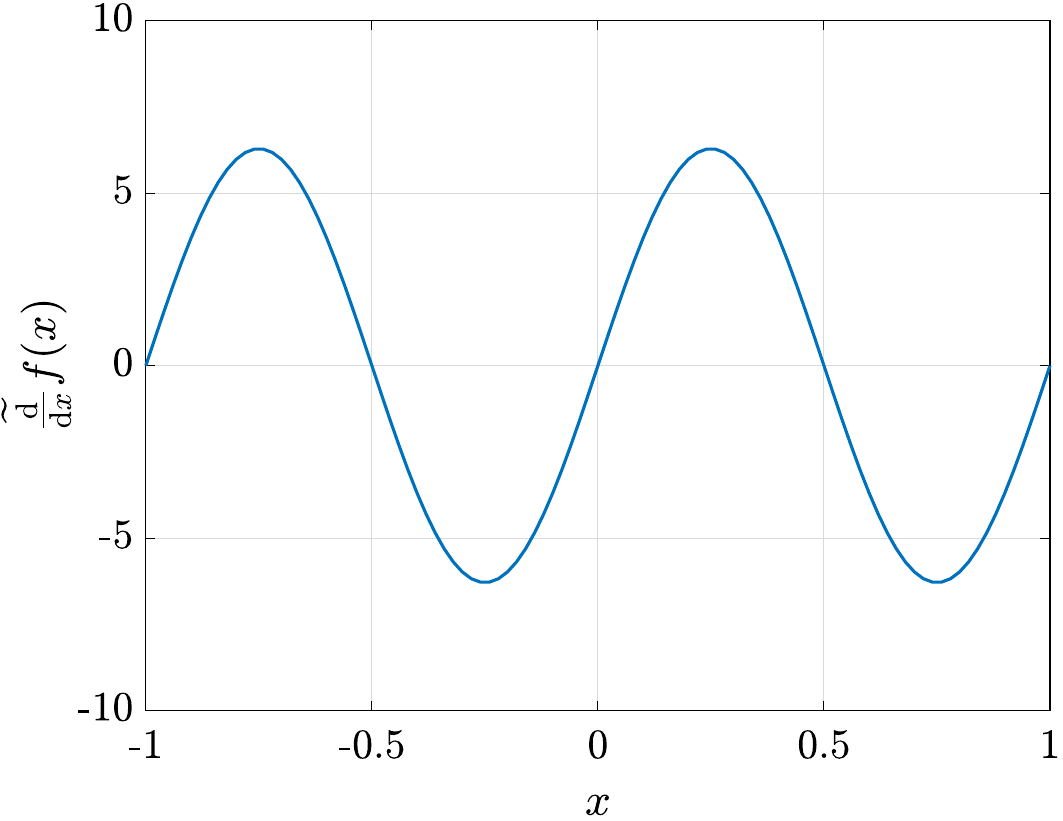}
		\caption{\VJ{Approximation of} differentation of the dual representation of $-\cos(2\pi x)$ for $K_{\varun{el}}=100$.}
	\end{subfigure}
	\caption{Application of the derivative applied to dual variables for $K_{\varun{el}}=100$ and $N=1$.} \label{fig:dual_derivative_K100}
\end{figure}

As an example of the derivative of dual \VJ{representations}, we take a piecewise constant function which approximates $- \cos(2\pi x)$ and take its derivative as described above for $K_{\varun{el}}=15$ in Figure~\ref{fig:dual_derivative_K15}, and $K_{\varun{el}}=100$ in Figure~\ref{fig:dual_derivative_K100}.
\VJ{In this case the boundary values $\hat{\phi}^h$ are set to the value of the function $- \cos(2\pi x)$ at the boundary of domain, $\varun{\widetilde{\mathcal{B}}^0\bb{\hat{\phi}}} = \bb{-1,-1}^\intercal$.}

%
%

\begin{remark}
	Although the Hodge star operator in differential geometry is well-defined, there is some freedom numerically as can be seen by the ever growing number of papers on discrete Hodge operators and Hodge matrices; primarily in the finite difference/finite volume literature. In the finite element community this is less common, but we can interpret the dual representation as a finite dimensional Hodge dual of the primal representation. In finite element formulations one generally employs inner-products for the weak formulation. An inner-product can also be written as a combination of a (topological) wedge product and a Hodge star operator using
	\[ \int_{\Omega} ( \omega,  \eta) \ \mathrm{d} \Omega = \int_{\Omega} \omega \wedge \star \eta \ \mathrm{d} \Omega \;,\]
	where $\omega$ and $\eta$ are two differential $k$-forms, see, for instance, \cite[\S4, pg. 320]{AFW10} or \cite[pg. 18]{framework2011}. 
	The dual representation is the finite element analogue of the Hodge star operator applied to the differential form $\eta$.
That the mass matrix, which converts primal to dual representations, acts as Hodge operator in a finite element setting was already observed by Tarhasaari et al., \cite{TKB99}.
\VJ{In this paper the mass matrices $\mathbb{M}^{(0)}$ and $\mathbb{M}^{(1)}$ (and later the mass matrices along the boundary, $\mathbb{M}_b^0$ and $\mathbb{M}_b^1$) play the role of Hodge operator.}
By contracting the mass matrix into a new variable, called here the dual variable, the inner-product is converted to a metric-free wedge product. 
%
\end{remark}
\subsection{Discrete de Rham sequence}
Using the function spaces defined in \secref{sec:multi_elem_1D} and \secref{sec:diff_dual_variables_1D}, we can write the primal de Rham sequence as
\begin{equation*}
\begin{tikzcd}[row sep=1.0cm, column sep = 2cm]
G\bb{\Omega} \arrow{r}{\varun{\frac{\d}{\d x}}} & S\bb{\Omega}
\end{tikzcd} \;,
\end{equation*}
and the dual de Rham sequence as
\begin{equation*}
\begin{tikzcd}[row sep=1.0cm, column sep = 2cm]
\widetilde{S}\bb{\Omega} \times \widetilde{G}\bb{\partial \Omega} \arrow{r}{\varun{\widetilde{\frac{\d}{\d x}}}} & \widetilde{G}\bb{\Omega} \times  \VJ{0}
\end{tikzcd} \;.
\end{equation*}
\section{Two-dimensional dual spaces}\label{sec:2D_dual_spaces}
In order to address more challenging problems, it is important to consider in detail the finite element spaces in 2D.
For $d = \dim \Omega = 2$, we have two sets of function spaces that obey the de Rham cohomology~\cite{AFW10,2002Hiptmair}
\begin{equation*}
\begin{tikzcd}[row sep=1.0cm, column sep = 1cm]
H\bb{\mathrm{curl};\Omega} \arrow{r}{\mathrm{curl}} & H\bb{\mathrm{div};\Omega} \arrow{r}{\mathrm{div}} & L^2 \bb{\Omega}
\end{tikzcd} \quad \mbox{and} \quad 
\begin{tikzcd}[row sep=1.0cm, column sep = 1cm]
H^1 \bb{\Omega} \arrow{r}{\mathrm{grad}} & H\bb{\mathrm{rot};\Omega} \arrow{r}{\mathrm{rot}} & L^2 \bb{\Omega}
\end{tikzcd} \;,
\end{equation*}
{where curl of a scalar field $\psi$ is the vector field $(\partial \psi/\partial y, - \partial \psi/\partial x)^{\intercal}$, the div applied to a vector field $\bm{u}=(u,v)^{\intercal}$ is the scalar field $\partial u/\partial x + \partial v/\partial y$, gradient of a scalar field $\phi$ is the vector field $\bb{\partial \phi / \partial x , \partial \phi / \partial y}^{\intercal}$ and the rot for a vector field $\bm{u}=(u,v)^{\intercal}$  is the scalar field given by $\bb{\partial v/\partial x - \partial u/\partial y}$}.
We will introduce {the} three finite element spaces, $C({\Omega}) \subset H(\mbox{curl};{\Omega})$, $D({\Omega}) \subset H(\mbox{div};{\Omega})$, $S({\Omega}) \subset L^2({\Omega})$, and the corresponding dual spaces $\widetilde{C} \bb{{\Omega}}$, $\widetilde{D} \bb{{\Omega}} $, $\widetilde{S}\bb{{\Omega}} $, respectively, such that they obey the discrete de Rham complex.

We will first present the construction of finite element spaces on \VJ{the} reference element $K = [-1,1]^2$, \VJ{and then on an arbitrary element $\Omega _k \subset \mathbb{R}^2$}.
This will be followed by construction of \VJ{finite element} spaces in case of multiple elements, and then the derivation of differential operators on the dual representations.
\subsection{The function space $C(K)\ \subset H\bb{\mathrm{curl};K}$}
\subsubsection{Primal finite element}
Let $\xi_{i},\eta_{i} \in [-1,1]$, $i=0,\dots,N$, be {the} \varun{GLL} points, and $\mathcal{P}$ denote the space of polynomials of degree $N$ defined on the interval $[-1,1]$, see Example \ref{ex:Lagrange_1D}. Consider now the polynomial tensor product space $C(K) := \mathcal{P} \otimes\mathcal{P}$. Given the set $\bm{x}$ of 2D nodes $\bm{x}_{k}$ defined as $\bm{x}:=\{\bm{x}_{i(N+1) + j}=(\xi_{i},\eta_{j}) \,\,|\,\, i,j = 0,\dots,N\}$, we can introduce for any $p^{{h}} \in C(K)$ the degrees of freedom as
			\begin{equation*}
				\mathcal{N}^{0}_{k}(p^{{h}}) := p^{{h}}(\bm{x}_{k}), \quad k = 0, \dots, (N+1)^{2}-1\;.
			\end{equation*}
		
\begin{figure}
			        \centering
			        \begin{subfigure}[b]{0.45\textwidth}
			                \includegraphics[width=\textwidth]{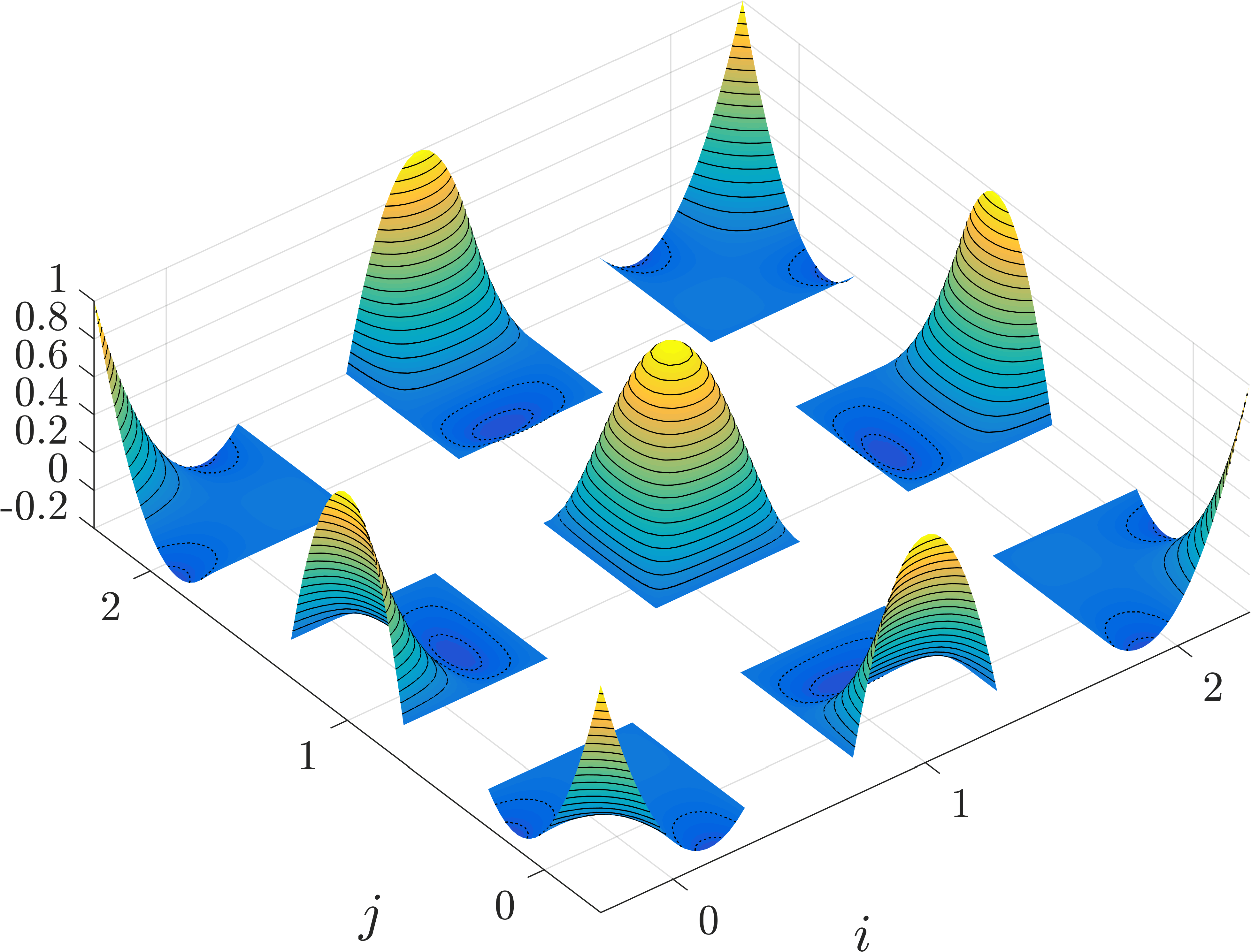}
			                \caption{The \VJ{primal} basis nodal polynomials}
			                \label{fig:H1_primal_basis}
			        \end{subfigure}%
			        \quad
			        ~ 
			        \begin{subfigure}[b]{0.45\textwidth}
			                \includegraphics[width=\textwidth]{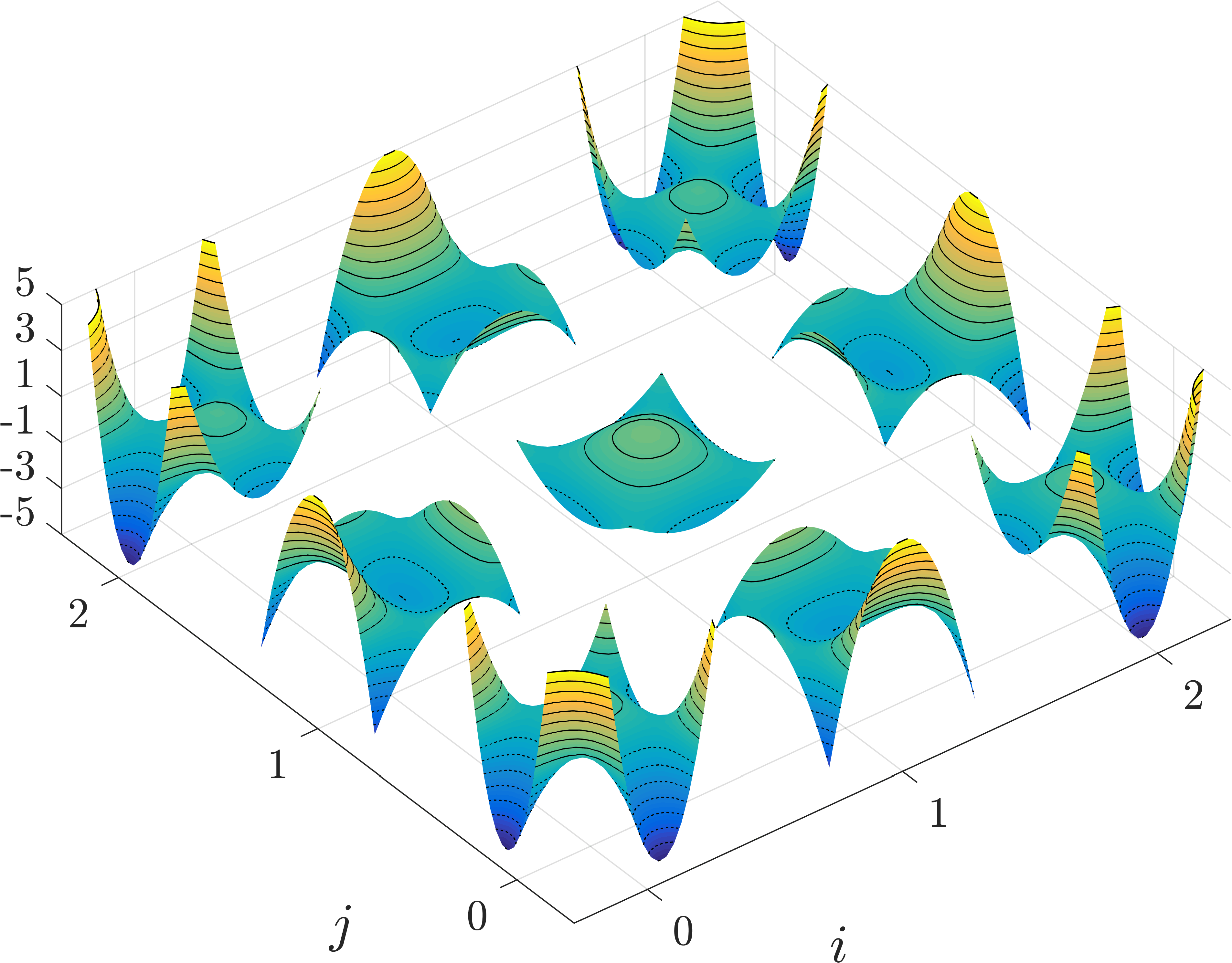}
			                \caption{The dual basis nodal polynomials}
			                \label{fig:H1_dual_basis}
			        \end{subfigure}
			        \caption{Visualization of primal, $\epsilon^{(0)}_{i(N+1)+j} (\xi,\eta)$, and dual, $\tilde{\epsilon}^{(0)}_{i(N+1)+j} (\xi,\eta)$, {nodal} basis functions of the spaces $\VJ{C}(K)$ and $\VJ{\widetilde{C}}(K)$ for polynomial degree $N=2$.}     					
			        \label{fig:H1_bases}
\end{figure}

			The basis which satisfies the Kronecker-delta property from Proposition~\ref{prop:basis_functions} is given by the  Lagrange (or nodal) polynomials, $\epsilon^{(0)}_{k}$, $k=0,\dots, (N+1)^2 -1$, through the two-dimensional GLL nodes $\boldsymbol{x}_{i(N+1) + j}=(\xi_{i},\eta_{j})$, $i,j = 0,\dots,N$, such that
			\begin{equation*}
				\epsilon^{(0)}_{i(N+1)+j} (\xi,\eta) := h_{i}(\xi) h_{j}(\eta),\quad i,j = 0, \dots, N\,,
			\end{equation*}
			where $h_{i}$ are the 1D nodal interpolants introduced in Example~\ref{ex:Lagrange_1D}. A visual representation of these basis functions for $N={2}$ is presented in \figref{fig:H1_primal_basis}.
\subsubsection{The dual finite element}
The construction of the dual basis functions follows the ideas presented in Section~\ref{sec:construction_dual_basis}. Here we outline the direct application to the 2D case of constructing the dual basis of the space $\widetilde{C} \bb{K}$.
\begin{definition}
The degrees of freedom of the dual element are given by
			\begin{equation*}
				\widetilde{\mathcal{N}}^{2}(p^{{h}}) :=  \mathbb{M}^{(0)}\mathcal{N}^{0}(p^{{h}})\;,
			\end{equation*}
where
\[ \Psi^0(\xi,\eta) = \left ( \ \epsilon_0^{(0)}(\xi,\eta) \; \dots \; \epsilon_{(N+1)^2-1}^{(0)}(\xi,\eta) \ \right ) \quad \mbox{and} \quad \mathbb{M}^{(0)} = \int_K \Psi^0(\xi,\eta)^\intercal \Psi^0(\xi,\eta) \ \mathrm{d}K \;.\]
\end{definition}
Since the dual basis functions $\widetilde{\epsilon}_{j}^{(2)}$ need to satisfy the Kronecker-delta property
			\begin{equation*}
				 \widetilde{\mathcal{N}}_i^{2}(\widetilde{\epsilon}_{j}^{(2)}) = \delta_{ij}\,,	
			\end{equation*}
by Corollary~\ref{cor:dual_nodal_basis_functions} we have that the dual basis functions can be expressed in terms of the primal basis functions as
			\begin{equation*}
				{\widetilde{\Psi}^2(\xi,\eta) :=}\left(
				\begin{array}{ccc}
					 \widetilde{\epsilon}_{0}^{(2)} & \dots & \widetilde{\epsilon}_{(N+1)^{2}-1}^{(2)}
				\end{array}
				\right)
				=
				\left(
				\begin{array}{ccc}
					\epsilon_{0}^{(0)} & \dots & \epsilon_{(N+1)^{2}-1}^{(0)}
				\end{array}
				\right)  \bb{\mathbb{M}^{(0)}}^{-1} = \Psi^0(\xi,\eta)\bb{\mathbb{M}^{(0)}}^{-1}\,.
			\end{equation*}
\VJ{A visualization of the dual basis functions in the reference domain} for $N={2}$ is presented in 	\figref{fig:H1_dual_basis}.
\subsection{The function space $D(K)\ \subset H\bb{\mathrm{div};K}$}
		\subsubsection{Primal finite element}
			Let $\xi_{i},\eta_{i} \in [-1,1]$, $i=0,\dots,N$, be the GLL points, and $\mathcal{P}$ and $\mathcal{Q}$ denote the space of polynomials of degree $N$ and $N-1$, respectively, defined on the interval $[-1,1]$, as in Example \ref{ex:Lagrange_1D} and \ref{ex:edge_1D}.
Consider now the polynomial tensor product spaces $D_{\xi} := \mathcal{P}\otimes\mathcal{Q}$ and $D_{\eta} := \mathcal{Q}\otimes\mathcal{P}$.
We introduce for any polynomial vector field $\boldsymbol{p}^{{h}}\in D\bb{K}  = D_{\xi}\times D_{\eta}$ the degrees of freedom as
\begin{equation*}
\left\lbrace
\begin{array}{lll}
\mathcal{N}^{1}_{iN + j}(\boldsymbol{p}^{{h}}) & := \int_{\eta_{j-1}}^{\eta_{j}} \boldsymbol{p}^{{h}}{\bb{\xi _i, \eta}} \cdot\boldsymbol{e}_{\xi}\,\mathrm{d}{\eta}, \quad \quad & i = 0, \dots, N \mbox{ and } j = 1, \dots, N\,,\\[1.5ex]
\mathcal{N}^{1}_{(i-1)(N+1) + j + 1 + N(N+1)}(\boldsymbol{p}^{{h}}) & := \int_{\xi_{i-1}}^{\xi_{i}} \boldsymbol{p}^{{h}} {\bb{\xi, \eta _j}} \cdot\boldsymbol{e}_{\eta}\,\mathrm{d}{\xi},\quad \quad & i = 1, \dots, N \mbox{ and } j = 0, \dots, N\,,
\end{array}
\right.
\end{equation*}
where $\boldsymbol{e}_{\xi}$, and $\boldsymbol{e}_{\eta}$ are the unit vectors in the $\xi$- and $\eta$- directions, respectively. In a polynomial vector space these integrals are well-defined.
			
\begin{figure}[!htbp]
				\centering
					\begin{subfigure}[b]{0.45\textwidth}
			             \includegraphics[width=\textwidth]{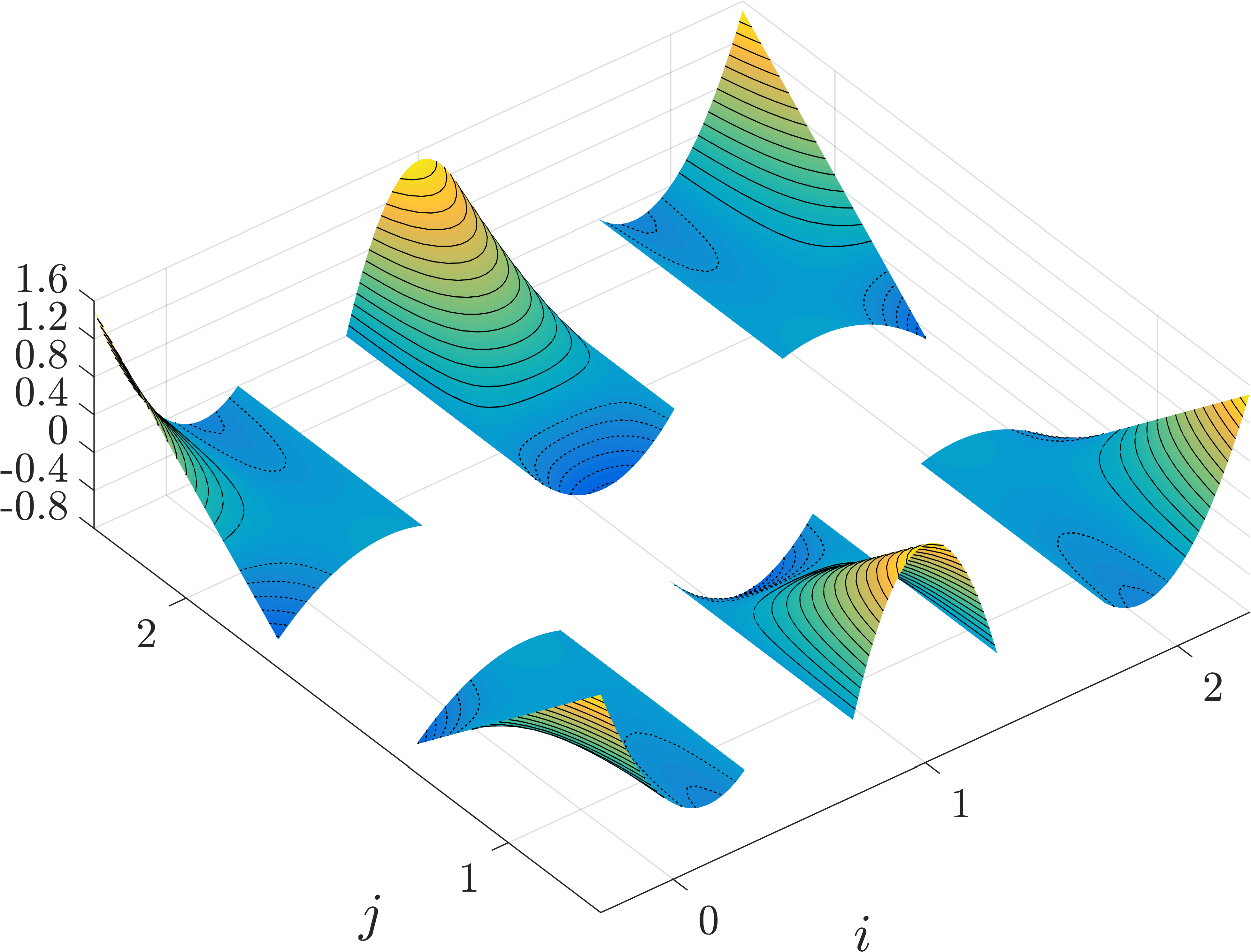}
			             \includegraphics[width=\textwidth]{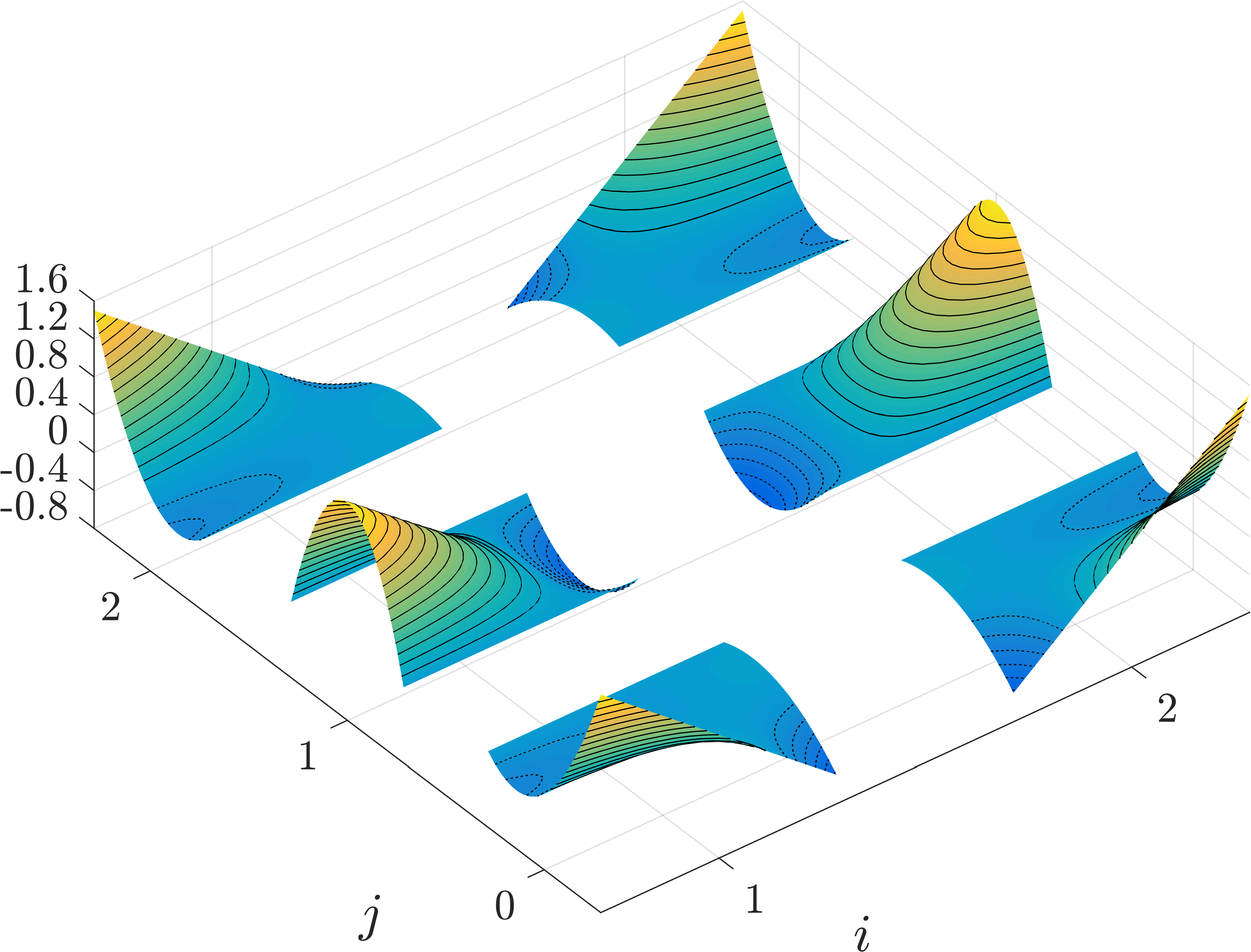}
			                		\caption{The primal basis edge polynomials}
			                		\label{fig:primal_Q1}
			        		\end{subfigure}%
					\quad
					\begin{subfigure}[b]{0.45\textwidth}
			             \includegraphics[width=\textwidth]{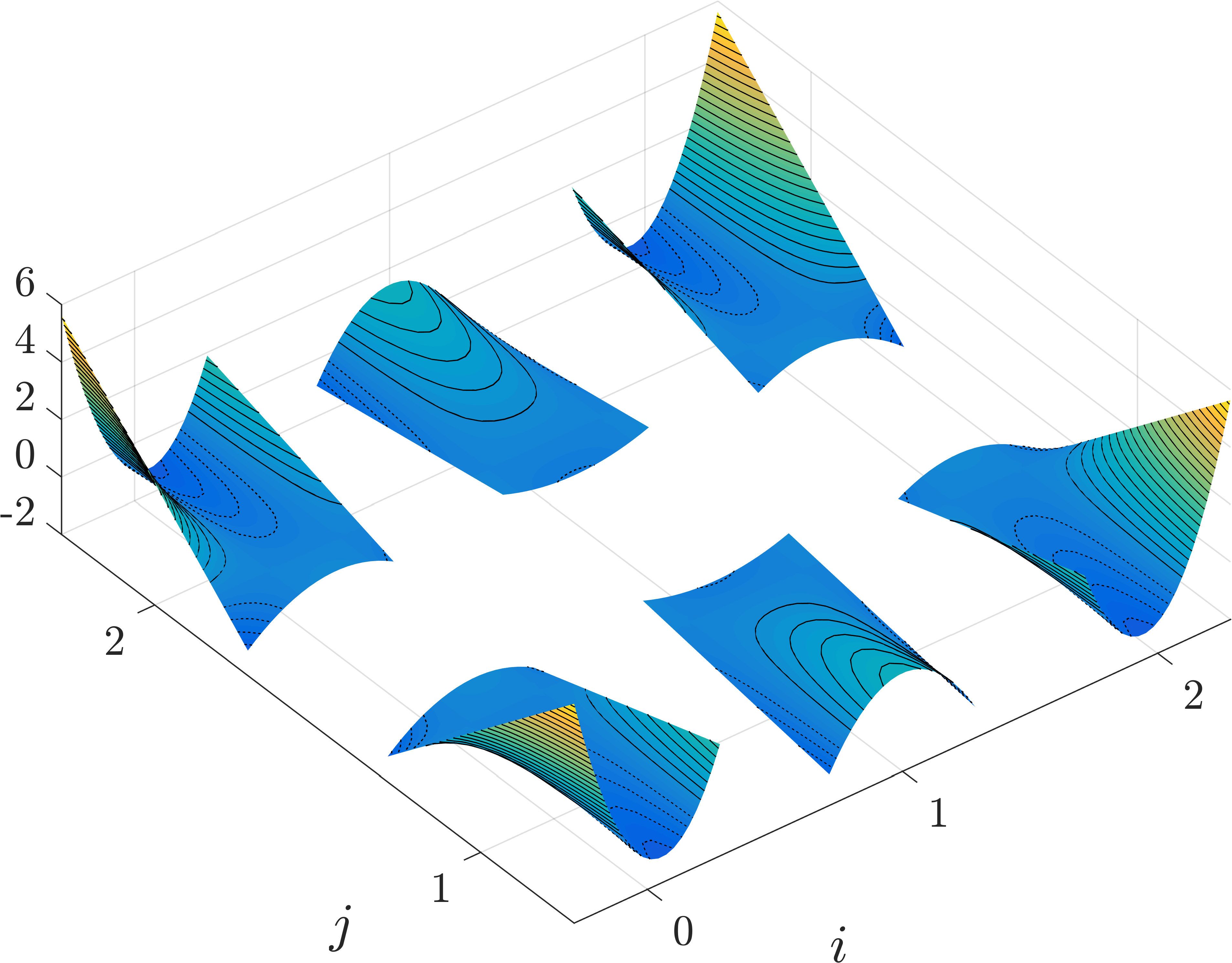}
						 \includegraphics[width=\textwidth]{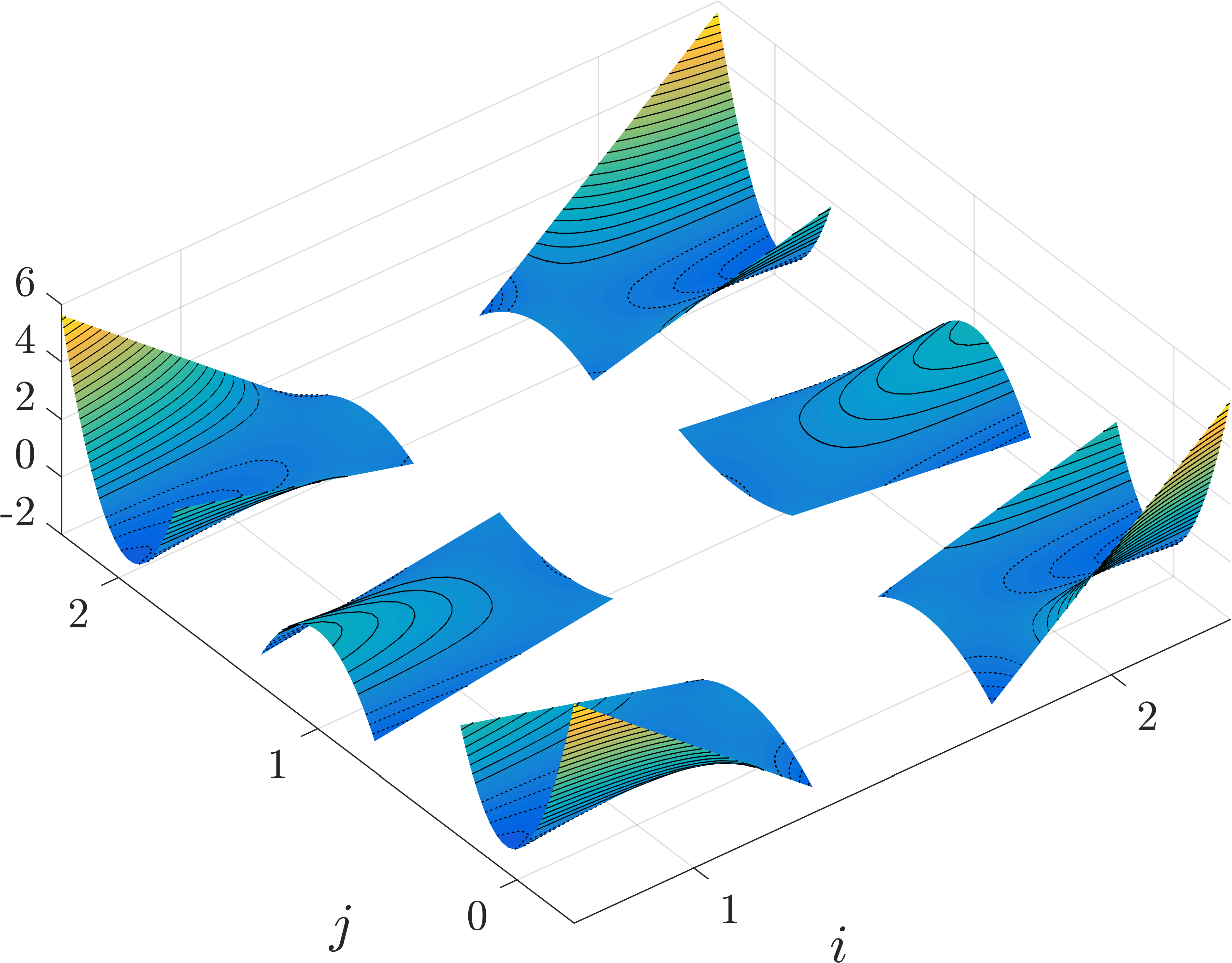}
			                		\caption{The dual basis edge polynomials}
			                		\label{fig:dual_Q1}
			        		\end{subfigure}%
			   		\caption{Visualization of one component of the primal basis functions, $\boldsymbol{\epsilon}^{(1)}_{iN+j} (\xi,\eta)$ (top left), $\boldsymbol{\epsilon}^{(1)}_{(i-1)(N+1)+j+1 + N(N+1)} (\xi,\eta)$ (bottom left), and dual basis functions, $\tilde{\boldsymbol{\epsilon}}^{(1)}_{iN+j} (\xi,\eta)$ (top right), and $\tilde{\boldsymbol{\epsilon}}^{(1)}_{(i-1)(N+1)+j+1 + N(N+1)} (\xi,\eta)$ (bottom right),  for the spaces $\VJ{D}(K)$ and $\VJ{\widetilde{D}}(K)$ with $N={2}$.}
			   		\label{fig:Hdiv_primal_dual_basis}
			 \end{figure}

It is possible to show, see \cite{Gerritsma11,Kreeft,Palha2014,BochevGerritsma}, that the basis functions which satisfy the Kronecker-delta property from Proposition~\ref{prop:basis_functions} are the 2D edge polynomials, $\boldsymbol{\epsilon}_{k}^{(1)}$, $k=1, \dots, 2N(N+1)$, defined as
			\begin{equation} \label{eq:edge_basis2D}
\Psi^1 \bb{\xi , \eta} = \left\lbrace
				\begin{array}{lll}					 \boldsymbol{\epsilon}_{iN+j}^{(1)}(\xi,\eta) & := h_{i}(\xi)\,e_{j}(\eta)\,\boldsymbol{e}_{\xi}\;,\quad \quad & i=0,\ldots,N \mbox{ and } j=1,\ldots,N  \;, \\[1.5ex]
					 \boldsymbol{\epsilon}_{(i-1)(N+1)+j+1+N(N+1)}^{(1)}(\xi,\eta) & := e_{i}(\xi)\,h_{j}(\eta)\,\boldsymbol{e}_{\eta}\;,\quad \quad & i=1,\ldots,N \mbox{ and } j=0,\ldots,N \;,
				\end{array}
				\right.
			\end{equation}
			where $h_{i}$ are the 1D nodal interpolants introduced in Example~\ref{ex:Lagrange_1D}, and $e_{j}$ are the 1D edge interpolants introduced in Example~\ref{ex:edge_1D}.
			A visual representation of these basis functions for $N={2}$ is presented in \figref{fig:primal_Q1}.
			
Let $\omega^h \in C \bb{K}$ be represented as
\[ \omega^h \bb{\xi, \eta} = \sum_{i=0}^N \sum_{j=0}^N \omega_{ij}\ h_i(\xi) \ h_j(\eta) \;,\]
then 
\[ \mbox{curl}\ \omega^h = \left ( \begin{array}{c}
\partial \omega/\partial \eta \\
-\partial \omega/\partial \xi
\end{array} \right ) \;,\]
using Corollary~\ref{cor:derivative_incidence}, we have
\begin{equation} \mbox{curl}\ \omega ^h = \left ( \begin{array}{c}
\sum_{i=0}^N \sum_{j=1}^N ( \omega_{i,j}-\omega_{i,j-1})h_i(\xi) e_j(\eta) 
\\
- \sum_{i=1}^N \sum_{j=0}^N ( \omega_{i,j}-\omega_{i-1,j})e_i(\xi) h_j(\eta) 
\end{array} \right) = \Psi^1(\xi,\eta) \mathbb{E}^{1,0} \mathcal{N}^0(\omega^h) \;,
\label{eq:curl_incidence}
\end{equation}
where $\mathbb{E}^{1,0}$ is the incidence matrix for \VJ{the edge and nodal degrees of freedom of} a 2D element\footnote{The incidence matrix $\mathbb{E}^{1,0}$ is not the same as \varun{used} in \eqref{eq:incidence_matrix_1D} which is for \VJ{degrees of freedom of} a 1D element.}.

If $\mathcal{R}(\mbox{curl};C(K))$ denotes the range of the curl operator applied to elements from $C(K)$, this implies that $\mathcal{R} \bb{\mbox{curl} ; C \bb{K}} \subset D(K)$.
This is a necessary requirement for $C(K)$ and $D(K)$ to form a finite dimensional de Rham sequence.

\subsubsection{Dual finite element}
The construction of the dual basis functions of the space $\widetilde{D}(K)$ is done in the same \VJ{way} as for the dual basis functions of the space $\widetilde{C}(K)$. In this case, the dual basis functions are expressed in terms of the primal basis functions as
\begin{align*}
\widetilde{\Psi}^1(\xi,\eta) & : =
\left [ \begin{array}{cccccc}
\widetilde{\boldsymbol{\epsilon}}_{1}^{(1)} & \ldots & \widetilde{\boldsymbol{\epsilon}}_{N(N+1)}^{(1)} & 0 & \ldots & 0 \\
0 & \ldots & 0 & \widetilde{\boldsymbol{\epsilon}}_{N(N+1)+1}^{(1)} & \ldots & \widetilde{\boldsymbol{\epsilon}}_{2N(N+1)}^{(1)}
\end{array} \right ] \nonumber \\[1ex]
 & =
\left [ \begin{array}{cccccc}
\boldsymbol{\epsilon}_{1}^{(1)} & \ldots & \boldsymbol{\epsilon}_{N(N+1)}^{(1)} & 0 & \ldots & 0 \\
0 & \ldots & 0 & \boldsymbol{\epsilon}_{N(N+1)+1}^{(1)} & \ldots & \boldsymbol{\epsilon}_{2N(N+1)}^{(1)}
\end{array} \right ] \bb{\mathbb{M}^{(1)}}^{-1} =: \Psi^1(\xi,\eta)\bb{\mathbb{M}^{(1)}}^{-1} \;,
\end{align*}
with
			\begin{equation*}
				\mathbb{M}^{(1)}_{ij} := \int_{K}\boldsymbol{\epsilon}^{(1)}_{i}(\xi,\eta)\cdot \boldsymbol{\epsilon}^{(1)}_{j}(\xi,\eta)\,\mathrm{d}K \,, \quad i,j = 1, \dots, 2N(N+1)\,.
			\end{equation*}
A visual representation of these basis functions for $N={2}$ is presented in \figref{fig:dual_Q1}.
\begin{definition}
The dual degrees of freedom for $\bm{p}^h \in \widetilde{D}\bb{K}$ are given by
\begin{equation*}
\widetilde{\mathcal{N}}^{1}(\boldsymbol{p}^{{h}}) :=  \mathbb{M}^{(1)}\mathcal{N}^{1}(\boldsymbol{p}^{{h}})\;.
\end{equation*}
\end{definition}
\subsection{The function space $S(K)\ \subset L^2\bb{K}$} \label{sec:2d_surface_basis}
\subsubsection{Primal finite element}
\begin{figure}[!htbp]
				 	\centering
					\begin{subfigure}[b]{0.45\textwidth}
			                 	 \includegraphics[width=\textwidth]{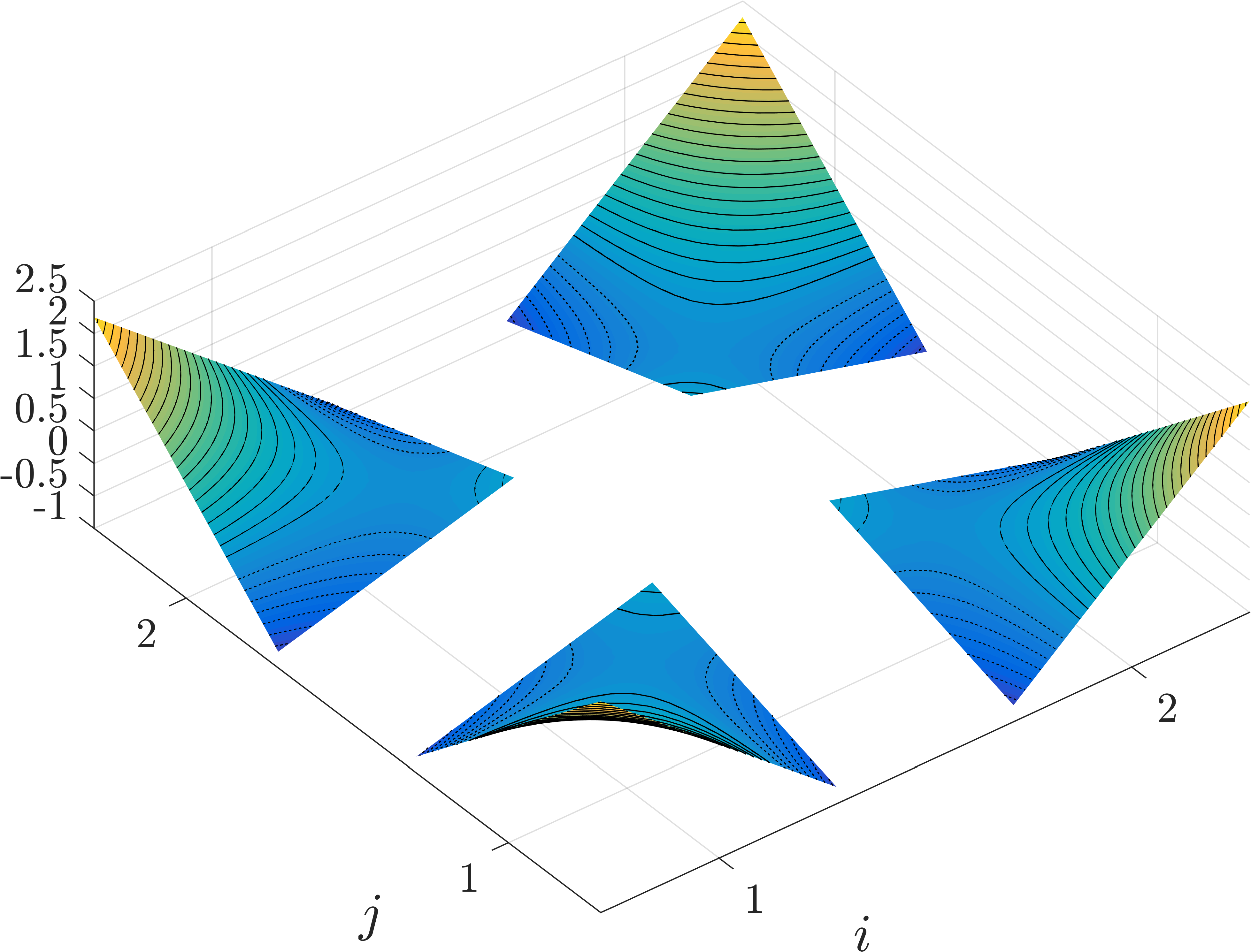}
						\caption{The primal basis surface polynomials}
			                		\label{fig:primal_Q2}
			        		\end{subfigure}%
					\quad
			   		\begin{subfigure}[b]{0.45\textwidth}
			                 	 \includegraphics[width=\textwidth]{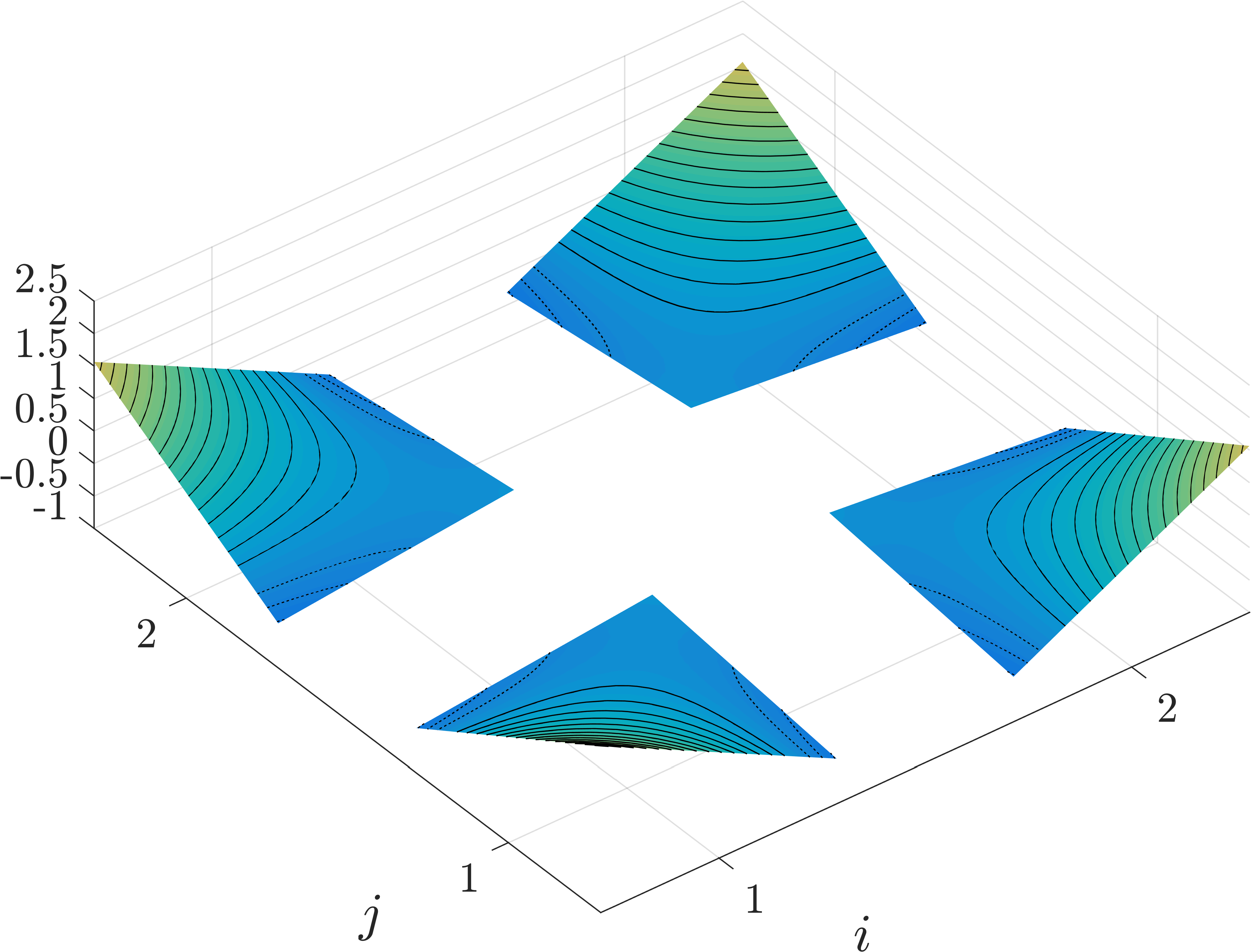}
						\caption{The dual basis surface polynomials}
			                		\label{fig:dual_Q2}
			        		\end{subfigure}%
			   		\caption{Visualization of primal basis functions {(left)}, $\epsilon^{(2)}_{{(i-1)}N+j} (\xi,\eta)$, and  dual basis functions {(right)}, $\tilde{\epsilon}^{(0)}_{{(i-1)}N+j} (\xi,\eta)$,  for the spaces $S(K)$ and $\widetilde{S}(K)$ with {$N=2$}.}
			   		\label{fig:L2_primal_dual_basis}
\end{figure}

Once again, let $\xi_{i}, \eta_{i}\in [-1,1]$, $i=0,\dots, N$, be the GLL nodes, and $\mathcal{Q}$ represent the space of polynomials of degree {$N-1$}, as in Example \ref{ex:edge_1D}.
Consider now the polynomial tensor product space ${S\bb{K}} := \mathcal{Q} \otimes \mathcal{Q}$.
The degrees of freedom for this finite element can be introduced for any polynomial  $p^h \in S\bb{K}$ as
	\begin{equation*} \label{eq:l2_primal_dof}
		\mathcal{N}^{2}(p^{{h}}) := \int_{\eta_{j-1}}^{\eta_{j}}\int_{\xi_{i-1}}^{\xi_{i}}p^{{h}} {(\xi, \eta)} \,\mathrm{d}\xi \mathrm{d}\eta\,,\quad i,j = 1, \dots, N\;.
	\end{equation*}
These integrals are well-defined in a polynomial space.
It is possible to demonstrate, see \cite{Gerritsma11,Kreeft,Palha2014,BochevGerritsma}, that the basis functions which satisfy the Kronecker-delta property from Proposition~\ref{prop:basis_functions} are the surface polynomials, $\epsilon_{k}^{(2)}$, $k=1,\dots,N^{2}$, defined as
			\begin{equation*}
				\epsilon^{(2)}_{(i-1)N + j}(\xi,\eta) := e_{i}(\xi)\,e_{j}(\eta), \quad i,j = 1,\dots,N\;,
			\end{equation*}
			where, as before, $e_{j}$ are the 1D edge interpolants introduced in Example~\ref{ex:edge_1D}.
{For ease of notation we will write these basis functions as
\begin{equation*}
\Psi ^2 \bb{\xi, \eta} := \bb{ \begin{array}{cccc}
\epsilon ^{(2)}_1\bb{\xi , \eta} & \epsilon ^{(2)}_2\bb{\xi , \eta} & \dots & \epsilon ^{(2)}_{N^2}\bb{\xi , \eta}
\end{array}
} \;.
\end{equation*}
}
A visual representation of these basis functions for {$N=2$} is presented in \figref{fig:primal_Q2}.

An element $\bm{q}^h \in D(K)$ can be represented as a linear combination of the basis \eqref{eq:edge_basis2D} as
\begin{equation} \label{eq:q_expansion_2D}
\bm{q}^h(\xi,\eta) = \left ( \begin{array}{c}
\sum_{i=0}^N \sum_{j=1}^N q^{\xi}_{i,j} h_i(\xi) e_j(\eta) 
\\[1ex]
\sum_{i=1}^N \sum_{j=0}^N q^{\eta}_{i,j} e_i(\xi) h_j(\eta) 
\end{array} \right ) \;.
\end{equation}
If we take the divergence of this vector field and use \eqref{eq:derivative_nodal_expansion} repeatedly, we have
\begin{equation}
\mbox{div}\,\bm{q}^h(\xi,\eta) = \sum_{i=1}^N \sum_{j=1}^N \left [ q^{\xi}_{i,j} - q^{\xi}_{i-1,j} + q^{\eta}_{i,j}-q^{\eta}_{i,j-1} \right ] e_i(\xi) e_j(\eta) \;.
\label{eq:div_expansion}
\end{equation}
So, we see that the divergence modifies the degrees of freedom (the expansion coefficients) and changes the basis functions from basis functions in $D(K)$ to basis functions for $S(K)$.
We can write this as
\begin{equation}
\mbox{div}\,\bm{q}^h(\xi,\eta) = \Psi^2(\xi,\eta) \mathbb{E}^{2,1} \mathcal{N}^1(\bm{q}^h) \;,
\label{eq:incidence_21}
\end{equation}	
where the incidence matrix $\mathbb{E}^{2,1}$ is a sparse matrix which only contains the non-zero entries, $-1$ and $1$, \VJ{that} can be obtained from \eqref{eq:div_expansion}.

Application of \eqref{eq:derivative_nodal_expansion} shows that $\mathcal{R}(\mbox{div};D(K)) \subseteq S(K)$, which is required for the spaces $D(K)$ and $S(K)$ to be part of the finite dimensional de Rham sequence.

Since we have that $\mbox{div}\ \mbox{curl} \ \omega^h \equiv 0$ for all $\omega^h \in C(K)$, we have, by combining (\ref{eq:curl_incidence}) and (\ref{eq:incidence_21}) that
\[ \mbox{div}\ \mbox{curl} \ \omega^h = \Psi^2(\xi,\eta) \mathbb{E}^{2,1} \mathbb{E}^{1,0} \mathcal{N}^0(\omega^h) \equiv 0 \;.\]
Since $\Psi^2(\xi,\eta)$ forms basis for $S(K)$ it implies that $\mathbb{E}^{2,1} \mathbb{E}^{1,0} \mathcal{N}^0(\omega^h) \equiv 0$. If, in addition, this needs to hold for all $\omega^h \in C(K)$, therefore we need to have $\mathbb{E}^{2,1} \mathbb{E}^{1,0} \equiv 0$, which is a well-known property of incidence matrices in mimetic methods.
		\subsubsection{Dual finite element}
			The construction of dual basis functions of the space $\VJ{\widetilde{S}}(K)$ \VJ{follows} the same steps as performed for the spaces  $\widetilde{C}(K)$ and $\widetilde{D}(K)$.
\begin{definition}
The dual degrees of freedom for $p^h \in \widetilde{S}\bb{K}$ are given by
\begin{equation*}
\widetilde{\mathcal{N}}^{0}(p^{{h}}) :=  \mathbb{M}^{(2)}\mathcal{N}^{2}(p^{{h}})\;,
\end{equation*}
with
			\begin{equation*}
				\mathbb{M}^{(2)}_{ij} := \int_{K}\epsilon^{(2)}_{i}(\xi,\eta)\,\epsilon^{(2)}_{j}(\xi,\eta)\,\mathrm{d}K\,, \quad i,j = 1, \dots, N^{2}\,.
			\end{equation*}
\end{definition}
The associated dual basis functions are expressed in terms of the primal basis functions as
			\begin{equation*}
				\widetilde{\Psi}^0(\xi,\eta) := \left(
				\begin{array}{ccc}
					\tilde{\epsilon}_{1}^{(0)} & \dots & \tilde{\epsilon}_{N^2}^{(0)}
				\end{array}
				\right)
				=
				\left(
				\begin{array}{ccc}
					\epsilon_{1}^{(2)} & \dots & \epsilon_{N^{2}}^{(2)}
				\end{array}
				\right)  \bb{\mathbb{M}^{(2)}}^{-1} =: \Psi^2(\xi,\eta) \bb{\mathbb{M}^{(2)}}^{-1} \;.
			\end{equation*}		
A visual representation of these basis functions for {$N=2$} is presented in \figref{fig:dual_Q2}.
\VJ{
\subsection{Transformation rules for 2D function spaces}
In the previous sections we have introduced \varun{the} construction of function spaces on a reference domain $K:=[-1,1]^2$.
In this section we will extend this construction to an arbitrary element $\Omega _k \subset \mathbb{R}^2$.
Let $\Phi _k$ be the diffeomorphism between the canonical domain $K$ and the arbitrary domain $\Omega _k$, and $\boldsymbol{\mathsf{J}}$ be its  Jacobian tensor such that
\[ \Phi_{k}: \bb{\xi , \eta} \in K \mapsto \bb{x,y} \in \Omega _k \qquad \mbox{and} \qquad \boldsymbol{\mathsf{J}} : = 
			\left[
				\begin{array}{cc}
					\frac{\partial\Phi_{k}^{x}}{\partial \xi} & \frac{\partial\Phi_{k}^{x}}{\partial \eta} \\
					\frac{\partial\Phi_{k}^{y}}{\partial \xi} & \frac{\partial\Phi_{k}^{y}}{\partial \eta}
				\end{array}
			\right] \;. \]
In this work, we construct the diffeomorphism using transfinite mapping~\cite{1973transfinite}. The mapping is performed element-wise such that the global map is continuous.
\subsubsection{\VJ{Transformation rules for spaces $C\bb{\Omega _k}$}}
Given a function $\bar{f} \in C(K)$ on $K$, we define $f \in C(\Omega _k)$ on $\Omega _k$ by
\begin{equation*}
f := \left(\Phi^{*}_{k}\right)^{-1} \left[\bar{f}\right] = \bar{f}\circ\Phi_{k}^{-1}\,,
\end{equation*}
where $\Phi^{*}_k$ is the pullback operator. We can then reverse the relation to obtain
\begin{equation*}
\bar{f} := \Phi^{*}_{k} \left[f\right] = f\circ\Phi_{k}\,.
\end{equation*}
Given two functions $f, g \in C(\Omega _k)$ on $\Omega _k$ their inner product can be computed by
\begin{equation*}
\bb{ f, g}_{\Omega _k} := \int_{\Omega _k} f g\, \mathrm{d}\Omega =  \int_{\Phi_{k}\left(K\right)} fg\,\mathrm{d}\Omega = \int_{K} \Phi^{*}_{k} \left[f\right] \Phi^{*}_{k} \left[g\right] \,\mathrm{det}\left(\boldsymbol{\mathsf{J}}\right)\mathrm{d}K= \int_K \bar{f} \bar{g}\,\mathrm{det}\left(\boldsymbol{\mathsf{J}}\right)\mathrm{d}K\,.
\end{equation*}
\subsubsection{Transformation rules for spaces $D\bb{\Omega _k}$}
Given a vector field $\bar{\boldsymbol{u}} \in D(K)$, the transformed vector field $\boldsymbol{u} \in D(\Omega _k)$ on $\Omega _k$ is given by
\begin{equation} \label{eq:u_with_mapping}
\boldsymbol{u} :=  \left(\Phi^{*}_{k}\right)^{-1} \left[\bar{\boldsymbol{u}}\right] = \frac{1}{\det\left(\boldsymbol{\mathsf{J}}\circ\Phi_{k}^{-1}\right)}\left(\boldsymbol{\mathsf{J}}\circ\Phi_{k}^{-1}\right)\left(\bar{\boldsymbol{u}}\circ\Phi_{k}^{-1}\right)\,.
\end{equation}
It is possible to compute the inverse of this transformation, resulting in
\begin{equation*}
\bar{\boldsymbol{u}} :=  \Phi^{*}_{k}\left[\boldsymbol{u}\right] = \det\left(\boldsymbol{\mathsf{J}}\right)\boldsymbol{\mathsf{J}}^{-1}\left(\boldsymbol{u}\circ\Phi_{k}\right)\,.
\end{equation*}
Given two vector fields $\boldsymbol{u}, \boldsymbol{v} \in D(\Omega _k)$ on $\Omega _k$ their inner product can be computed by
\begin{equation*}
\bb{ \boldsymbol{u}, \boldsymbol{v}}_{\Omega _k} := \int_{\Omega _k} \boldsymbol{u}^{\intercal} \boldsymbol{v}\, \mathrm{d}\Omega  = \int_K \left(\Phi^{*}_{k} \left[\boldsymbol{u}\right]\right)^{\intercal} \boldsymbol{\mathsf{J}}^{\intercal}\boldsymbol{\mathsf{J}}\, \Phi^{*}_{k} \left[\boldsymbol{v}\right] \frac{1}{\mathrm{det}\left(\boldsymbol{\mathsf{J}}\right)}\,\mathrm{d}K= \int_K \bar{\boldsymbol{u}}^{\intercal} \, \boldsymbol{\mathsf{J}}^{\intercal}\boldsymbol{\mathsf{J}}\,  \bar{\boldsymbol{v}}\,\frac{1}{\mathrm{det}\left(\boldsymbol{\mathsf{J}}\right)}\,\mathrm{d}K\,.
\end{equation*}
\subsubsection{Transformation rules for spaces $S\bb{\Omega _k}$}
Given a function $\bar{g} \in S(K)$ on $K$, we define $g \in S(\Omega _k)$ on $\Omega _k$ by
\begin{equation*}
g := \left(\Phi^{*}_{k}\right)^{-1} \left[\bar{g}\right] = \frac{1}{\det\left(\boldsymbol{\mathsf{J}}\circ\Phi_{k}^{-1}\right)}(\bar{g}\circ\Phi_{k}^{-1})\,.
\end{equation*}
The inverse relation can be computed in a similar fashion yielding
\begin{equation*}
\bar{g} = \Phi^{*}_{k} \left[g\right] = \det\left(\boldsymbol{\mathsf{J}}\right)(g\circ\Phi_{k}) \,.
\end{equation*}
Given two functions $f, g \in S(\Omega _k)$ on $\Omega _k$ their inner product can be computed by
\begin{equation*}
\bb{ f, g}_{\Omega _k} := \int_{\Omega _k} f g\, \mathrm{d}\Omega =  \int_{\Phi_{k}\left(K\right)} fg\,\mathrm{d}\Omega = \int_K \Phi^{*}_{k} \left[f\right] \Phi^{*}_{k} \left[g\right] \,\frac{1}{\mathrm{det}\left(\boldsymbol{\mathsf{J}}\right)}\mathrm{d}K= \int_K\bar{f} \bar{g}\,\frac{1}{\mathrm{det}\left(\boldsymbol{\mathsf{J}}\right)}\mathrm{d}K\,.
\end{equation*}
\subsubsection{Transformation rules for dual function spaces}
The construction of dual spaces for an arbitrary element $\Omega _k \subset \mathbb{R}^2$ follows the same procedure as that for the canonical domain $K$ .
}

\varun{Let $\Psi^k \bb{\bm{x}}$ for $k=0,1,2$ be the transformed basis functions in space $C\bb{\Omega}$, $D\bb{\Omega}$ and $S\bb{\Omega}$, respectively, and let the associated mass matrix be given by
\[ \mathbb{M}^{(k)} = \int _\Omega {\Psi ^k \bb{\bm{x}}}^\intercal \Psi ^k \bb{\bm{x}} \mathrm{d}\Omega \;. \] 
Then we have the transformed dual basis functions and the dual degrees of freedom given by
\[ \widetilde{\Psi}^{d-k}\bb{\bm{x}} = \Psi^k \bb{\bm{x}} \bb{\mathbb{M}^{(k)}}^{(-1)} \qquad \mbox{and} \qquad \widetilde{\mathcal{N}}^{d-k}\bb{p^h} = \mathbb{M}^{(k)} \mathcal{N}^k\bb{p^h} \quad \mbox{for} \quad d=2 \;, k =0,1,2 \;. \]}
\subsection{Multi-element case}\label{sec:Multi-elem_2D}
Let a two dimensional, bounded domain $\Omega$ with Lipschitz continuous boundary be partitioned in $K_{el}$ non-overlapping conforming elements. 
In the multi-element case, the global dual representation is constructed similarly to the one dimensional case described in Section~\ref{sec:multi_elem_1D}.
{We now consider the global degrees of freedom and the dual degrees of freedom are obtained by pre-multiplying the global degrees of freedom with the global \VJ{mass} matrix.}
The corresponding global dual basis functions are constructed by post-multiplying the global primal basis functions by the inverse of the global mass matrices.
Once again, these operations are far too expensive to perform explicitly, but in practice we do not need to construct these dual degrees of freedom and basis functions explicitly and merely use their properties.

Global representations of piecewise polynomials in $C(\Omega) \subset H(\mbox{curl};\Omega)$ are in $C^0(\Omega)$; the nodal degrees of freedom along an element boundary are shared by neighbouring elements. Since the dual representation is a linear combination of the primal basis functions $\widetilde{C}(\Omega)$ also contains continuous functions.

For the global piecewise polynomials in $D(\Omega)\subset H(\mbox{div};\Omega)$, the degrees of freedom at the boundary, which represent integrated fluxes, are shared by adjacent elements.
The normal component of these vector fields is continuous between adjacent elements, while the tangential component is discontinuous.
The dual representation has different degrees of freedom and different basis functions, but being a linear combination of basis functions in $D(\Omega)$, these functions also have continuous normal component and discontinuous tangential components between neighbouring elements.

Finally, the piecewise polynomials in $S(\Omega) \subset L^2(\Omega)$ are discontinuous between elements and therefore  so are the elements of $\widetilde{S}(\Omega)$.

\subsection{Vector operations on dual variables} \label{sec:dula_operators_2d}
Just as in the one dimensional case we will define weak derivatives for which we use the dual presentation. 

\subsubsection{The gradient acting on $\widetilde{S}(\Omega) \times \widetilde{D}(\partial \Omega)$}\label{def:weak_gradient}

\begin{definition} \label{def:weak_grad_2d}
We define the \varun{operator} $\VJ{\widetilde{\mathrm{grad}}} : \widetilde{S}(\Omega) \times \widetilde{D}(\partial \Omega)\rightarrow \widetilde{D}(\Omega)$ such that, for $(s^h,\hat{s}^h) \in \widetilde{S}(\Omega)\times \widetilde{D}(\partial \Omega)$
\begin{equation}
\int_{\Omega} \VJ{\widetilde{\mathrm{grad}}}\ (s^h,\hat{s}^h) \cdot \VJ{\bm{q}}^h \ \mathrm{d}\Omega = \int_{\Omega} s^h \left ( -\VJ{\mathrm{div}}\ \VJ{\bm{q}}^h \right ) \ \mathrm{d} \Omega + \int_{\partial \Omega} \hat{s}^h \left ( \VJ{\bm{q}}^h \cdot \VJ{\bm{n}} \right ) \ \mathrm{d} \Gamma \;, \quad \quad \forall \VJ{\bm{q}}^h \in D(\Omega)\;.
\label{eq:weak_grad_2D}
\end{equation}
\end{definition}
Note that the standard gradient cannot be applied to elements of $\widetilde{S}(\Omega)$ because these elements lack the required smoothness. We see that in (\ref{eq:weak_grad_2D}) the integral on the left hand side is a bilinear form on $\widetilde{D}(\Omega)\times D(\Omega)$, the first integral on the right hand side is a bilinear form on $\widetilde{S}(\Omega)\times S(\Omega)$ and the boundary integral on the right is a bilinear form on $\widetilde{D}(\partial \Omega)\times D(\partial \Omega)$. 

\MIGnew{Just as in the one-dimensional case we use standard finite element assembly of local basis functions and degrees of freedom to form global basis functions and degrees of freedom, which we will denote by ${\Psi}^1(x,y)$ and ${\mathcal{\bm{N}}}^1$, respectively. An element $\bm{q}^h \in D(\Omega)$ can be expanded as
	\begin{equation}
	\bm{q}^h(\bm{x}) =   \Psi^1(\bm{x}) \mathcal{N}^1(\bm{q}^h) \;.
	\label{eq:global_expansion_flux}
	\end{equation}
Since $\widetilde{\mbox{grad}}$ maps into $\widetilde{D}(\Omega)$ we can expand $\widetilde{\mbox{grad}}\ (s^h,\hat{s}^h)$ as
\begin{equation}
\widetilde{\mbox{grad}}(s^h,\hat{s}^h) = \widetilde{\Psi}^1(\bm{x}) \widetilde{\mathcal{N}}^1(\widetilde{\mbox{grad}}(s^h,\hat{s}^h)) \;.
\end{equation}
}
This shows that the integral on the left hand side of \eqref{eq:weak_grad_2D} evaluates to
\MIGnew{\[ 
	\int_{\Omega} \widetilde{\mbox{grad}}\ (s^h,\hat{s}^h) \cdot \bm{q}^h \ \mathrm{d}\Omega =
	{\mathcal{N}^1(\bm{q}^h)}^{\intercal}    \widetilde{\mathcal{N}}^1(\widetilde{\mbox{grad}}\ (s^h,\hat{s}^h))  \;.
	\]}
\MIGnew{If the expansion of $\bm{q}^h$ is given by (\ref{eq:global_expansion_flux}), div$\ \bm{q}^h$ is expanded as
\begin{equation}
\mbox{div}\ \bm{q}^h(\bm{x}) = \Psi^2(\bm{x}) \mathbb{E}^{2,1} \mathcal{N}^1(\bm{q}^h) \;,
\end{equation}
and $s^h \in \widetilde{S}(\Omega)$ as
\begin{equation}
s^h(\mathbf{x}) =  \widetilde{\Psi}^0(\mathbf{x}) \widetilde{\mathcal{N}}^0(s^h) \;,
\end{equation}}
\MIGnew{The first integral on the right hand side \varun{of \eqref{eq:weak_grad_2D}} gives
	\MIGnew{\[ 
		\int_{\Omega} s^h \left ( -\mbox{div}\ \bm{q}^h \right ) \ \mathrm{d} \Omega = - {\mathcal{N}^1(\bm{q}^h)}^{\intercal} {\mathbb{E}^{2,1}}^{\intercal} \widetilde{\mathcal{N}}^0(s^h)   \;, 
		\]}
	where once again, due to the bi-orthogonality of the primal and dual bases, the integral can be expressed in terms of the degrees of freedom only (no integration points or integration weights are required). }

\MIGnew{The restriction of $\bm{q}^h$ to the boundary $\partial \Omega$, is given by
\begin{equation}
\bb{\bm{q}^h\cdot \bm{n}}(\bm{x}) = \Psi^1(\bm{x}) \mathbb{N}_1 \mathbb{N}_1^{\intercal} \mathcal{N}^1(\bm{q}^h) \;,\quad \quad \mbox{for } \bm{x} \in \partial \Omega \;.
\label{eq:expansion_trace_q}
\end{equation}
\varun{In this case $\mathbb{N}_1$ is the 2D inclusion matrix} which maps degrees of freedom on the boundary, $D(\partial \Omega)$, to the global degrees of freedom in $D(\Omega)$, where the subscript $'1'$ on the matrix denotes that degrees of freedom are associated to edges (geometric dimension $'1'$).
Just as its 1D analogue, see for instance (\ref{eq:inclusion_matrix_1D}), this sparse matrix only has $-1$ and $1$ as non-zero entries.
\varun{See Example~\ref{ex:2D_inclusion_matrix_N1} on how to construct $\mathbb{N}_1$.}
The basis functions $\Psi^1_b(\bm{x}) = \Psi^1(\bm{x}) \mathbb{N}_1$ form a basis for the trace space and the degrees of freedom $\mathbb{N}_1^{\intercal} \mathcal{N}^1(\bm{q}^h)$ the associated degrees of freedom. The mass matrix of this trace space is given by
\begin{equation}
\mathbb{M}_b^{(1)} = \int_{\partial \Omega} {\Psi_b^1(\bm{x})}^{\intercal} \Psi_b^1(\bm{x}) \ \mathrm{d}\Gamma \;. 
\end{equation}
Using this mass matrix, we can set up the dual basis functions in the trace space given by $\widetilde{\Psi}_b^0(\bm{x})=\Psi_b^1(\bm{x}) {\mathbb{M}_b^{(1)}}^{-1}$. With these dual basis functions, we can express $\hat{s}^h \in \widetilde{D}(\partial \Omega)$ as
\begin{equation}
\hat{s}^h(\bm{x}) = \widetilde{\Psi}^0_b(\bm{x}) \widetilde{\mathcal{B}}^0(\hat{s}^h) \;, \quad \quad \mbox{for } \bm{x} \in \partial \Omega \;,
\label{eq:expansion_trace_hat_s}
\end{equation}
where, using Definition \ref{def:nodal_sammpling_L2} we have
\begin{equation}
\widetilde{\mathcal{B}}^0(\hat{s}^h) = \int_{\partial \Omega} \Psi^0_b(\bm{x}) \hat{s}^h(\bm{x}) \ \mathrm{d}\Gamma \;.
\end{equation}
With the expansions (\ref{eq:expansion_trace_q}) and (\ref{eq:expansion_trace_hat_s}) we can write the boundary integral in (\ref{eq:weak_grad_2D}) as
\begin{equation}
\int_{\partial \Omega} \hat{s} \left (\bm{q}^h \cdot \bm{n} \right ) \ \mathrm{d} \Gamma = {\mathcal{N}^1(\bm{q}^h)}^{\intercal} \mathbb{N}_1 \widetilde{\mathcal{B}}^0(\hat{s}^h) \;.
\end{equation}}


\MIGnew{Collecting all integrals in (\ref{eq:weak_grad_2D}) we have
\begin{equation}
{\mathcal{N}^1(\bm{q}^h)}^{\intercal}    \widetilde{\mathcal{N}}^1(\widetilde{\mbox{grad}}\ (s^h,\hat{s}^h)) = - {\mathcal{N}^1(\bm{q}^h)}^{\intercal} {\mathbb{E}^{2,1}}^{\intercal} \widetilde{\mathcal{N}}^0(s^h) + {\mathcal{N}^1(\bm{q}^h)}^{\intercal} \mathbb{N}_1 \widetilde{\mathcal{B}}^0(\hat{s}^h) \;. 
\end{equation}	
This equation needs to be satisfied for all $\bm{q}^h$, i.e. $\mathcal{N}^1(\bm{q}^h)$, so we need to have
\begin{equation}
    \widetilde{\mathcal{N}}^1(\widetilde{\mbox{grad}}\ (s^h,\hat{s}^h)) = -  {\mathbb{E}^{2,1}}^{\intercal} \widetilde{\mathcal{N}}^0(s^h) +  \mathbb{N}_1 \widetilde{\mathcal{B}}^0(\hat{s}^h) \;. 
    \label{eq:topological_grad}
\end{equation} 
}
%

\MIGnew{\begin{example}\label{ex:2D_inclusion_matrix_N1}
\begin{figure}[!htbp]
	\centering
		\includegraphics[width=.55\textwidth]{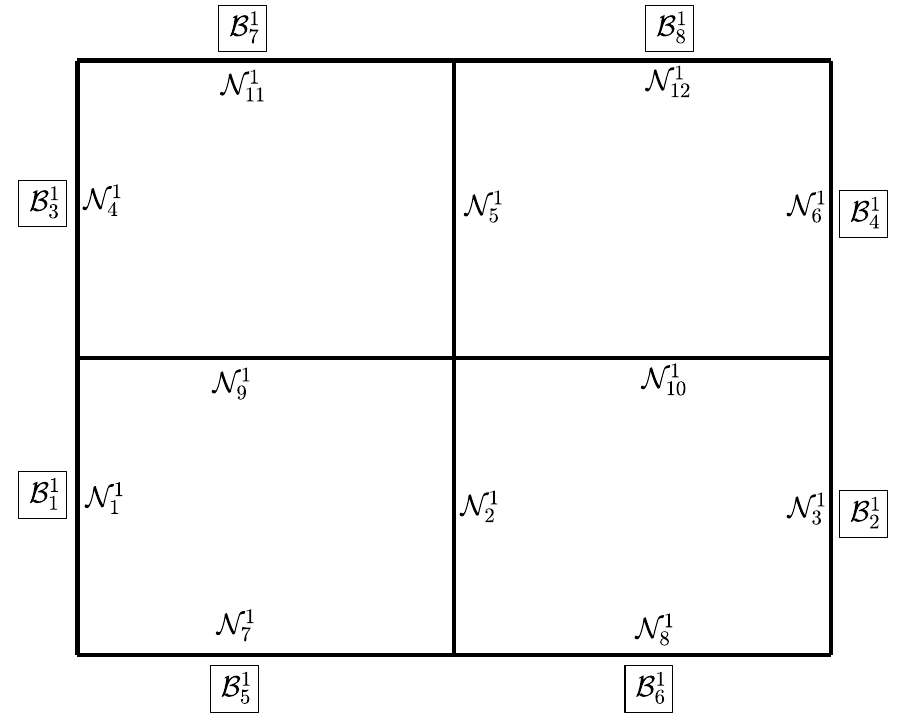}
		\caption{Global degrees of freedom $\mathcal{N}_i^1$ for vector fields in $D(\Omega)$ and the degrees of freedom $\mathcal{B}^1_j$ for functions in $D(\partial \Omega)$.}
		\label{fig:DOFs_and_BDOFs}
\end{figure}
In Figure~\ref{fig:DOFs_and_BDOFs} a small mesh is shown and the degrees of freedom $\mathcal{N}^1_i$ for elements in $D(\Omega)$ and the boundary degrees of freedom $\mathcal{B}_j^1$ for functions in $D(\partial \Omega)$. The matrix $\mathbb{N}_1$, which plays a role in (\ref{eq:topological_grad}), in this case is given by
\[
\mathbb{N}_1 = \left[ \begin{array}{cccccccc}
-1 & 0 & 0 & 0 & 0 & 0 & 0 & 0 \nl
0 & 0 & 0 & 0 & 0 & 0 & 0 & 0 \nl
0 & 1 & 0 & 0 & 0 & 0 & 0 & 0 \nl
0 & 0 & -1 & 0 & 0 & 0 & 0 & 0 \nl
0 & 0 & 0 & 0 & 0 & 0 & 0 & 0 \nl
0 & 0 & 0 & 1 & 0 & 0 & 0 & 0 \nl
0 & 0 & 0 & 0 & -1 & 0 & 0 & 0 \nl
0 & 0 & 0 & 0 & 0 & -1 & 0 & 0 \nl
0 & 0 & 0 & 0 & 0 & 0 & 0 & 0 \nl
0 & 0 & 0 & 0 & 0 & 0 & 0 & 0 \nl
0 & 0 & 0 & 0 & 0 & 0 & 1 & 0 \nl
0 & 0 & 0 & 0 & 0 & 0 & 0 & 1
\end{array} \right] \;.
\]
$\mathbb{N}_1$ is a $\mbox{dim}(D(\Omega))\times \mbox{dim}(D(\partial \Omega))$ matrix. On the eastern and northern boundary the outward unit normal points to the {\em positive}  $x$- and $y$-direction, respectively, and therefore the boundary degrees of freedom $\mathcal{B}_j^1$ equal the corresponding degrees of freedom along the boundary, $\mathcal{N}_i^1$. While for the western and southern boundary the outward unit normal points in the {\em negative} $x$- and $y$-direction, respectively, and therefore the boundary degree of freedom $\mathcal{B}_j^1$ is {\em minus} the degree of freedom $\mathcal{N}_i^1$. 
\end{example}	
}

\subsubsection{The rot acting on $\widetilde{D}(\Omega)\times \widetilde{C}(\partial \Omega)$}
The space of continuous, piecewise \MIGnew{functions} $C(\Omega)$ can be restricted to the boundary $\partial \Omega$, which we will refer to as $C(\partial \Omega)$. Using the mass matrix in this restricted function space, we can set up the dual representation for \VJ{this one dimensional space} as discussed in Section~\ref{sec:multi_elem_1D}. This dual representation will be denoted by $\widetilde{C}(\partial \Omega)$

As we have seen in Section ~\ref{sec:Multi-elem_2D} the piecewise polynomial vector fields in $\widetilde{D}(\Omega)$ only have continuity of the normal component over element boundaries; the tangential component will be discontinuous. So $\widetilde{D}\bb{\Omega}$ lacks the smoothness to apply the conventional curl operator and therefore we will define the curl as, $\widetilde{\mbox{rot}}:\widetilde{D}(\Omega) \times \widetilde{C}(\partial \Omega)\rightarrow \widetilde{C}(\Omega)$.
\begin{definition}
For $(\bm{d}^h,\hat{d}^h) \in \widetilde{D}(\Omega) \times \widetilde{C}(\partial \Omega)$ we have
\begin{equation}
\int_{\Omega} \widetilde{\VJ{\mathrm{rot}}} \ (\bm{d}^h,\hat{d}^h)\cdot c^h \ \mathrm{d}\Omega = \int_{\Omega} \bm{d}^h \cdot \VJ{\mathrm{curl}} \ c^h \ \mathrm{d} \Omega \ \varun{+} \int_{\partial \Omega} \hat{d}^h \cdot c^h \ \mathrm{d} \Gamma \qquad \qquad  \VJ{\forall c^h \in C\bb{\Omega}} \;.
\label{def:weak_rot}
\end{equation}
\end{definition}
Here we see that the integral on the left is a bilinear for on $\widetilde{C}(\Omega)\times C(\Omega)$, the first integral on the right hand side is a bilinear form on $\widetilde{D}(\Omega)\times D(\Omega)$, while the boundary integral is a bilinear form over $\widetilde{C}(\partial \Omega)\times C(\partial \Omega)$.

Expanding all variables in their appropriate basis functions reveals that
\begin{equation}
\widetilde{\mathcal{N}}^2(\widetilde{\mbox{rot}}\ (\bm{d}^h,\hat{d}^h))  =   {\mathbb{E}^{1,0}}^T \widetilde{\mathcal{N}}^1(\bm{d}^h) \ \varun{+} \ \MIGnew{\mathbb{N}_{0}} \VJ{\widetilde{\mathcal{B}}}^1(\hat{d}^h)   \;,
\label{eq:topological_rot}
\end{equation}
\VJ{where \MIGnew{$\mathbb{N}_{0}$} is the sparse inclusion matrix containing only the values $-1$, $0$ and $1$ that maps the boundary nodal degrees of freedom to the global nodes.
The subscript $'0'$ denotes the fact that this matrix is acting on the nodes (geometric dimension $'0'$).}

\MIGnew{\begin{example}\label{ex:2D_inclusion_matrix_N0}
		\begin{figure}[!htbp]
			\centering
			\includegraphics[width=.55\textwidth]{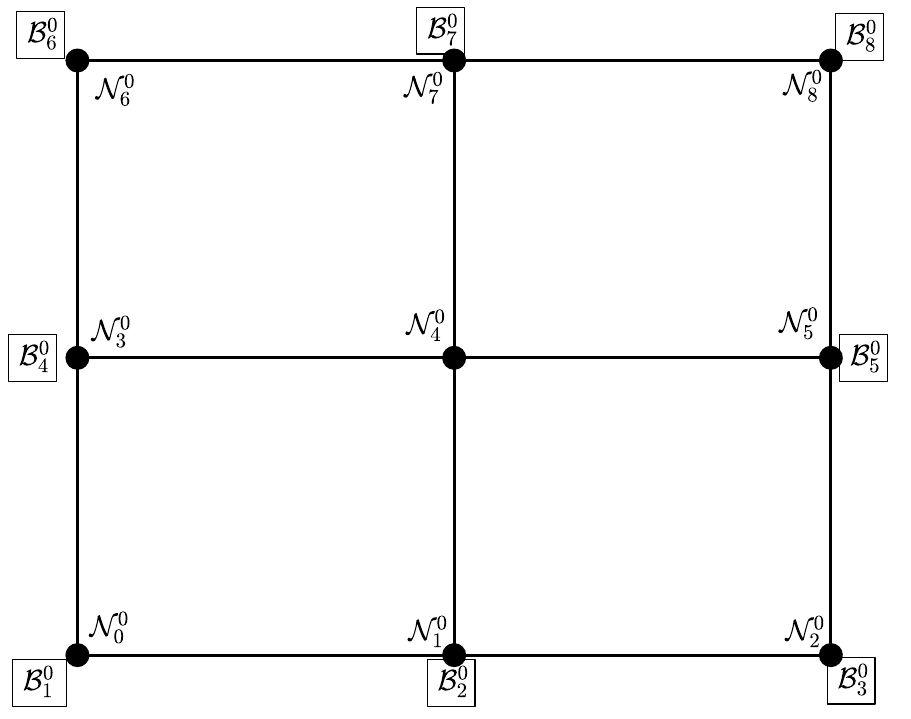}
			\caption{Global nodal degrees of freedom $\mathcal{N}_i^0$ for scalar fields in $C(\Omega)$ and the degrees of freedom $\mathcal{B}^0_j$ for functions in $C(\partial \Omega)$.}
			\label{fig:nodal_DOFs_and_BDOFs}
		\end{figure}
		In Figure~\ref{fig:nodal_DOFs_and_BDOFs} the small mesh is shown again but this time with the nodal degrees of freedom $\mathcal{N}^0_i$ for elements in $C(\Omega)$ and the boundary degrees of freedom $\mathcal{B}_j^0$ for functions in $C(\partial \Omega)$. The matrix $\mathbb{N}_0$, which plays a role in (\ref{eq:topological_rot}), in this case is given by
		\[ 
		\mathbb{N}_0 = \left[ \begin{array}{cccccccc} 
		1 & 0 & 0 & 0 & 0 & 0 & 0 & 0 \nl
		0 & 1 & 0 & 0 & 0 & 0 & 0 & 0 \nl
		0 & 0 & 1 & 0 & 0 & 0 & 0 & 0 \nl
		0 & 0 & 0 & 1 & 0 & 0 & 0 & 0 \nl
		0 & 0 & 0 & 0 & 0 & 0 & 0 & 0 \nl
		0 & 0 & 0 & 0 & 1 & 0 & 0 & 0 \nl
		0 & 0 & 0 & 0 & 0 & 1 & 0 & 0 \nl
		0 & 0 & 0 & 0 & 0 & 0 & 1 & 0 \nl
		0 & 0 & 0 & 0 & 0 & 0 & 0 & 1
		\end{array} \right]  \;.
		\]
		$\mathbb{N}_0$ is a $\mbox{dim}(C(\Omega))\times \mbox{dim}(C(\partial \Omega))$ matrix. Just like the edges at the outer boundary in Example~\ref{ex:2D_inclusion_matrix_N1} had an orientation given by the direction of the outward unit normal,  all points in the domain have an orientation given by the sense of rotation around the points which we have chosen to be anti-clockwise. Since we endow the boundary degrees of freedom with the same orientation, the matrix $\mathbb{N}_0$ only contains the non-zero value $1$.
	\end{example}	
}
\varun{
\begin{lemma}\label{lem:identity_in_boundary}
	For the inclusion matrices $\mathbb{N}_{k-1}$, $\mathbb{N}_{k}$ and the incidence matrix $\mathbb{E}^{k,k-1}$, $0<k<d$, we have the identity
	\begin{equation}
	\mathbb{N}_{k-1} \mathbb{N}_{k-1}^\intercal {\mathbb{E}^{k,k-1}}^\intercal \mathbb{N}_k = {\mathbb{E}^{k,k-1}}^\intercal \mathbb{N}_k \;. 
	\end{equation}
\begin{proof}
	The matrix $\mathbb{N}_{k-1} \mathbb{N}_{k-1}^\intercal$ maps global $(k-1)$-dimensional degrees of freedom in the primal representation to global $(k-1)$-dimensional degrees of freedom in the primal representation. It is an identity for the degrees of freedom on the boundary and sets the internal degrees of freedom to zero. Consequently, the matrix $\mathbb{I} - \mathbb{N}_{k-1} \mathbb{N}_{k-1}^\intercal$ sets the degrees of freedom along the boundary to zero and leaves the internal degrees of freedom unchanged. Then $\mathbb{E}^{k,k-1}$ maps the $(k-1)$-dimensional degrees of freedom  to the $k$-dimensional degrees of freedom. The $k$-dimensional degrees of freedom $\mathbb{E}^{k,k-1} \left ( \mathbb{I} - \mathbb{N}_{k-1} \mathbb{N}_{k-1}^\intercal \right )\mathcal{N}^{k-1}(\phi^h)$ will be zero along the boundary, since these degrees of freedom are linear combination of $(k-1)$-dimensional degrees of freedom along the boundary which are zero. The internal $k$-dimensional degrees of freedom will, in general, be non-zero. Premultiplication with $\mathbb{N}_k^\intercal$ restricts the $k$-dimensional degrees of freedom to the boundary and therefore
	\begin{equation*}
	\mathbb{N}_k^\intercal \mathbb{E}^{k,k-1} \left ( \mathbb{I} - \mathbb{N}_{k-1} \mathbb{N}_{k-1}^\intercal \right )\mathcal{N}^{k-1}(\phi^h) \equiv 0 \;,\quad \quad \forall \mathcal{N}^{k-1}(\phi^h) \;, 
	\end{equation*}
	which means that
	\begin{equation*}
	\mathbb{N}_k^\intercal \mathbb{E}^{k,k-1} \left ( \mathbb{I} - \mathbb{N}_{k-1} \mathbb{N}_{k-1}^\intercal \right ) \equiv 0 \;.
	\end{equation*}
	Taking the transpose of these matrices on both sides of the equality sign gives the desired result.
\end{proof}
\end{lemma}
A more constructive proof of this result will be given in Appendix~B.
}
\subsection{Extended derivatives of dual variables -- the dual de Rham sequence}
Straightforward calculation reveals that while $\mbox{div}\ \mbox{curl}\ g^h \equiv 0$ for $g^h \in G(\Omega)$, $\widetilde{\mbox{rot}}\ ( \widetilde{\mbox{grad}}\ (s^h,\hat{s}^h),\hat{d}^h )$ will in general not be identically zero for $s^h \in \widetilde{S}(\Omega)$, $\hat{s}^h \in \widetilde{D}(\partial \Omega)$ and $\hat{d}^h \in \widetilde{C}(\partial \Omega)$. Because for all $c^h \in C(\Omega)$ we have
\begin{eqnarray}
\int_{\Omega} \widetilde{\mbox{rot}}\ ( \widetilde{\mbox{grad}}\ (s^h,\hat{s}^h),\hat{d}^h ) \cdot c^h \ \mathrm{d} \Omega & \stackrel{(\ref{def:weak_rot})}{=} & \int_{\Omega} \widetilde{\mbox{grad}}\ (s^h,\hat{s}^h) \cdot \mbox{curl}\ c^h \ \mathrm{d}\Omega \ \varun{+} \ \int_{\MIGnew{\partial}\Omega} \hat{d}^h \cdot c^h \ \mathrm{d} \Gamma \nonumber \\
 & \stackrel{(\ref{eq:weak_grad_2D})}{=} & -\int_{\Omega} s^h \ \mbox{div}\ \mbox{curl} \ c^h \ \mathrm{d}\Omega + \int_{\partial \Omega} \hat{s}^h (\mbox{curl} \ c^h \cdot \bm{n}) \ \mathrm{d} \Gamma \ \varun{+} \ \int_{\MIGnew{\partial} \Omega} \hat{d}^h \ c^h \ \mathrm{d} \Gamma \label{eq:weak_rot_grad}\\
 & = & \int_{\partial \Omega} \left [ \hat{s}^h (\mbox{curl} \ c^h \cdot \bm{n}) \ \varun{+} \ \hat{d}^h \ c^h \right ] \ \mathrm{d} \Gamma \;, \nonumber
\end{eqnarray}
because the $\mbox{div}\ \mbox{curl} \equiv 0$.

Equation (\ref{eq:weak_rot_grad}) shows that the range of $\widetilde{\mbox{grad}}$ does not necessarily map into the null space of $\widetilde{\mbox{rot}}$. As a consequence the dual spaces and their vector operations do not form a de Rham sequence. We can remedy this by requiring that
\begin{equation}
\int_{\partial \Omega} \left [ \hat{s}^h (\mbox{curl} \ c^h \cdot \bm{n}) \ \varun{+} \ \hat{d}^h \ c^h \right ] \ \mathrm{d} \Gamma = 0 \;,\quad \quad \forall c^h \in C(\Omega) \;.
\end{equation}
Note that the integral of $\hat{s}^h$ and $(\mbox{curl} \ c^h \cdot \bm{n})$ is a bilinear form over $\widetilde{D}(\partial \Omega)\times D(\partial \Omega)$, while the integral of $\hat{d}^h$ and $c^h$ is a bilinear form over $\widetilde{C}(\partial \Omega)\times C(\partial \Omega)$. Since the integrals involve only mutual dual representations, they can be evaluated without reference to the basis functions and only in terms of the degrees of freedom. This gives
\begin{equation}
\int_{\partial \Omega} \left [ \hat{s}^h (\mbox{curl} \ c^h \cdot \bm{n} ) \ \varun{+} \ \hat{d}^h \ c^h \right ] \ \mathrm{d} \Gamma =
\MIGnew{{\mathcal{N}^0(c^h)}^{\intercal} {\mathbb{E}^{1,0}}^{\intercal} \mathbb{N}_1 \widetilde{\mathcal{B}}^0(\hat{s}^h) \ \varun{+} \ {\mathcal{N}^0(c^h)}^{\intercal} \mathbb{N}_0 \widetilde{\mathcal{B}}^1(\hat{d}^h) = 0 \;, \quad \quad \forall c^h \in C(\Omega) \;.}
\label{eq:deRham_in_boundary}
\end{equation}
\MIGnew{This implies that we need to have
\begin{equation} 
\mathbb{N}_0 \widetilde{\mathcal{B}}^1(\hat{d}^h)  = \ \varun{-} \ {\mathbb{E}^{1,0}}^{\intercal} \mathbb{N}_1 \widetilde{\mathcal{B}}^0(\hat{s}^h)\;.
\label{eq:rot_in_boundary}  
\end{equation} 
Since $\mathbb{N}_0^\intercal \mathbb{N}_0 = \mathbb{I}$, we have
\[  \widetilde{\mathcal{B}}^1(\hat{d}^h)  = \ \varun{-} \ \mathbb{N}_0^\intercal {\mathbb{E}^{1,0}}^{\intercal} \mathbb{N}_1 \widetilde{\mathcal{B}}^0(\hat{s}^h)\;.
\]}
%
%
\VJ{\begin{remark}
The above expression is synonymous to the operation of 1D \varun{dual} gradient on the elements of $\partial \Omega$ which is a 1D manifold.
This means that the trace spaces on $\partial \Omega$ also satisfy a \varun{1D} de Rham sequence.
\end{remark}}
\VJ{\begin{definition} \label{def:dual_GRAD_2d}
We now extend the definition of $\widetilde{\mbox{grad}}$ in \eqref{eq:weak_grad_2D} with the condition (\ref{eq:rot_in_boundary}) to a map $\VJ{\widetilde{\mbox{GRAD}}}: \widetilde{S}(\Omega) \times \widetilde{D}(\partial \Omega) \rightarrow \widetilde{D}(\Omega) \times \widetilde{C}(\partial \Omega)$.
Let $\bb{s^h, \hat{s}^h} \in \dualSb$ and $\bb{\bm{d}^h, \hat{d}^h} \in \dualDb$, then
\begin{eqnarray}
\widetilde{\mathrm{GRAD}}\bb{s^h, \hat{s}^h} & := & \bb{ \begin{array}{cc} \dualgrad \bb{s^h, \hat{s}^h} \;, & \widetilde{\mathrm{grad}}_b  \bb{\hat{s}^h} \end{array} } \nonumber \nl
& = & \left ( \begin{array}{cc}
\widetilde{\Psi}^1\bb{-{\mathbb{E}^{2,1}}^\intercal \widetilde{\mathcal{N}}^0(s^h)   +  \mathbb{N}_{\MIGnew{1}} \widetilde{\mathcal{\MIGnew{B}}}^0(\hat{s}^h)}, & \ \varun{-} \ \widetilde{\Psi}_b^1\ \mathbb{N}_0^\intercal {\mathbb{E}^{1,0}}^{\intercal} \mathbb{N}_1 \widetilde{\mathcal{\MIGnew{B}}}^0(\hat{s}^h)
\end{array} \right ) \nl
& = & \bb{\bm{d}^h, \hat{d}^h} \;, \nonumber
\end{eqnarray}
where $\widetilde{\Psi}_b^1 = \Psi^0 \mathbb{N}_0 \left ( \mathbb{M}_b^{0} \right )^{-1}$ are the dual basis functions defined over the boundary $\partial \Omega$.
\end{definition}}
\VJ{\begin{definition} \label{def:dual_ROT_2d}
We extend the $\widetilde{\mbox{rot}}$ operator in \eqref{eq:topological_rot} to a map $\dualDb \rightarrow \widetilde{C}\bb{\Omega} \times 0$.
Let $\bb{\bm{d}^h, \hat{d}^h} \in \dualDb$ and $\bb{f^h,0} \in \widetilde{C}\bb{\Omega} \times 0$, then
\begin{eqnarray}
\widetilde{\mathrm{ROT}}\bb{\bm{d}^h, \hat{d}^h} & := & \bb{ \begin{array}{cc} \widetilde{\mathrm{rot}} \bb{\bm{d}^h, \hat{d}^h} \;, & 0 \end{array} } \nonumber \nl
& = & \left ( \begin{array}{cc}
\widetilde{\Psi}^2\bb{{\mathbb{E}^{1,0}}^\intercal \widetilde{\mathcal{N}}^1(\bm{d}^h)   \ \varun{+} \  \mathbb{N}_{\MIGnew{0}} \widetilde{\mathcal{\MIGnew{B}}}^1(\hat{d}^h)},& 0
\end{array} \right )\label{eq:extended_rot_2d} \nl
& = & \bb{f^h,0} \;. \nonumber
\end{eqnarray} 
\end{definition}}
Using Definition~\ref{def:dual_GRAD_2d} and Definition~\ref{def:dual_ROT_2d} we have
\MIGnew{\begin{equation}
\widetilde{\mathcal{N}}^2(\widetilde{\mbox{ROT}} (\widetilde{\mbox{GRAD}}(s^h,\hat{s}^h)))  =  \left ( - {\mathbb{E}^{1,0}}^\intercal {\mathbb{E}^{2,1}}^\intercal \widetilde{\mathcal{N}}^0(s^h) + {\mathbb{E}^{1,0}}^\intercal \mathbb{N}_{1} \widetilde{\mathcal{B}}^0(\hat{s}^h) - \mathbb{N}_{0} \mathbb{N}_{0}^\intercal {\mathbb{E}^{1,0}}^\intercal \mathbb{N}_1 \widetilde{\mathcal{B}}^0(\hat{s}^h) \ , \  0 \right ) 
 =  \left ( 0\ ,\ 0 \right ) \;. 
\end{equation}}
Here we used the fact that $\mathbb{E}^{2,1} \mathbb{E}^{1,0} \equiv 0$ and \varun{Lemma~\ref{lem:identity_in_boundary}} that
\[\MIGnew{
\mathbb{N}_{0} \mathbb{N}_{0}^\intercal {\mathbb{E}^{1,0}}^\intercal \mathbb{N}_1 \widetilde{\mathcal{B}}^0(\hat{s}^h) = {\mathbb{E}^{1,0}}^\intercal \mathbb{N}_{1} \widetilde{\mathcal{B}}^0(\hat{s}^h) \;.}
\]
Therefore the $\mathcal{R}(\VJ{\widetilde{\mbox{GRAD}}};\widetilde{S}(\Omega)\times \widetilde{D}(\partial \Omega)) \subset \mathcal{K}(\VJ{\widetilde{\mbox{ROT}}};\widetilde{D}(\Omega)\times \widetilde{C}(\partial \Omega))$, where $\mathcal{K}(\VJ{\widetilde{\mbox{ROT}}};\widetilde{D}(\Omega)\times \widetilde{C}(\partial \Omega))$ denotes here the null space of the $\VJ{\widetilde{\mbox{ROT}}}$ operator applied to the space $\widetilde{D}(\Omega)\times \widetilde{C}(\partial \Omega)$.
So in order for dual operators to form a de Rham sequence, not only the derivative of the dual variables needs to be defined as was done in (\ref{eq:topological_grad}), but also the derivative of the trace variables, (\ref{eq:rot_in_boundary}).
\begin{definition} \label{def:norms_2d}
	\VJ{
		For $\widetilde{s}^h = \bb{s^h \;, \hat{s}^h} \in \widetilde{S}\bb{\Omega} \times \widetilde{D}\bb{\partial \Omega}$, we have
		\begin{equation}
		\| \widetilde{s}^h \|^2_{H \bb{\widetilde{\mathrm{grad}}; \Omega}} := 
		\| s^h \|^2_{L^2\bb{\Omega}} + \| \widetilde{\mathrm{grad}} \bb{s^h, \hat{s}^h} \|^2 _{L^2\bb{\Omega}}
		\end{equation}
		For $\widetilde{\bm{d}}^h \in \widetilde{D}\bb{\Omega} \times \widetilde{C}\bb{\partial \Omega}$, we have 
		\begin{equation}
		\| \widetilde{\bm{d}}^h \|^2_{H\bb{\widetilde{\mathrm{rot}} ;\Omega}} := 
		\| \bm{d}^h \|^2_{L^2\bb{\Omega}} + \| \widetilde{\mathrm{rot}} \bb{\bm{d}^h, \hat{d}^h} \|^2 _{L^2\bb{\Omega}}
		\end{equation}}
\MIGnew{\begin{corollary}
	Using the norms in Definition~\ref{def:norms_2d} and the degrees of freedom for $\widetilde{\mathrm{grad}}(s^h,\hat{s}^h)$ in (\ref{eq:topological_grad}) we have
	\[ 
	\| \widetilde{s}^h \|_{H(\widetilde{\mathrm{grad}};\Omega)}^2 = \widetilde{\mathcal{N}}^0(s^h)^\intercal \widetilde{\mathbb{M}}^{(0)} \widetilde{\mathcal{N}}^0(s^h) + \left ( \widetilde{\mathcal{N}}^0(s^h)^\intercal \mathbb{E}^{2,1} - \widetilde{\mathcal{B}}^0(\hat{s}^h )^\intercal \mathbb{N}_1^\intercal \right ) \widetilde{\mathbb{M}}^{(1)} \left ( {\mathbb{E}^{2,1}}^\intercal \widetilde{\mathcal{N}}^0(s^h) - \mathbb{N}_1 \widetilde{\mathcal{B}}^0(\hat{s}^h )  \right ) \;.
	\]
	Using the norms in Definition~\ref{def:norms_2d} and the degrees of freedom for $\widetilde{\mathrm{rot}}(\bm{d}^h,\hat{d}^h)$ in (\ref{eq:topological_rot}) we have
	\[ 
	\| \widetilde{\bm{d}}^h \|_{H(\widetilde{\mathrm{rot}};\Omega)}^2 = \widetilde{\mathcal{N}}^1(\bm{d}^h)^\intercal \widetilde{\mathbb{M}}^{(1)} \widetilde{\mathcal{N}}^1(\bm{d}^h) + \left ( \widetilde{\mathcal{N}}^1(\bm{d}^h)^\intercal \mathbb{E}^{1,0} - \widetilde{\mathcal{B}}^1(\hat{d}^h )^\intercal \mathbb{N}_0^\intercal \right ) \widetilde{\mathbb{M}}^{(2)} \left ( {\mathbb{E}^{1,0}}^\intercal \widetilde{\mathcal{N}}^1(\bm{d}^h) - \mathbb{N}_0 \widetilde{\mathcal{B}}^1(\hat{d}^h )  \right ) \;.
	\]
\end{corollary}}
\end{definition}
\VJ{\begin{definition}
		We define the norm in the dual of the trace spaces following \cite[eq. 1.16]{1986Girault} as: \\
		For $\hat{s}^h \in \widetilde{D}\bb{\partial \Omega}$, and $\widetilde{s}^h = \bb{s^h \;, \hat{s}^h} \in \widetilde{S}\bb{\Omega} \times \widetilde{D}\bb{\partial \Omega}$we have
		\begin{equation} \label{eq:grad_norm_2d}
		\| \hat{s}^h \|_{\widetilde{D}(\partial \Omega)} := \inf_{\varun{s^h \in \widetilde{S}(\Omega)}} \| \widetilde{s}^h \|_{H\bb{\widetilde{\mathrm{grad}}; \Omega}} \;.
		\end{equation}
		For $\hat{\bm{d}}^h \in \widetilde{C}\bb{\partial \Omega}$, and $\widetilde{\bm{d}}^h = \bb{\bm{d}^h \;, \hat{d}^h} \in \dualDb $ we have
		\begin{equation}
		\| \hat{\bm{d}}^h \|_{\widetilde{C}(\partial \Omega)} := \inf_{ \varun{ \bm{d}^h \in \widetilde{D}(\Omega)} } \| \widetilde{\bm{d}}^h \|_{H\bb{\widetilde{\mathrm{rot}};\Omega}} \;.
		\end{equation}
\end{definition}}
\subsection{Discrete de Rham sequence}
Using the function spaces defined in \secref{sec:Multi-elem_2D} and \secref{sec:dula_operators_2d}, the de Rham sequence for primal spaces is given as
\begin{equation}
\begin{tikzcd}[row sep=1.0cm, column sep = 2cm]
C\bb{\Omega} \arrow{r}{\mathrm{curl}} & D\bb{\Omega} \arrow{r}{\mathrm{div}} & S\bb{\Omega}
\end{tikzcd} \;.
\end{equation}
The de Rham sequence for the dual spaces is given as
\begin{equation}
\begin{tikzcd}[row sep=1.0cm, column sep = 2cm]
\widetilde{S}\bb{\Omega} \times \widetilde{D}\bb{\partial \Omega} \arrow{r}{\VJ{\widetilde{\mathrm{GRAD}}}} & \widetilde{D}\bb{\Omega} \times \widetilde{C}\bb{\partial \Omega} \arrow{r}{\VJ{\widetilde{\mathrm{ROT}}}} & \widetilde{C}\bb{\Omega} \times \VJ{0}
\end{tikzcd} \;.
\end{equation}
\section{Three-dimensional dual spaces}\label{sec:3D_dual_spaces}
Similar to \secref{sec:2D_dual_spaces} we can define primal and dual spaces for three dimensional problems.
For $\Omega \subset \mathbb{R}^3$ the de Rham sequence is given by
\begin{equation*}
\begin{tikzcd}[row sep=1.0cm, column sep = 2cm]
H^1 \bb{{\Omega}} \arrow{r}{\mathrm{grad}} & H \bb{\mathrm{curl};{\Omega}} \arrow{r}{\mathrm{curl}} & H \bb{\mathrm{div};{\Omega}} \arrow{r}{\mathrm{div}} & L^2 \bb{{\Omega}}
\end{tikzcd} \;.
\end{equation*}
We will {first} define the four finite element spaces for {$K \in [-1,1]^3$}, $G (K) \subset H^1(K)$,
$C(K) \subset H(\mathrm{curl};K)$, $D(K) \subset H(\mathrm{div};K)$, $S(K) \subset L^2(K)$ and the corresponding dual spaces $\widetilde{G} (K)$, $\widetilde{C} (K)$, $\widetilde{D} (K)$, $\widetilde{S} (K)$.
\VJ{We will then extend this construction to an arbitrary $\Omega _k \in \mathbb{R}^3$.
	This will be followed by an explanation on the treatment of} multiple element case and derivation of differential operators on the dual spaces.

Let $\xi_i, \eta _j, \zeta _k$, $i,j,k = 0,\hdots,N$, be the GLL \varun{nodes}. 
Let $\mathcal{P}$ be the space of polynomials of degree $N$ as introduced in \exampleref{ex:Lagrange_1D} and $\mathcal{Q}$ be the space of polynomials of degree $N-1$ as introduced in \exampleref{ex:edge_1D}.

\subsection{The function space $\VJ{G}(K) \subset H^1(K)$} \label{sec:ghk_3D}
Let $G(K) := \mathcal{P} \otimes \mathcal{P} \otimes \mathcal{P}$ be tensor product space.
Any element $p^{{h}} \in G(K)$ can be represented by
\begin{equation*}
p^{{h}}(\xi, \eta, \zeta) = \Psi ^0 \bb{\xi, \eta, \zeta} \mathcal{N}^0(p^{{h}}) \;,
\end{equation*}
where $N^0(p^{{h}})$ are the degrees of freedom defined at set of 3D nodes $\bm{x}_{i(N+1)^2 + j(N+1) + k} = \bb{\xi _i, \eta _j, \zeta _k}$, $i, j, k = 0, \hdots, N$ given by,
\begin{equation*}
\mathcal{N}^0 \bb{p^{{h}}} := p^{{h}}\bb{\bm{x}_k} \;, \qquad k = 0, \hdots , \bb{N+1}^3 - 1 \;,
\end{equation*}
and $\Psi^0\bb{\xi, \eta, \zeta} $ are the Lagrange (or nodal) basis through $\bm{x}_k$ given by
\begin{equation*}
\Psi ^0 \bb{\xi, \eta, \zeta} = \bb{\begin{array}{cccc}
	\epsilon ^{(0)}_0\bb{\xi, \eta, \zeta}  & \epsilon ^{(0)}_1\bb{\xi, \eta, \zeta} & \hdots & \epsilon ^{(0)}_{(N+1)^3 -1}\bb{\xi, \eta, \zeta}
	\end{array}} \;,
\end{equation*}
where
\begin{equation*}
\epsilon ^{(0)}_{i(N+1)^2 + j(N+1) + k} \bb{\xi, \eta, \zeta} := h_i\bb{\xi} h_j \bb{\eta} h_k \bb{\zeta} \quad i,j,k = 0, \hdots, N \;.
\end{equation*}
\begin{definition}
	The dual degrees of freedom for $p^h \in \widetilde{G}(K)$ are given by
	\begin{equation*}
	\widetilde{\mathcal{N}}^{3}(p^h) :=  \mathbb{M}^{(0)}\mathcal{N}^{0}(p^h) \qquad \mbox{where} \qquad \mathbb{M}^{(0)} = \int _K \Psi ^0 \bb{\xi, \eta, \zeta} ^\intercal \Psi ^0 \bb{\xi, \eta, \zeta} \d K\;,
	\label{eq:dual_volume_DOF_3D}
	\end{equation*}
	and the associated dual basis, following Corollary \ref{cor:dual_nodal_basis_functions}, are given by
	\begin{equation*}
	\widetilde{\Psi}^{{3}}\bb{\xi, \eta, \zeta} = \Psi ^0\bb{\xi, \eta, \zeta} \bb{\mathbb{M}^{(0)}}^{-1} \;.
	\end{equation*}
\end{definition}
\subsection{The function space $C(K) \subset H\bb{\mathrm{curl};K}$} \label{sec:chk_3D}
Consider the polynomial tensor spaces given by, $C_{\xi} := \mathcal{Q} \otimes \mathcal{P} \otimes \mathcal{P}$, $C_{\eta} := \mathcal{P} \otimes \mathcal{Q} \otimes \mathcal{P}$ and $C_{\zeta} := \mathcal{P} \otimes \mathcal{P} \otimes \mathcal{Q}$.
We define the finite element space in 3D as $C(K) := C_{\xi} \times C_{\eta} \times C_{\zeta}$.
Any vector field $\bm{c}^h \in C(K)$ can be represented as
\begin{equation*}
\bm{c}^{{h}} \bb{\xi, \eta, \zeta} = \Psi ^1 \bb{\xi, \eta, \zeta} \VJ{\mathcal{N}}^1(\bm{c}^{{h}}) \;,
\end{equation*}
where the the degrees of freedom defined at edges are given by
\begin{equation*}
\mathcal{N}^1 (\bm{c}^{{h}}) = \left\lbrace \begin{array}{lll}
\VJ{\mathcal{N}^1}_{(i-1)(N+1)^2 + j(N+1) + k} \bb{\bm{c}^{{h}}} & := {\int_{\xi _{i-1}}^{\xi _i}} \bb{\bm{c}^{{h}}{({\xi , \eta_ j, \zeta _k})} \cdot \bm{e}_{\xi }}\ \d \xi, & i = 1, \hdots,N; j,k = 0, \hdots, N \nl
\VJ{\mathcal{N}^1}_{N(N+1)^2+iN(N+1) + (j-1)(N+1) + k} \bb{\bm{c}^{{h}}} & := {\int_{\eta_{j-1}}^{\eta_j} \bb{\bm{c}^{{h}}\bb{\xi _i, \eta, \zeta_k} \cdot \bm{e}_{\eta}}}\ \d \eta, & i,k = 0, \hdots,N; j = 1, \hdots, N \nl
\VJ{\mathcal{N}^1}_{2N(N+1)^2+iN(N+1) + jN + k-1} \bb{\bm{c}^{{h}}} & :={ \int_{\zeta_{k-1}}^{\zeta_k} \bb{\bm{c}^{{h}}\bb{\xi _i, \eta _j, \zeta} \cdot \bm{e}_{\zeta}}}\ \d \zeta, & i,j = 0, \hdots,N; k = 1, \hdots,N
\end{array}  \right. \;,
\end{equation*}
and the basis that satisfies the Kronecker delta property, in the sense of integrals over the edges, are given by
\begin{equation*}
\Psi ^1 \bb{\xi, \eta, \zeta} = \bb{\begin{array}{cccc}
	\epsilon ^{(1)}_0\bb{\xi, \eta, \zeta}  & \epsilon ^{(1)}_1\bb{\xi, \eta, \zeta} & \hdots & \epsilon ^{(1)}_{6N(N+1)^2 -1}\bb{\xi, \eta, \zeta}
	\end{array}} \;,
\end{equation*}
where,
\begin{equation*}
\left\lbrace \begin{array}{lll}
\bm{\epsilon}^{(1)}_{(i-1)(N+1)^2 + j(N+1) + k} \bb{\xi, \eta, \zeta} & := e_i\bb{\xi} h_j \bb{\eta} h_k \bb{\zeta} \VJ{\bm{e}_\xi} & i = 1, \hdots, N; j,k = 0, \hdots, N \nl
\bm{\epsilon}^{(1)}_{N(N+1)^2+iN(N+1) + (j-1)(N+1) + k} \bb{\xi, \eta, \zeta} & := h_i\bb{\xi} e_j \bb{\eta} h_k \bb{\zeta} \VJ{\bm{e}_\eta} & i,k = 0, \hdots, N; j = 1, \hdots, N \nl
\bm{\epsilon}^{(1)}_{2N(N+1)^2+iN(N+1) + jN + k-1} \bb{\xi, \eta, \zeta} & := h_i\bb{\xi} h_j \bb{\eta} e_k \bb{\zeta} \VJ{\bm{e}_\zeta}& i,j = 0, \hdots, N; k = 1, \hdots, N
\end{array} \right. \;.
\end{equation*}
Let $p^h \in G(K)$, then following Corollary\ref{cor:derivative_incidence}  $\mathrm{grad}\ p^h \in C\bb{K}$ is given by
\begin{equation*}
\mbox{grad}\,p^h(\xi,\eta, \zeta) = \Psi^1(\xi,\eta,\zeta) \mathbb{E}^{1,0} \mathcal{N}^0(p^h) \;,
\end{equation*}
expanded in the basis functions of $C(K)$,
\varun{where}
{the incidence matrix $\mathbb{E}^{1,0}$ is for \VJ{degrees of freedom of} a 3D element\footnote{{The incidence matrix $\mathbb{E}^{1,0}$ is not the same as \varun{used} earlier in \eqref{eq:incidence_matrix_1D} for \VJ{degrees of freedom of a} 1D element or in \eqref{eq:curl_incidence} for \VJ{degrees of freedom of a} 2D element.}}.} 
This implies that $\mathcal{R}\bb{\mathrm{grad};G(K)} \subset C(K)$.
\begin{definition}
	The dual degrees of freedom in function space $\widetilde{C}(K)$ are given by
	\begin{equation*} \label{eq:dual_dof_N2_3d}
	\widetilde{\mathcal{N}}^{2}(\bm{c}^{{h}}) :=  \mathbb{M}^{(1)}\mathcal{N}^{1}(\bm{c}^{{h}}) \qquad \mbox{where} \qquad \mathbb{M}^{(1)} = \int _K \Psi ^1 \bb{\xi, \eta, \zeta} ^\intercal \Psi ^1 \bb{\xi, \eta, \zeta} \d K \;,
	\end{equation*}
	and the dual basis are given by, {Corollary \ref{cor:dual_nodal_basis_functions}},
	\begin{equation*}
	\widetilde{\Psi}^{{2}}\bb{\xi, \eta, \zeta} {:=} \Psi ^1\bb{\xi, \eta, \zeta} \bb{\mathbb{M}^{(1)}}^{-1} \;.
	\end{equation*}
\end{definition}
\subsection{The function space $D(K)  \subset H\bb{\mathrm{div};K}$} \label{sec:dhk_3D}
We define the finite element space as $D(K) := D_{\xi} \times D_{\eta} \times D_{\zeta}$, where, $D_{\xi} := \mathcal{P} \otimes \mathcal{Q} \otimes \mathcal{Q}$, $D_{\eta} := \mathcal{Q} \otimes \mathcal{P} \otimes \mathcal{Q}$ and $D_{\zeta} := \mathcal{Q} \otimes \mathcal{Q} \otimes \mathcal{P}$.
Any vector field $\bm{q}^{{h}} \in D(K)$ can be expressed as
\begin{equation} \label{eq:expansion_q}
\bm{q}^{{h}} \bb{\xi, \eta, \zeta} = \Psi ^2 \bb{\xi, \eta, \zeta} \mathcal{N}^2 \bb{\bm{q}^{{h}}} \;,
\end{equation}
where, the degrees of freedom defined at the surfaces are given by
\begin{equation*}
\mathcal{N}^2 (\bm{q}^{{h}}) = \left\lbrace \begin{array}{lll}
\mathcal{N}^2_{iN^2 + (j-1)N + k-1} \bb{\bm{q}^{{h}}} & :={ \int_{\eta _{j-1}}^{\eta _j} \int_{\zeta _{k-1}}^{\zeta _k} \bb{\bm{q}^{{h}}\bb{\xi	_i, \eta, \zeta} \varun{\cdot} \ \bm{e}_{\xi}}}  \d \eta\ \d \zeta & i = 0, \hdots, N; j,k = 1, \hdots, N \nl
\mathcal{N}^2_{N^2(N+1)+(i-1)N(N+1) + jN + k-1} \bb{\bm{q}^{{h}}} & := {\int_{\xi _{i-1}}^{\xi _i} \int_{\zeta _{k-1}}^{\zeta _k} \bb{\bm{q}^{{h}}\bb{\xi, \eta _j, \zeta} \varun{\cdot} \ \bm{e}_{\eta}}}  \d \xi\ \d \zeta & i,k = 1, \hdots, N; j = 0, \hdots, N \nl
\mathcal{N}^2_{2N^2(N+1)+(i-1)N(N+1) + (j-1)(N+1) + k} \bb{\bm{q}^{{h}}} & := {\int_{\xi _{i-1}}^{\xi _i} \int_{\eta _{j-1}}^{\eta _j} \bb{\bm{q}^{{h}}\bb{\xi	, \eta, \zeta _k} \varun{\cdot} \ \bm{e}_{\zeta}}}  \d \xi\  \d \eta & i,j = 1, \hdots, N; k = 0, \hdots, N 
\end{array}  \right. \;,
\end{equation*}
and the basis function that satisfy the Kronecker delta property are given by
\begin{equation*}
\Psi ^2 \bb{\xi, \eta, \zeta} = \bb{\begin{array}{cccc}
	\VJ{\bm{\epsilon}} ^{(2)}_0\bb{\xi, \eta, \zeta}  & \VJ{\bm{\epsilon}} ^{(2)}_1\bb{\xi, \eta, \zeta} & \hdots & \VJ{\bm{\epsilon}} ^{(2)}_{3N^2(N+1)-1}\bb{\xi, \eta, \zeta}
	\end{array}} \;,
\end{equation*}
where,
\begin{equation*}
\left\lbrace \begin{array}{lll}
\bm{\epsilon}^{(2)}_{iN^2 + (j-1)N + k-1} \bb{\xi, \eta, \zeta} & := h_i\bb{\xi} e_j \bb{\eta} e_k \bb{\zeta} \varun{\bm{e}_\xi} & i = 0, \hdots, N; j,k = 1, \hdots, N \nl
\bm{\epsilon}^{(2)}_{N^2(N+1)+(i-1)N(N+1) + jN + k-1} \bb{\xi, \eta, \zeta} & := e_i\bb{\xi} h_j \bb{\eta} e_k \bb{\zeta} \varun{\bm{e}_\eta} & i,k = 1, \hdots, N; j = 0, \hdots, N \nl
\bm{\epsilon}^{(2)}_{2N^2(N+1)+(i-1)N(N+1) + (j-1)(N+1) + k} \bb{\xi, \eta, \zeta} & := e_i\bb{\xi} e_j \bb{\eta} h_k \bb{\zeta} \varun{\bm{e}_\zeta} & i,j = 1, \hdots, N; k = 0, \hdots, N
\end{array} \right. \;.
\end{equation*}
Let $\bm{c}^h \in C(K)$, then the $\mathrm{curl}$ operation on $\bm{c}^h$, following Corollary \ref{cor:derivative_incidence}, gives
\begin{equation*}
\mbox{curl}\,\bm{{c}}^h(\xi,\eta, \zeta) = \Psi^2(\xi,\eta,{\zeta}) \mathbb{E}^{2,1} \mathcal{N}^1(\bm{{c}}^h) \;,
\end{equation*}
expanded in the basis functions of $D(K)$.
{Here the incidence matrix $\mathbb{E}^{2,1}$ is defined for \VJ{degrees of freedom of} a 3D element\footnote{{The incidence matrix $\mathbb{E}^{2,1}$ is not the same as \varun{used earlier} in \eqref{eq:incidence_21} for \VJ{degrees of freedom of a} 2D element.}}.}
This implies that $\mathcal{R}\bb{\mathrm{curl};C(K)} \subset D(K)$.
\begin{definition}
	The dual degrees of freedom and dual basis functions for $\bm{p}^h \in \widetilde{D}\bb{K}$ are given by
	\begin{equation} \label{eq:mass_matrix_M2_3d}
	\widetilde{\mathcal{N}}^{1}(\bm{p}^h) :=  \mathbb{M}^{(2)}\mathcal{N}^{2}(\bm{p}^h) \qquad \mbox{where} \qquad \mathbb{M}^{(2)} = \int _K \Psi ^2 \bb{\xi, \eta, \zeta} ^\intercal \Psi ^2 \bb{\xi, \eta, \zeta} \d K\;,
	\end{equation}
	\begin{equation*}
	\widetilde{\Psi}^{{1}}\bb{\xi, \eta, \zeta} {:=} \Psi ^2\bb{\xi, \eta, \zeta} \bb{\mathbb{M}^{(2)}}^{-1} \;.
	\end{equation*}
\end{definition}
\subsection{The function space $S(K)$}  \label{sec:shk_3D}
In this case, we define the finite element space in 3D as $S(K) := \mathcal{Q} \times \mathcal{Q} \times \mathcal{Q}$.
We can express any polynomial $f^{{h}}\bb{\xi, \eta, \zeta} \in S(K)$ as
\begin{equation*}
f^{{h}} \bb{\xi, \eta, \zeta} = \Psi ^3\bb{\xi, \eta, \zeta} \mathcal{N}^3(f^{{h}}) \;,
\end{equation*}
where the degrees of \VJ{freedom} evaluated over a volume are given by
\begin{equation} \label{eq:dof_f}
\mathcal{N}^3_{(i-1)N^2 + (j-1)N + k-1} {(f^h)} := \int_{\xi _{i-1}}^{\xi _i} \int_{\eta _{j-1}}^{\eta _j} \int _{\zeta _{k-1}}^{\zeta _k} f^{{h}} {\bb{\xi, \eta, \zeta}} \  \d \xi\ \d \eta\ \d \zeta  \qquad \mbox{for} \qquad i,j,k = 1, \hdots, N \;,
\end{equation}
and the basis functions that satisfy the Kronecker-delta property are given by
\begin{equation*}
\Psi ^3 \bb{\xi, \eta, \zeta} = \bb{\begin{array}{cccc}
	\epsilon ^{(3)}_0\bb{\xi, \eta, \zeta}  & \epsilon ^{(3)}_1\bb{\xi, \eta, \zeta} & \hdots & \epsilon ^{(3)}_{N^3 -1}\bb{\xi, \eta, \zeta}
	\end{array}} \;,
\end{equation*}
where,
\begin{equation*}
\epsilon _{(i-1)N^2 + (j-1)N + k-1} \bb{\xi, \eta, \zeta} := e_i\bb{\xi} e_j \bb{\eta} e_k \bb{\zeta} \qquad i,j,k = 1, \hdots, N \;.
\end{equation*}
Let $\bm{q}^h \in D(K)$, then the $\mathrm{div}$ operation on $\bm{q}^h$ gives
\begin{equation}
\mbox{div}\,\bm{q}^h(\xi,\eta, \zeta) = \Psi^3(\xi,\eta, \zeta) \mathbb{E}^{3,2} \mathcal{N}^2(\bm{q}^h) \;,
\label{eq:incidence_21_3D}
\end{equation}
expressed in terms of the basis functions of $S(K)$.
This implies that $\mathcal{R}\bb{\mathrm{div};D(K)} = S(K)$.
\begin{definition}
	The dual degrees of freedom and the dual basis for $f^h \in \widetilde{S}\bb{K}$ are given by,
	\begin{equation} \label{eq:mass_matrix_M3_3d}
	\widetilde{\mathcal{N}}^{0}(f^h) :=  \mathbb{M}^{(3)}\mathcal{N}^{3}(f^h) \qquad \mbox{where} \qquad \mathbb{M}^{(3)} = \int _K \Psi ^3 \bb{\xi, \eta, \zeta} ^\intercal \Psi ^3 \bb{\xi, \eta, \zeta} \d K \;,
	\end{equation}
	\begin{equation*} \label{eq:dual_volumes}
	\widetilde{\Psi}^{{0}}\bb{\xi, \eta, \zeta} {:=} \Psi ^3\bb{\xi, \eta, \zeta} \bb{\mathbb{M}^{(3)}}^{-1}  \;.
	\end{equation*}
\end{definition}
So far we have introduced finite element spaces on the primal and the dual complex, and the differential operators on the primal complex for a reference element $K \in [-1,1]^3$.
\VJ{In the next sections we will present transformation rules \varun{for} spaces \varun{of} an arbitrary $\Omega _k \in \mathbb{R}^3$}.
\VJ{
	\subsection{Transformation rules for 3D function spaces}	
	In this section we introduce the function spaces \varun{for an arbitrary} $\Omega _k \subset \mathbb{R}^3$. The transformations are performed element-by-element provided that the global transformation is continuous.
	Let $\Phi _k$ be the diffeomorphism and $\boldsymbol{\mathsf{J}}$ be the Jacobian tensor such that,
	\[ \Phi_{k}: \bb{\xi , \eta , \zeta} \in K \mapsto \bb{x,y,z} \in \Omega _k \qquad \mbox{and} \qquad \boldsymbol{\mathsf{J}} : = 
	\left[
	\begin{array}{ccc}
	\frac{\partial\Phi_{k}^{x}}{\partial \xi} & \frac{\partial\Phi_{k}^{x}}{\partial \eta} & \frac{\partial\Phi_{k}^{x}}{\partial \zeta} \\
	\frac{\partial\Phi_{k}^{y}}{\partial \xi} & \frac{\partial\Phi_{k}^{y}}{\partial \eta} & \frac{\partial\Phi_{k}^{y}}{\partial \zeta} \\
	\frac{\partial\Phi_{k}^{z}}{\partial \xi} & \frac{\partial\Phi_{k}^{z}}{\partial \eta} & \frac{\partial\Phi_{k}^{z}}{\partial \zeta}
	\end{array}
	\right]. \]
	\subsubsection{\VJ{Transformation rules for spaces $G\bb{\Omega_k}$}}
	Given a function $\bar{f} \in G(K)$ on $K$, we define $f \in G(\Omega _k)$ on $\Omega _k$ by
	\begin{equation*}
	f := \left(\Phi^{*}_{k}\right)^{-1} \left[\bar{f}\right] = \bar{f}\circ\Phi_{k}^{-1}\,,
	\end{equation*}
	where $\Phi^{*}_k$ is the pullback operator. We can then reverse the relation to obtain
	\begin{equation*}
	\bar{f} := \Phi^{*}_{k} \left[f\right] = f\circ\Phi_{k}\,.
	\end{equation*}
	Given two functions $f, g \in G(\Omega _k)$ on $\Omega _k$ their inner product can be computed by
	\begin{equation*}
	\bb{ f, g}_{\Omega _k} := \int_{\Omega _k} f g\, \mathrm{d}\Omega =  \int_{\Phi_{k}\left(K\right)} fg\,\mathrm{d}\Omega = \int_{K} \Phi^{*}_{k} \left[f\right] \Phi^{*}_{k} \left[g\right] \,\mathrm{det}\left(\boldsymbol{\mathsf{J}}\right)\mathrm{d}K= \int_K \bar{f} \bar{g}\,\mathrm{det}\left(\boldsymbol{\mathsf{J}}\right)\mathrm{d}K\,.
	\end{equation*}
	\subsubsection{Transformation rules for spaces $C\bb{\Omega _k}$}
	Given a vector field $\bar{\boldsymbol{v}} \in C(K)$, the transformed vector field $\boldsymbol{v} \in C(\Omega _k)$ on $\Omega _k$ by the following expression
	\begin{equation*}
	\boldsymbol{v} :=  \left(\Phi^{*}_{k}\right)^{-1} \left[\bar{\boldsymbol{v}}\right] = \left(\boldsymbol{\mathsf{J}}^{\intercal}\circ\Phi_{k}^{-1}\right)^{-1}\left(\bar{\boldsymbol{v}}\circ\Phi_{k}^{-1}\right)\,.
	\end{equation*}
	Following a similar approach as before, the inverse of this transformation is given by 
	\begin{equation*}
	\bar{\boldsymbol{v}} :=  \Phi^{*}_{k}\left[\boldsymbol{v}\right] = \boldsymbol{\mathsf{J}}^{\intercal}\left(\boldsymbol{v}\circ\Phi_{k}\right)\,.
	\end{equation*}
	Given two vector fields $\boldsymbol{u}, \boldsymbol{v} \in C(\Omega _k)$ on $\Omega _k$ their inner product can be computed by
	\begin{equation*}
	\bb{ \boldsymbol{u}, \boldsymbol{v}}_{\Omega _k} := \int_{\Omega _k} \boldsymbol{u}^{\intercal} \boldsymbol{v}\, \mathrm{d}\Omega  = \int_{K} \left(\Phi^{*}_{k} \left[\boldsymbol{u}\right]\right)^{\intercal} \left(\boldsymbol{\mathsf{J}}^{\intercal}\boldsymbol{\mathsf{J}}\right)^{-1} \Phi^{*}_{k} \left[\boldsymbol{v}\right] \,\mathrm{det}\left(\boldsymbol{\mathsf{J}}\right)\mathrm{d}K= \int_{K} \bar{\boldsymbol{u}}^{\intercal}  \left(\boldsymbol{\mathsf{J}}^{\intercal}\boldsymbol{\mathsf{J}}\right)^{-1}  \bar{\boldsymbol{v}}\,\mathrm{det}\left(\boldsymbol{\mathsf{J}}\right)\mathrm{d}K\,.
	\end{equation*}
	\subsubsection{Transformation rules for spaces $D\bb{\Omega _k}$}
	Given a vector field $\bar{\boldsymbol{u}} \in D(K)$, the transformed vector field $\boldsymbol{u} \in D(\Omega _k)$ on $\Omega _k$ \VJ{is given} by the following expression
	\begin{equation*}
	\boldsymbol{u} :=  \left(\Phi^{*}_{k}\right)^{-1} \left[\bar{\boldsymbol{u}}\right] = \frac{1}{\det\left(\boldsymbol{\mathsf{J}}\circ\Phi_{k}^{-1}\right)}\left(\boldsymbol{\mathsf{J}}\circ\Phi_{k}^{-1}\right)\left(\bar{\boldsymbol{u}}\circ\Phi_{k}^{-1}\right)\,.
	\end{equation*}
	It is possible to compute the inverse of this transformation, resulting in
	\begin{equation*}
	\bar{\boldsymbol{u}} :=  \Phi^{*}_{k}\left[\boldsymbol{u}\right] = \det\left(\boldsymbol{\mathsf{J}}\right)\boldsymbol{\mathsf{J}}^{-1}\left(\boldsymbol{u}\circ\Phi_{k}\right)\,.
	\end{equation*}
	Given two vector fields $\boldsymbol{u}, \boldsymbol{v} \in D(\Omega _k)$ on $\Omega _k$ their inner product can be computed by
	\begin{equation*}
	\bb{\boldsymbol{u}, \boldsymbol{v}}_{\Omega _k} := \int_{\Omega _k} \boldsymbol{u}^{\intercal} \boldsymbol{v}\, \mathrm{d}\Omega  = \int_K \left(\Phi^{*}_{k} \left[\boldsymbol{u}\right]\right)^{\intercal} \boldsymbol{\mathsf{J}}^{\intercal}\boldsymbol{\mathsf{J}}\, \Phi^{*}_{k} \left[\boldsymbol{v}\right] \frac{1}{\mathrm{det}\left(\boldsymbol{\mathsf{J}}\right)}\,\mathrm{d}K= \int_K \bar{\boldsymbol{u}}^{\intercal} \, \boldsymbol{\mathsf{J}}^{\intercal}\boldsymbol{\mathsf{J}}\,  \bar{\boldsymbol{v}}\,\frac{1}{\mathrm{det}\left(\boldsymbol{\mathsf{J}}\right)}\,\mathrm{d}K\,.
	\end{equation*}
	\subsubsection{\VJ{Transformation rules for spaces $S\bb{\Omega _k}$}}
	Given a function $\bar{g} \in S(K)$ on $K$, we define $g \in S(\Omega _k)$ on $\Omega _k$ by
	\begin{equation*}
	g := \left(\Phi^{*}_{k}\right)^{-1} \left[\bar{g}\right] = \frac{1}{\det\left(\boldsymbol{\mathsf{J}}\circ\Phi_{k}^{-1}\right)}(\bar{g}\circ\Phi_{k}^{-1})\,.
	\end{equation*}
	The inverse relation can be computed in a similar fashion yielding
	\begin{equation*}
	\bar{g} = \Phi^{*}_{k} \left[g\right] = \det\left(\boldsymbol{\mathsf{J}}\right)(g\circ\Phi_{k}) \,.
	\end{equation*}
	Given two functions $f, g \in S(\Omega _k)$ on $\Omega _k$ their inner product can be computed by
	\begin{equation*}
	\bb{ f, g}_{\Omega _k} := \int_{\Omega _k} f g\, \mathrm{d}\Omega =  \int_{\Phi_{k}\left(K\right)} fg\,\mathrm{d}\Omega = \int_K \Phi^{*}_{k} \left[f\right] \Phi^{*}_{k} \left[g\right] \,\frac{1}{\mathrm{det}\left(\boldsymbol{\mathsf{J}}\right)}\mathrm{d}K= \int_K\bar{f} \bar{g}\,\frac{1}{\mathrm{det}\left(\boldsymbol{\mathsf{J}}\right)}\mathrm{d}K\,.
	\end{equation*}
	\subsubsection{Transformation rules for dual function spaces}
	The construction of dual function spaces for an arbitrary $\Omega_k \in \mathbb{R}^3$ follows the same steps as that for the reference element.}
	
\varun{Let $\Psi ^k\bb{\bm{x}}$ for $k=0,1,2,3$ be the transformed basis functions for spaces $G\bb{\Omega}$, $C\bb{\Omega}$, $D\bb{\Omega}$ and $S\bb{\Omega}$, respectively.
Let the associated mass matrix be given by
\[ \mathbb{M}^{(k)} = \int _\Omega {\Psi ^k\bb{\bm{x}}}^\intercal {\Psi ^k\bb{\bm{x}}} \mathrm{d}\Omega \qquad \qquad \mbox{for} \quad k=0,1,2,3 \;. \]
Then the dual basis and the dual degrees of freedom are given by
\[ \widetilde{\Psi}^{d-k}\bb{\bm{x}} = \Psi ^k\bb{\bm{x}} {\mathbb{M}^{(k)}}^{-1} \qquad \mbox{and} \qquad \widetilde{\mathcal{N}}^{d-k}\bb{p^h} = \mathbb{M}^{(k)} \mathcal{N}^k\bb{p^h} \qquad \mbox{for} \quad k=0,1,2,3 \;. \]
}

\VJ{As mentioned earlier, we do not explicitly need to construct these basis functions.
	In practice we simply use their property, that the dual basis are bi-orthogonal to the primal basis.
	To post-process the dual degrees of freedom, they are first converted to primal degrees of freedom by solving a linear system of equations and then reconstructed using primal basis functions.}
\subsection{Finite element spaces for multiple elements in 3D} \label{sec:multiple_elements_3d}
In this section we introduce the finite element spaces for multiple elements in 3D case.
We follow the same procedure as for the multiple elements in 1D and 2D \varun{cases}.


\MIGnew{The g}lobal \MIGnew{basis functions} in $G\bb{\Omega} \subset H^1  \bb{\Omega}$ are in $C^0\bb{\Omega}$.
The nodal degrees of freedom at the boundaries are shared between the elements.
Since the dual representation is a linear combination of primal basis functions, $\widetilde{G}\bb{\Omega}$ also contains continuous functions.

\MIGnew{The global basis functions} in $C\bb{\Omega} \subset H\bb{\mathrm{curl;\Omega}}$ have continuous tangential components and discontinuous normal component between the neighbouring elements.
As the dual representations in $\widetilde{C}\bb{\Omega}$ are linear combinations of basis functions in $C\bb{\Omega}$, these functions also have continuous tangential component and discontinuous normal component between the neighbouring elements.

Global representation of vector fields in $D\bb{\Omega} \subset H\bb{\mathrm{div};\Omega}$ have a continuous normal component and discontinuous tangential components between the neighbouring elements.
As the dual polynomials in $\widetilde{D}\bb{\Omega}$ are linear combinations of basis functions in $D\bb{\Omega}$ they also have continuous normal component and discontinuous tangential components between the neighbouring elements.

Lastly, piecewise polynomials in $S\bb{\Omega} \subset L^2\bb{\Omega}$ are discontinuous between elements and therefore the dual representations in $\widetilde{S}\bb{\Omega}$ are also discontinuous between the elements.

\subsection{Vector operations on dual variables} \label{sec:dual_operators_3d}
In this section we will define the differential operators for the dual representations.
We will follow the procedure of 1D and 2D cases and present direct application of the ideas in 3D case.
\subsubsection{The gradient acting on $\widetilde{S} \bb{\Omega} \times \widetilde{D} \bb{\partial \Omega} $}
Let the restriction of vector fields in $D\bb{\Omega}$ to the domain boundary $\partial \Omega$ be denoted by $D\bb{\partial \Omega}$.
Using \secref{sec:2d_surface_basis} we can define the dual representation of this space.
We will denote this space as $\widetilde{D}\bb{\partial \Omega}$.
\begin{definition}
	We define gradient operation on dual space as $\widetilde{\mathrm{grad}} : \widetilde{S}\bb{\Omega} \times \widetilde{D}\bb{\partial \Omega} \mapsto \widetilde{D}\bb{\Omega}$, such that, for $\bb{s^h,\hat{s}^h} \in \widetilde{S}\bb{\Omega} \times \widetilde{D}\bb{\partial \Omega}$
	\begin{equation} \label{eq:3d_dual_gad}
	\int _{\Omega} \widetilde{\mathrm{grad}} \bb{s^h,\hat{s}^h} \bm{q}^h \d \Omega = - \int _\Omega s^h \bb{\mathrm{div}\bm{q}^h} \d \Omega + \int _{\partial \Omega} \hat{s}^h \bb{\bm{q}^h \cdot \bm{n}} \d \Gamma \qquad \forall\ \bm{q}^h \in D\bb{\Omega}\;.
	\end{equation}
\end{definition}
Since the gradient maps into $\widetilde{D}\bb{\Omega}$, we have that \MIGnew{the left hand side} is a bilinear form in $\widetilde{D}\bb{\Omega} \times D\bb{\Omega}$.
The first term on the \MIGnew{right hand side} is a bilinear form in $\widetilde{S} \bb{\Omega} \times S\bb{\Omega}$, and the boundary integral term is the bilinear form $\widetilde{D}\bb{\partial \Omega} \times D\bb{\VJ{\partial}\Omega}$.

Therefore, we can write \eqref{eq:3d_dual_gad} as
\MIGnew{ 
	\begin{equation*}
	\mathcal{N}^2 \bb{\bm{q}^h}^{\intercal} \widetilde{\mathcal{N}}^1 \bb{\widetilde{\mathrm{grad}}\bb{s^h, \hat{s}^h}} = -  \mathcal{N}^2\bb{\bm{q}^h}^{\intercal} {\mathbb{E}^{3,2}}^{\intercal} \widetilde{\mathcal{N}}^0 \bb{s^h}  +  \mathcal{N}^2 \bb{\bm{q}^h}^{\intercal} \mathbb{N}_{2} \widetilde{\mathcal{B}}^0 \bb{\hat{s}^h}  \;,
	\end{equation*}}
where \MIGnew{$\mathbb{N}_{2}$} is \MIGnew{a} sparse \MIGnew{inclusion} matrix which only contains the non-zero values $-1$ and $1$.
\varun{The subscript $'2'$ indicates that it maps boundary degrees of freedom associated to the surfaces (geometric dimension $'2'$) to the global degrees of freedom.}
Using the fact that the above equation needs to hold for all $\bm{q}^h \in D\bb{\Omega}$, we have
\begin{equation} \label{eq:3d_dual_gradient}
\widetilde{\VJ{\mathcal{N}}}^1 \bb{\widetilde{\mathrm{grad}}\bb{s^h, \hat{s}^h}} = -{ \mathbb{E}^{3,2}}^\intercal \widetilde{\mathcal{N}}^0 \bb{s^h}  +  \mathbb{N}_{\MIGnew{2}} \widetilde{\mathcal{B}}^0 \bb{\hat{s}^h}  \;. 
\end{equation}
\begin{definition} \label{def:extended_grad_3d}
	In order to construct a de Rham sequence for the dual representations, we are going to extend $\widetilde{\mathrm{grad}}$ in \eqref{eq:3d_dual_gradient} as $\widetilde{\mathrm{GRAD}}: \VJ{\widetilde{S}}(\Omega) \times \VJ{\widetilde{D}}(\partial \Omega) \rightarrow \VJ{\widetilde{D}}(\Omega) \times \VJ{\widetilde{C}}(\partial \Omega)$
	\MIGnew{
		\begin{eqnarray}
		\widetilde{\mathrm{GRAD}}\bb{s^h, \hat{s}^h} & := & \bb{ \begin{array}{cc} \dualgrad \bb{s^h,\ \hat{s}^h} , & \dualgrad _b \bb{\hat{s}^h} \end{array} } \;, \nl
		& = & \bb{\begin{array}{cc} \widetilde{\Psi}^1\bb{- {\mathbb{E}^{3,2}}^T \widetilde{\mathcal{N}}^0(s^h) + \mathbb{N}_{\MIGnew{2}} \widetilde{\mathcal{B}}^0(\hat{s}^h)}, &  \varun{-} \widetilde{\Psi}_b^1 
			\bb{{\mathbb{N}_1^\intercal \mathbb{E}^{2,1}}^T \mathbb{N}_2 \widetilde{\mathcal{B}}^0(\hat{s}^h)}
			\end{array}} \;,
		\end{eqnarray}}
\MIGnew{where $\mathbb{N}_1$ is the inclusion matrix from the boundary degrees of freedom $C(\partial \Omega)$ to the degrees of freedom of $C(\Omega)$ and $\mathbb{N}_2$ is the inclusion matrix for the degrees of freedom of $D(\partial \Omega)$ and $D(\Omega)$, respectively
\varun{and $\widetilde{\Psi}_b^1$ are the dual basis on the boundary}.}
\end{definition}
\subsubsection{The curl acting on $\widetilde{D} \bb{\Omega} \times \widetilde{C} \bb{\partial \Omega} $}
Let the restriction of vector fields in $C\bb{\Omega}$ to the domain boundary be denoted as $C\bb{\partial \Omega}$, and the corresponding dual representation be given by $\widetilde{C}\bb{\partial \Omega}$.
\begin{definition} 
	We define the curl operator on the dual space as $\widetilde{\mathrm{curl}} : \widetilde{D} \bb{\Omega} \times \widetilde{C}\bb{\partial \Omega} \mapsto \widetilde{C}\bb{\Omega}$, such that for $\bb{\bm{q}^h,\hat{\bm{q}}^h} \in \widetilde{D}\bb{\Omega} \times \widetilde{C}\bb{\partial \Omega} $ we have
	\begin{equation} \label{eq:3d_dual_curl}
	\int _{\Omega} \widetilde{\mathrm{curl}} \bb{\bm{q}^h,\hat{\bm{q}}^h} \VJ{\cdot} \bm{c}^h \d \Omega = \int _\Omega \bm{q}^h \VJ{\cdot} \bb{\mathrm{curl}\ \bm{c}^h} \d \Omega + \int _{\partial \Omega} \hat{\bm{q}}^h \VJ{\cdot} \bb{\bm{c}^h \times \bm{n}} \d \Gamma \qquad \forall \ \bm{c}^h \in C\bb{\Omega} \;.
	\end{equation}
\end{definition}
The $\widetilde{\mathrm{curl}}$ maps into $\widetilde{C}\bb{\Omega}$, therefore the \MIGnew{left hand side} term is a bilinear form in $\widetilde{C}\bb{\Omega} \times C\bb{\Omega}$.
The first term on \MIGnew{the right hand side} is a bilinear form in $\widetilde{D}\bb{\Omega} \times D\bb{\Omega}$, and the boundary integral term is a bilinear form in $\widetilde{C}\bb{\partial \Omega} \times C \bb{\partial \Omega}$.
Using assembled degrees of freedom, and appropriate expansions, we can write \eqref{eq:3d_dual_curl} as
\VJ{\begin{equation*}
	\widetilde{\mathrm{curl}}\bb{\bm{q}^h, \hat{\bm{q}}^h} = \widetilde{\Psi}^2 \bb{{\mathbb{E}^{2,1}}^\intercal \widetilde{\mathcal{N}}^1 \bb{\bm{q}^h}  +  \mathbb{N}_{\MIGnew{1}} \widetilde{\mathcal{B}}^1 \bb{\hat{\bm{q}}^h}}  \;,
	\end{equation*}}
where $\mathbb{N}_{\MIGnew{1}}$
\MIGnew{is the inclusion matrix from $C(\partial \Omega)$ to $C(\Omega)$} \varun{with dimensions $dim\bb{C(\Omega)} \times dim\bb{C(\partial\Omega)}$}.
\begin{definition} \label{def:extended_curl_3d}
	We extend the curl applied to the dual representation in \eqref{eq:3d_dual_curl} to include the boundary \varun{elements} as: $\widetilde{\mathrm{CURL}}:\VJ{\widetilde{D}}(\Omega)\times \VJ{\widetilde{C}}(\partial \Omega) \rightarrow \VJ{\widetilde{C}}(\Omega) \times \VJ{\widetilde{G}}(\partial \Omega)$
	\begin{eqnarray}
	\widetilde{\mathrm{CURL}} \bb{\bm{q}^h, \hat{\bm{q}}^h} & := & \bb{\begin{array}{cc}
		\dualcurl \bb{\bm{q}^h,\ \hat{\bm{q}}^h}, & \dualcurl _b \bb{\hat{\bm{q}}^h}
		\end{array} } \nl
	& = & \left ( \begin{array}{cc} \widetilde{\Psi}^2 \bb{{\mathbb{E}^{2,1}}^T \widetilde{\mathcal{N}}^1(\bm{q}^h) + \mathbb{N}_{\MIGnew{1}} \widetilde{\mathcal{\MIGnew{B}}}^1(\hat{q}^h)}, & \widetilde{\Psi}^2_b \bb{\mathbb{N}_0^\intercal {\mathbb{E}^{1,0}}^T \mathbb{N}_1 \widetilde{\mathcal{\MIGnew{B}}}^1(\hat{q}^h)}
	\end{array} \right ) \;,
	\end{eqnarray}
\end{definition}
\VJ{Using Definition~\ref{def:extended_grad_3d} and Definition~\ref{def:extended_curl_3d}, we have}
\begin{eqnarray}
\widetilde{\mathcal{N}}^2 \VJ{\left( \widetilde{\mbox{CURL}} ( \widetilde{\mbox{GRAD}}(s^h,\hat{s}^h) \right)} & = & \VJ{\widetilde{\mathcal{N}}^2 \left( \widetilde{\mbox{CURL}} \left( \begin{array}{cc}
	- {\mathbb{E}^{\VJ{3,2}}}^\intercal \widetilde{\mathcal{N}}^0(s^h) + \mathbb{N}_{2} \widetilde{\mathcal{B}}^0(\hat{s}^h) \;,&  \varun{-} \mathbb{N}_1^\intercal {\mathbb{E}^{2,1}}^\intercal \mathbb{N}_2 \widetilde{\mathcal{B}}^0(\hat{s}^h)
	\end{array}
	\right) \right) \nonumber} \\
& = & \VJ{\widetilde{\mathcal{N}}^2} \left ( \begin{array}{cc}
\MIGnew{-}{\mathbb{E}^{2,1}}^T {\mathbb{E}^{3,2}}^T \widetilde{\mathcal{N}}^0(s^h) + {\mathbb{E}^{2,1}}^T \mathbb{N}_{\MIGnew{2}} \widetilde{\VJ{\mathcal{B}}}^0(\hat{s}^h) \varun{-} \mathbb{N}_{\MIGnew{1}} \MIGnew{\mathbb{N}_1^\intercal }{\mathbb{E}^{2,1}}^T \MIGnew{\mathbb{N}_2} \widetilde{\VJ{\mathcal{B}}}^0(\hat{s}^h) \;, & \varun{-}\mathbb{N}_0^\intercal {\mathbb{E}^{1,0}}^\intercal \mathbb{N}_1 \mathbb{N}_1^\intercal {\mathbb{E}^{2,1}}^\intercal \mathbb{N}_2 \widetilde{\mathcal{B}}^0(\hat{s}^h) 
\end{array} \right ) \nonumber \\
& = & \left ( \begin{array}{cc}
0, & 0 \end{array} \right ) \nonumber \;.
\end{eqnarray}
\MIGnew{Where we used that $\mathbb{E}^{3,2} \mathbb{E}^{2,1} \equiv 0$ and \varun{Lemma~\ref{lem:identity_in_boundary} that}
\[
\mathbb{N}_{\MIGnew{1}} \MIGnew{\mathbb{N}_1^\intercal }{\mathbb{E}^{2,1}}^T \MIGnew{\mathbb{N}_2} \widetilde{\VJ{\mathcal{B}}}^0(\hat{s}^h) =  {\mathbb{E}^{2,1}}^T \MIGnew{\mathbb{N}_2} \widetilde{\VJ{\mathcal{B}}}^0(\hat{s}^h) \;,
\]
in which case 
\[
\mathbb{N}_0^\intercal {\mathbb{E}^{1,0}}^\intercal \mathbb{N}_1 \mathbb{N}_1^\intercal {\mathbb{E}^{2,1}}^\intercal \mathbb{N}_2 \widetilde{\mathcal{B}}^0(\hat{s}^h) = \mathbb{N}_0^\intercal {\mathbb{E}^{1,0}}^\intercal  {\mathbb{E}^{2,1}}^\intercal \mathbb{N}_2 \widetilde{\mathcal{B}}^0(\hat{s}^h) = 0 \;,
\]
since $\mathbb{E}^{2,1} \mathbb{E}^{1,0} \equiv 0$.}
\MIGnew{This} proves that \MIGnew{the dual} curl of \MIGnew{the dual} gradient of scalar field is zero.
\subsubsection{The divergence acting on $\widetilde{C} \bb{\Omega} \times \widetilde{G} \bb{\partial \Omega} $}
Let the \MIGnew{restriction of the }space $G\bb{\Omega}$ to the domain boundary, $\partial \Omega$, be denoted as $G\bb{\partial \Omega}$, and the corresponding dual space as $\widetilde{G}\bb{\partial \Omega}$.
\begin{definition}
	We define the divergence operator on dual space as $\widetilde{\mathrm{div}} : \widetilde{C}\bb{\Omega} \times \widetilde{G}\bb{\partial \Omega} \mapsto \widetilde{G}\bb{\Omega}$, such that for $\bb{\bm{c}^h, \hat{\bm{c}}^h} \in \widetilde{C}\bb{\Omega} \times \widetilde{G}\bb{\partial \Omega} $ we have
	\begin{equation} \label{eq:3d_dual_div}
	\int _{\Omega} \widetilde{\mathrm{div}} \bb{\bm{c}^h,\hat{\bm{c}}^h} p^h \d \Omega = - \int _\Omega \MIGnew{\bm{c}}^h \bb{\mathrm{grad}p^h} \d \Omega + \int _{\partial \Omega} \bb{\hat{\bm{c}}^h \cdot \bm{n}} p^h \d \Gamma \qquad \forall \ p^h \in G\bb{\Omega}\;.
	\end{equation}
\end{definition}
\MIGnew{Since} $\widetilde{\mathrm{div}}$ maps into $\widetilde{G}\bb{\Omega}$, the \MIGnew{left hand side} term is a bilinear \MIGnew{form} in $\widetilde{G}\bb{\Omega} \times G\bb{\Omega}$.
The first term on the \MIGnew{right hand side} is a bilinear \MIGnew{form} on $\widetilde{C}\bb{\Omega} \times C\bb{\Omega}$, and the boundary integral term is a bilinear \MIGnew{form} on $\widetilde{G}\bb{\partial \Omega} \times G\bb{\partial \Omega}$.

\MIGnew{U}sing assembled degrees of freedom, we can write \eqref{eq:3d_dual_div} as
\begin{equation*}
\widetilde{\mathrm{div}}\bb{\bm{c}^h, \hat{\bm{c}}^h} = \widetilde{\Psi}^3\bb{- {\mathbb{E}^{1,0}}^\intercal \widetilde{\mathcal{N}}^2 \bb{\bm{c}^h}  +  \mathbb{N}_{\MIGnew{0}} \widetilde{\mathcal{\MIGnew{B}}}^2 \bb{\hat{\bm{c}}^h}} \;,
\end{equation*}
where, $\mathbb{N}_{\MIGnew{0}}$ is a sparse \MIGnew{inclusion} matrix only containing the non-zero values $-1$ and $1$, which maps  the degrees of freedom \MIGnew{$C(\partial \Omega)$} to \MIGnew{degrees of freedom of $C(\Omega)$}.
\begin{definition} \label{def:extended_div_3d}
	We extend the divergence applied to the dual representation as $\widetilde{\mathrm{DIV}}:\VJ{\widetilde{C}}(\Omega) \times \VJ{\widetilde{G}}(\partial \Omega) \rightarrow \VJ{\widetilde{G}}(\Omega)\times 0$ as
	\begin{eqnarray}
	\widetilde{\mathrm{DIV}}\bb{\bm{c}^h \;, \hat{\bm{c}}^h} & := & \bb{\begin{array}{cc}
		\dualdiv \bb{\bm{c}^h \;, \hat{\bm{c}}^h} \;, & 0
		\end{array}} \nl
	& = & \left ( \begin{array}{cc} \widetilde{\Psi}^3 \bb{ - {\mathbb{E}^{1,0}}^T \widetilde{\mathcal{N}}^2(\VJ{\bm{c}}^h) + \mathbb{N}_{\MIGnew{0}} \widetilde{\VJ{\mathcal{B}}}^2(\hat{\VJ{\bm{c}}}^h)} , &
	0
	\end{array} \right )
	\end{eqnarray}
\end{definition}
\VJ{Using Definition~\ref{def:extended_curl_3d} and Definition~\ref{def:extended_div_3d} we have}
\begin{eqnarray}
\widetilde{\mathcal{N}}^3 \left( \widetilde{\mbox{DIV}} \left( \widetilde{\mbox{CURL}}(\bm{q}^h,\hat{q}^h) \right) \right) & = & \widetilde{\mathcal{N}}^3 \left( \widetilde{\mbox{DIV}} \left(
\begin{array}{cc}
 {\mathbb{E}^{2,1}}^T \widetilde{\mathcal{N}}^1(\bm{q}^h) + \mathbb{N}_{\MIGnew{1}} \widetilde{\VJ{\mathcal{B}}}^1(\hat{q}^h) \;, & \MIGnew{\mathbb{N}_0^\intercal} {\mathbb{E}^{1,0}}^T \MIGnew{\mathbb{N}_1} \widetilde{\VJ{\mathcal{B}}}^1(\hat{q}^h)
\end{array} \right) \right) \nonumber \\
& = & \VJ{\widetilde{\mathcal{N}}^3} \left (
\begin{array}{cc}
- {\mathbb{E}^{1,0}}^T {\mathbb{E}^{2,1}}^T \widetilde{\mathcal{N}}^1(\bm{q}^h) \VJ{-} {\mathbb{E}^{1,0}}^T \mathbb{N}_{\MIGnew{1}} \widetilde{\VJ{\mathcal{B}}}^1(\hat{q}^h) \ \varun{+} \ \mathbb{N}_{\MIGnew{0}} \MIGnew{\mathbb{N}_0^\intercal} {\mathbb{E}^{1,0}}^T \MIGnew{\mathbb{N}_1} \widetilde{\VJ{\mathcal{B}}}^1(\hat{q}^h) \;, & \VJ{0} \end{array} \right ) \nonumber \\
& = & \left ( \begin{array}{cc}
0 \;, & 0 \end{array} \right ) \nonumber \;,
\end{eqnarray}
\MIGnew{since $\mathbb{E}^{2,1}\mathbb{E}^{1,0} \equiv 0$ and \varun{using Lemma~\ref{lem:identity_in_boundary} we have that}
\[
\mathbb{N}_{\MIGnew{0}} \MIGnew{\mathbb{N}_0^\intercal} {\mathbb{E}^{1,0}}^T \MIGnew{\mathbb{N}_1} \widetilde{\VJ{\mathcal{B}}}^1(\hat{q}^h) =  {\mathbb{E}^{1,0}}^T \MIGnew{\mathbb{N}_1} \widetilde{\VJ{\mathcal{B}}}^1(\hat{q}^h) \;.
\]}
\MIGnew{This} proves that \MIGnew{the dual} divergence of \MIGnew{the dual} curl of a vector field is zero.
\VJ{\begin{definition}
		For $\widetilde{s}^h = \bb{s^h \;, \hat{s}^h} \in \widetilde{S}\bb{\Omega} \times \widetilde{D}\bb{\partial \Omega}$, we have
		\begin{equation}
		\| \widetilde{s}^h \|^2_{H \bb{\widetilde{\mathrm{grad}}; \Omega}} := 
		\| s^h \|^2_{L^2\bb{\Omega}} + \| \widetilde{\mathrm{grad}} \bb{s^h, \hat{s}^h} \|^2 _{L^2\bb{\Omega}}
		\end{equation}
		For $\widetilde{\bm{d}}^h = \bb{\bm{d}^h \;, \hat{\bm{d}}^h} \in \widetilde{D}\bb{\Omega} \times \widetilde{C}\bb{\partial \Omega}$, we have 
		\begin{equation}
		\| \widetilde{\bm{d}}^h \|^2_{H\bb{\widetilde{\mathrm{curl}} ;\Omega}} := 
		\| \bm{d}^h \|^2_{L^2\bb{\Omega}} + \| \widetilde{\mathrm{curl}} \bb{\bm{d}^h, \hat{\bm{d}}^h} \|^2 _{L^2\bb{\Omega}}
		\end{equation}
		For $\widetilde{\bm{c}}^h = \bb{\bm{c}^h \;, \hat{{c}}^h} \in \widetilde{C}\bb{\Omega} \times \widetilde{G}\bb{\partial \Omega}$, we have 
		\begin{equation}
		\| \widetilde{\bm{c}}^h \|^2_{H\bb{\widetilde{\mathrm{div}} ;\Omega}} := 
		\| \bm{c}^h \|^2_{L^2\bb{\Omega}} + \| \widetilde{\mathrm{div}} \bb{\bm{c}^h, \hat{{c}}^h} \|^2 _{L^2\bb{\Omega}}
		\end{equation}
\end{definition}}
\VJ{\begin{definition}
		We can define the norm in the dual of the trace spaces using \cite[eq. 1.16]{1986Girault} as: \\
		For $\hat{s}^h \in \widetilde{D}\bb{\partial \Omega}$, and $\widetilde{s}^h = \bb{s^h \;, \hat{s}^h} \in \dualSb $ we have
		\begin{equation}
		\| \hat{s}^h \|_{\widetilde{D}(\partial \Omega)} := \inf_{\varun{s^h \in \widetilde{S}(\Omega) }} \| \widetilde{s}^h \|_{H\bb{\widetilde{\mathrm{grad}}; \Omega}} \;.
		\end{equation}
		For $\hat{\bm{d}}^h \in \widetilde{C}\bb{\partial \Omega}$, and $\widetilde{\bm{d}}^h \in \bb{\bm{d}^h \;, \hat{d}^h} \in \dualDb $ we have
		\begin{equation}
		\| \hat{\bm{d}}^h \|_{\widetilde{C}(\partial \Omega)} := \inf_{\varun{\bm{d}^h \in \widetilde{D}(\Omega)}} \| \widetilde{\bm{d}}^h \|_{H\bb{\widetilde{\mathrm{curl}};\Omega}} \;.
		\end{equation}
		For $\hat{{c}}^h \in \widetilde{G}\bb{\partial \Omega}$, and $\widetilde{\bm{c}}^h = \bb{\bm{c}^h \;, \hat{c}^h} \in \widetilde{C}\bb{\Omega} \times \widetilde{G}\bb{\partial \Omega}$ we have
		\begin{equation}
		\| \hat{{c}}^h \|_{\widetilde{G}(\partial \Omega)} := \inf_{\varun{\bm{c}^h \in \widetilde{C}(\Omega)}} \| \widetilde{\bm{c}}^h \|_{H\bb{\widetilde{\mathrm{div}};\Omega}} \;.
		\end{equation}
\end{definition}}
\subsection{Discrete de Rham sequence}
Using the function spaces defined in \secref{sec:multiple_elements_3d} and \secref{sec:dual_operators_3d}, we can define the primal de Rham's sequence as
\begin{equation*}
\begin{tikzcd}[row sep=1.0cm, column sep = 2cm]
G\bb{\Omega} \arrow{r}{\mathrm{grad}} & C\bb{\Omega} \arrow{r}{\mathrm{curl}} & D\bb{\Omega} \arrow{r}{\mathrm{div}} & S\bb{\Omega}
\end{tikzcd} \;.
\end{equation*}
and using the function spaces defined in \secref{sec:dual_operators_3d} we can define the dual de Rham sequence as
\begin{equation*}
\begin{tikzcd}[row sep=1.0cm, column sep = 2cm]
\widetilde{S}\bb{\Omega} \times \widetilde{D}\bb{\partial \Omega} \arrow{r}{\widetilde{\mathrm{GRAD}}} & \widetilde{D}\bb{\Omega} \times \widetilde{C}\bb{\partial \Omega} \arrow{r}{\widetilde{\mathrm{CURL}}} & \widetilde{C}\bb{\Omega} \times \widetilde{G}\bb{\partial \Omega} \arrow{r}{\widetilde{\mathrm{DIV}}} & \widetilde{G}\bb{\Omega} \times \VJ{0}
\end{tikzcd} \;.
\end{equation*}
\section{Mixed formulation of the Poisson equation}\label{sec:mixed_formulation}
In this section we want to assess the use of dual \varun{spaces} for the mixed formulation of the Poisson problem.
We present \varun{an} application of {3D finite element} spaces to a constrained minimization problem of the Poisson equation.
We will compare the results from the two formulations: 1) with primal spaces only, and 2) with primal and dual spaces.
In this application it will be shown that the use of dual spaces can give much sparser systems with a lower condition number without compromising on the accuracy of the solution.
Let ${\Omega} \subset \mathbb{R}^3$, then for $\phi \in L^2({\Omega})$ and $\bm{q} \in H(\mbox{div};{\Omega})$ we define the functional
\begin{equation*}
\mathcal{L}(\phi,\bm{q};f, {\hat{\phi}}) := \int_{{\Omega}} \frac{1}{2} |\bm{q} |^2 \, \mathrm{d}{\Omega} + \int_{{\Omega}} \phi \left ( \mbox{div}\,\bm{q}-f \right ) \, \mathrm{d}{\Omega}\ {- \int_{\partial {\Omega}} \hat{\phi}\ (\bm{q} \cdot \bm{n})\,  \,\mathrm{d}\Gamma} \;,
\end{equation*}
for prescribed functions $f \in L^2({\Omega})$ and $\hat{\phi} \in H^{1/2}\bb{\partial {\Omega}}$. The optimality conditions for this functional are given by
\begin{equation} \label{eq:var_mixed}
\left \{ \begin{array}{lll}
(\bm{p}, \bm{q} )_{{\Omega}} + ( \mbox{div}\,\bm{p}, \phi )_{{\Omega}} &= \int_{\partial {\Omega}} (\bm{p} \cdot \bm{n})\, \hat{\phi} \,\mathrm{d}\Gamma \quad \quad & \forall\ \bm{p} \in H(\mbox{div};{\Omega}) \nl
( \varphi, \mbox{div} \,\bm{q})_{{\Omega}} &= (\varphi,f)_{{\Omega}} & \forall\ \varphi \in L^2({\Omega})
\end{array} \right. \;.
\end{equation}
This corresponds to a Poisson equation for $\phi$ with Dirichlet boundary condition $\phi = \hat{\phi}$ along the boundary.
We will consider two different discretizations for this problem.
For the first approximation we choose $(\bm{q}^h,\phi^h) \in D({\Omega})\times S({\Omega})$, we will call this primal-primal formulation, while in the second case we approximate the solution as $(\bm{q}^h,\phi^h) \in D({\Omega})\times \widetilde{S}({\Omega})$, we will call this primal-dual formulation.
\subsection{Primal-primal formulation}
Let $\bm{q}^h \in D\bb{\Omega}$ be represented as 
\begin{equation*} 
\bm{q}^h({x,y,z}) = \Psi ^{2} ({x,y,z}) \mathcal{N}^{2}(\bm{q}^h) \;.
\end{equation*}
Then, using \eqref{eq:incidence_21_3D}, the divergence is given by
\[ \mbox{div} \, \bm{q}^h ({x,y,z}) = \Psi^{3}({x,y,z}) \mathbb{E}^{3,2} \mathcal{N}^{2}(\bm{q}^h) \;.\]
If we use this in the variational formulation \eqref{eq:var_mixed}, we get
\begin{equation}
\left ( \begin{array}{cc}
\mathbb{M}^{(2)} & {\mathbb{E}^{3,2}}^\intercal \mathbb{M}^{(3)} \nl
\mathbb{M}^{(3)} \mathbb{E}^{3,2} & 0
\end{array} \right )
\left ( \begin{array}{c}
\mathcal{N}^{2}(\bm{q}^h) \nl
\mathcal{N}^{3}(\phi^h)
\end{array} \right ) =
\left ( \begin{array}{c}
\varun{\mathbb{N}_2 \widetilde{\mathcal{B}}^0\bb{\hat{\phi}^h} } \nl
\mathbb{M}^{\bb{3}}\mathcal{N}^{3}(f^h)
\end{array} \right ) \;,
\label{eq:primal_primal_system}
\end{equation}
where the degrees of freedom of $f$ and the boundary integral term are evaluated \varun{using \varun{Definition~\ref{def:nodal_sammpling_L2}}}.

The incidence matrix $\mathbb{E}^{3,2}$ is a sparse topological matrix.
All metric properties are contained in the mass matrices ${\mathbb{M}}^{(2)}$ and ${\mathbb{M}}^{(3)}$.
For high order methods, these matrices are dense matrices with large full blocks that destroy the sparsity of the incidence matrix with which they are multiplied in (\ref{eq:primal_primal_system}).
We refer to this formulation as the {\em primal-primal} formulation, because both $\bm{q}^h$ and $\phi^h$ are expanded in primal basis functions.
If the mesh is deformed, all sub-matrices in \eqref{eq:primal_primal_system} will change because the {mass} matrices ${\mathbb{M}}^{\bb{k}}$ will change and need to be recomputed.
\subsection{Primal-dual formulation}
Alternatively, we may approximate $\phi^h \in \widetilde{S}(K)$.
In this case the discrete system will be
\begin{equation}
\left ( \begin{array}{cc}
{\mathbb{M}}^{(2)} & {{\mathbb{E}}^{3,2}}^\intercal \nl
{\mathbb{E}}^{3,2} & 0
\end{array} \right )
\left ( \begin{array}{c}
\mathcal{N}^{2}(\bm{q}^h) \nl
\widetilde{\mathcal{N}}^0(\phi^h)
\end{array} \right ) =
\left ( \begin{array}{c}
\varun{\mathbb{N}_2 \widetilde{\mathcal{B}}^0\bb{\hat{\phi}^h} } \nl
\mathcal{N}^{3}(f^h)
\end{array} \right ) \;,
\label{eq:primal_dual_system}
\end{equation}
{where the degrees of freedom of $f^h$ and the prescribed boundary conditions $\hat{\phi}^h$ are the same as evaluated earlier in \eqref{eq:primal_primal_system}.}
We see that if we expand $\phi^h$ in terms of dual polynomials, the discrete divergence and gradient blocks in \eqref{eq:primal_dual_system} are sparse and no longer depend on the metric of the mesh geometry.
It is only the mass matrix ${\mathbb{M}}^{(2)}$ that needs to be constructed for each element separately, unless all elements have the same size, shape \MIGnew{and polynomial degree}.
\begin{remark}
	We can immediately convert \eqref{eq:primal_primal_system} to \eqref{eq:primal_dual_system}.
	The mass matrix ${\mathbb{M}}^{(3)}$ in the second row of \eqref{eq:primal_primal_system} can be cancelled on both sides of the equation, while the mass matrix ${\mathbb{M}}^{(3)}$ in the first row can be contracted with the degrees of freedom $\mathcal{N}^{3}(\phi^h)$ to give ${\mathbb{M}}^{(3)} \mathcal{N}^{3}(\phi^h)$, but these new unknowns are just the dual degrees of freedom $\widetilde{\mathcal{N}}^0(\phi^h)$ according to \eqref{eq:mass_matrix_M3_3d}.
\end{remark}

\subsection{The discrete inf-sup condition}
At the continuous level one establishes well-posedness by showing that the divergence operator, $\mathrm{div}$, from $H\bb{\mathrm{div}; \Omega}$ into $L^2(\Omega)$ is surjective in which case the Closed Range Theorem states that the Hilbert adjoint of the divergence operator is bounding and therefore injective, \cite{2010Boffi}. At the continuous level surjectivity of the divergence operator is proven through the auxiliary problem: For an arbitrary $p \in L^2(\Omega)$ find $\psi \in H_0^1(\Omega)$ such that
\begin{equation*}
\int_{\Omega} \mbox{grad}\ \psi \cdot \mbox{grad} \ \phi \ \mathrm{d}\Omega = - \int_{\Omega} p \phi \ \mathrm{d}\Omega \;,\quad \quad \forall \phi \in H_0^1(\Omega) \;.
\end{equation*}
The Lax-Milgram lemma ensures uniqueness of $\psi$. If we set $\bm{u}_p = \mbox{grad}\ \psi$ then we have $\mbox{div}\ \bm{u}_p = p$ and since $p \in L^2(\Omega)$ was arbitrary, this shows surjectivity of the divergence operator.

At the finite dimensional level the discrete inf-sup condition is derived in exactly the same way. Let $\MIGnew{p}^h \in \widetilde{S}({\Omega})$.
Now prove that there is a $\VJ{\bm{q}}^h \in D({\Omega})$ which is mapped by the divergence operator onto $\MIGnew{p}^h$.
Just as in the continuous setting, we use an auxiliary problem: Find $\VJ{\psi^h \in \widetilde{S}({\Omega})}$ such that
\begin{equation*}
\int_{{\Omega}} \VJ{\widetilde{\mbox{grad}}\ \bb{\psi^h,0} \cdot \widetilde{\mbox{grad}} \ \bb{\phi^h,0}} \ \mathrm{d}{\Omega} = - \int_{{\Omega}} \MIGnew{p}^h \phi^h \ \mathrm{d}{\Omega} \;,\quad \quad \forall \phi^h \in \VJ{\widetilde{S}}({\Omega}) \;.
\end{equation*}
Using the dual basis functions this translates into
\[
\widetilde{\mathcal{N}}^0(\phi^h)^{\intercal} \mathbb{E}^{{3,2}} {\mathbb{M}^{({2})}}^{-1} {\mathbb{E}^{{3,2}}}^{\intercal} \widetilde{\mathcal{N}}^0(\psi^h) = - \widetilde{\mathcal{N}}^0(\phi^h)^{\intercal} {\mathbb{M}^{({3})}}^{-1} \widetilde{\mathcal{N}}^0(p^h) \;,
\]
which has to hold for all vectors $\widetilde{\mathcal{N}}^0(\phi^h)$ and therefore we have
\begin{equation}
\mathbb{E}^{{3,2}} {\mathbb{M}^{({2})}}^{-1} {\mathbb{E}^{{3,2}}}^{\intercal} \widetilde{\mathcal{N}}^0(\psi^h) = -  {\mathbb{M}^{({3})}}^{-1} \widetilde{\mathcal{N}}^0(p^h) \;.
\label{eq:discrete_auxiliary}
\end{equation}
Using \eqref{eq:mass_matrix_M2_3d} and \eqref{eq:3d_dual_gradient} we define $\mathcal{N}^{{2}}({\bm{u}}^h) = - {\mathbb{M}^{({2})}}^{-1} {\mathbb{E}^{{3,2}}}^{\intercal} \widetilde{\mathcal{N}}^0(\psi^h) \in D(\VJ{\Omega})$ in which case (\ref{eq:discrete_auxiliary}) can be written as
\[
\mathbb{E}^{{3,2}} \mathcal{N}^{{2}}({\bm{u}}^h) =   {\mathbb{M}^{({3})}}^{-1} \widetilde{\mathcal{N}}^0(p^h) \;,
\]
or
\begin{equation*}
\mathbb{M}^{({3})} \mathbb{E}^{{3,2}} \mathcal{N}^{{2}}({\bm{u}}^h) =   \widetilde{\mathcal{N}}^0(p^h) \;.
\end{equation*}
This equation states that for all vectors $\widetilde{\mathcal{N}}^0(p^h)$ there exists a vector $\mathcal{N}^{2}(\VJ{\bm{u}}^h)$ which is mapped by the discrete divergence, $\mathbb{E}^{{3,2}}$ to degrees of freedom in $S({\Omega})$ and then mapped by the mass matrix $\mathbb{M}^{({3})}$ onto $\VJ{\widetilde{S}}({\Omega})$.
\subsection{Test case}
\begin{figure}
	\centering
	\includegraphics[scale=0.5]{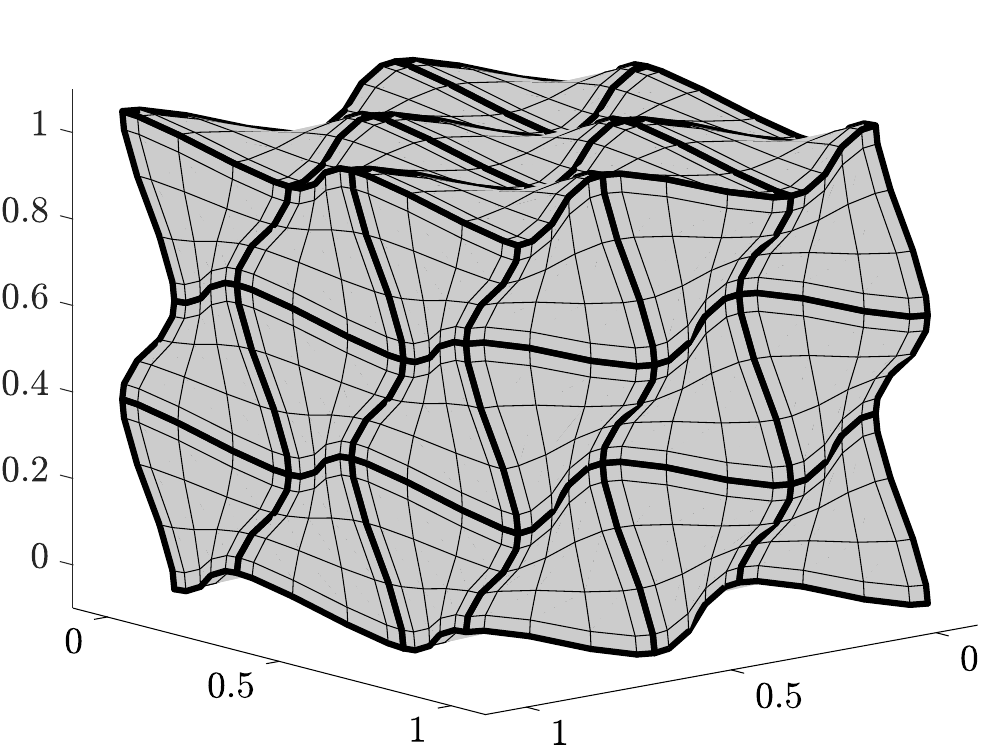}
	\caption{Test domain for number of elements $\VJ{K_{el}}=3 \times 3 \times 3$ and polynomial order {$N=5$}. \VJ{The bold lines show the element domains and the thin lines show the mesh generated by GLL quadrature points.}}
	\label{fig:domain_3D}
\end{figure}
In the following test case we compare the primal-primal formulation with the primal-dual formulation.
Using the primal-dual formulation we obtain sparse algebraic formulation and lower condition number without compromising on quantative results.
Consider the domain (taken from~\cite{2012Wheeler}) shown in \figref{fig:domain_3D}.
The deformed mesh coordinates $(x,y,z) \in \Omega$ are obtained by transforming the {reference} coordinates $(\xi,\eta, \zeta) \in {[-1,1]^3}$ with the mapping
\begin{equation*} \label{eq:mesh}
\left\lbrace
\begin{array}{l}
x = \hat{x} + 0.03 \cos \bb{3 \pi \hat{x}} \cos \bb{3 \pi \hat{y}} \cos \bb{3 \pi \hat{z}}  \\[1.5ex]
y = \hat{y} {- 0.04} \cos \bb{3 \pi \hat{x}} \cos \bb{3 \pi \hat{y}} \cos \bb{3 \pi \hat{z}}  \\[1.5ex]
z = \hat{z} {+ 0.05} \cos \bb{3 \pi \hat{x}} \cos \bb{3 \pi \hat{y}} \cos \bb{3 \pi \hat{z}}
\end{array} \right. \;, \mbox{where} \quad
\left\lbrace
\begin{array}{l}
\hat{x} = 0.5 \bb{1 + \xi}  \\[1.5ex]
\hat{y} = 0.5 \bb{1 + \eta}  \\[1.5ex]
\hat{z} = 0.5 \bb{1 + \zeta}
\end{array} \right. \;.
\end{equation*}
We compare both the formulations with a manufactured solution $\phi _{ex} = \sin(2\pi x)\sin(2 \pi y)\sin(2 \pi z)$ which gives
\[ f_{ex} = -\mbox{div} (\mbox{grad}\, \phi _{ex}) \;, \]
and {we use} Dirichlet boundary conditions over entire domain given by $\hat{\phi} = \phi _{ex} |_{\partial {\Omega}}$.
\begin{figure}[!htb]
	\centering
	\begin{subfigure}[b]{0.4\textwidth}
		\includegraphics[width=0.92\textwidth]{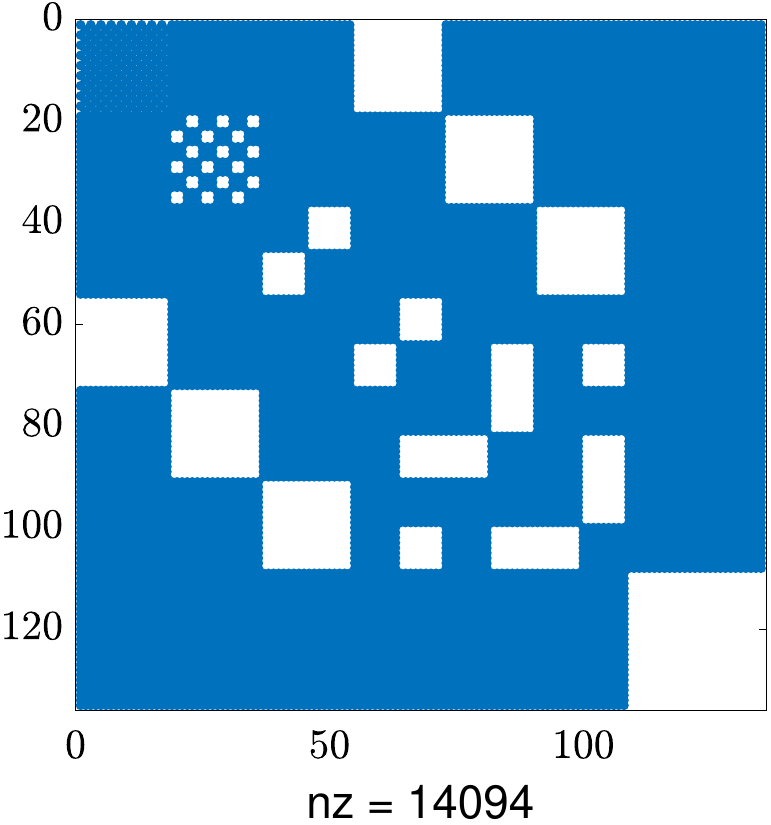}
	\end{subfigure}%
	\quad
	~ 
	\qquad
	\begin{subfigure}[b]{0.4\textwidth}
		\includegraphics[width=0.92\textwidth]{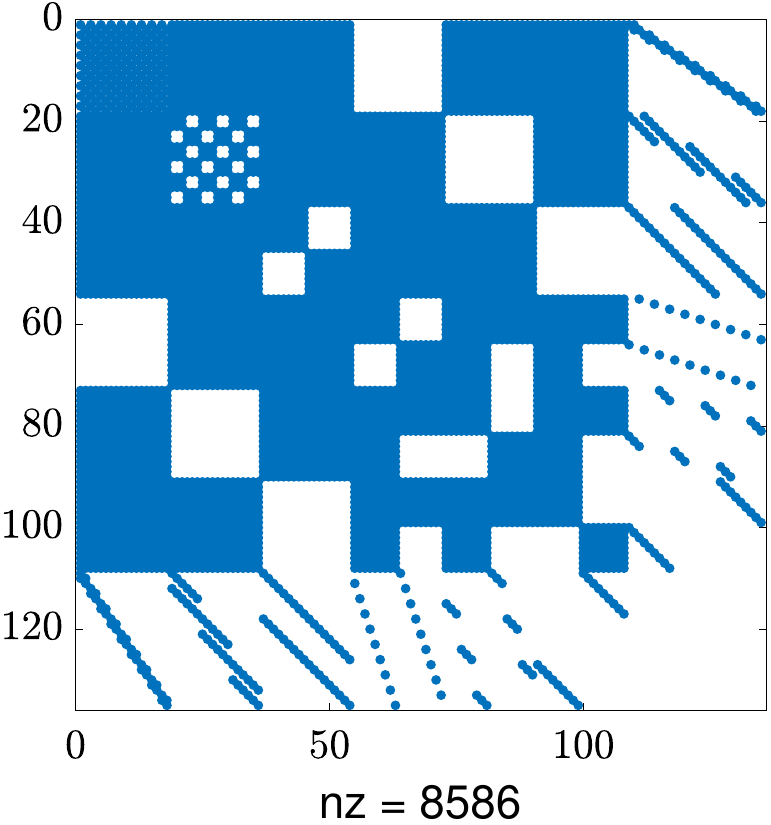}
	\end{subfigure}
	\begin{subfigure}[b]{0.4\textwidth}
		\includegraphics[width=0.92\textwidth]{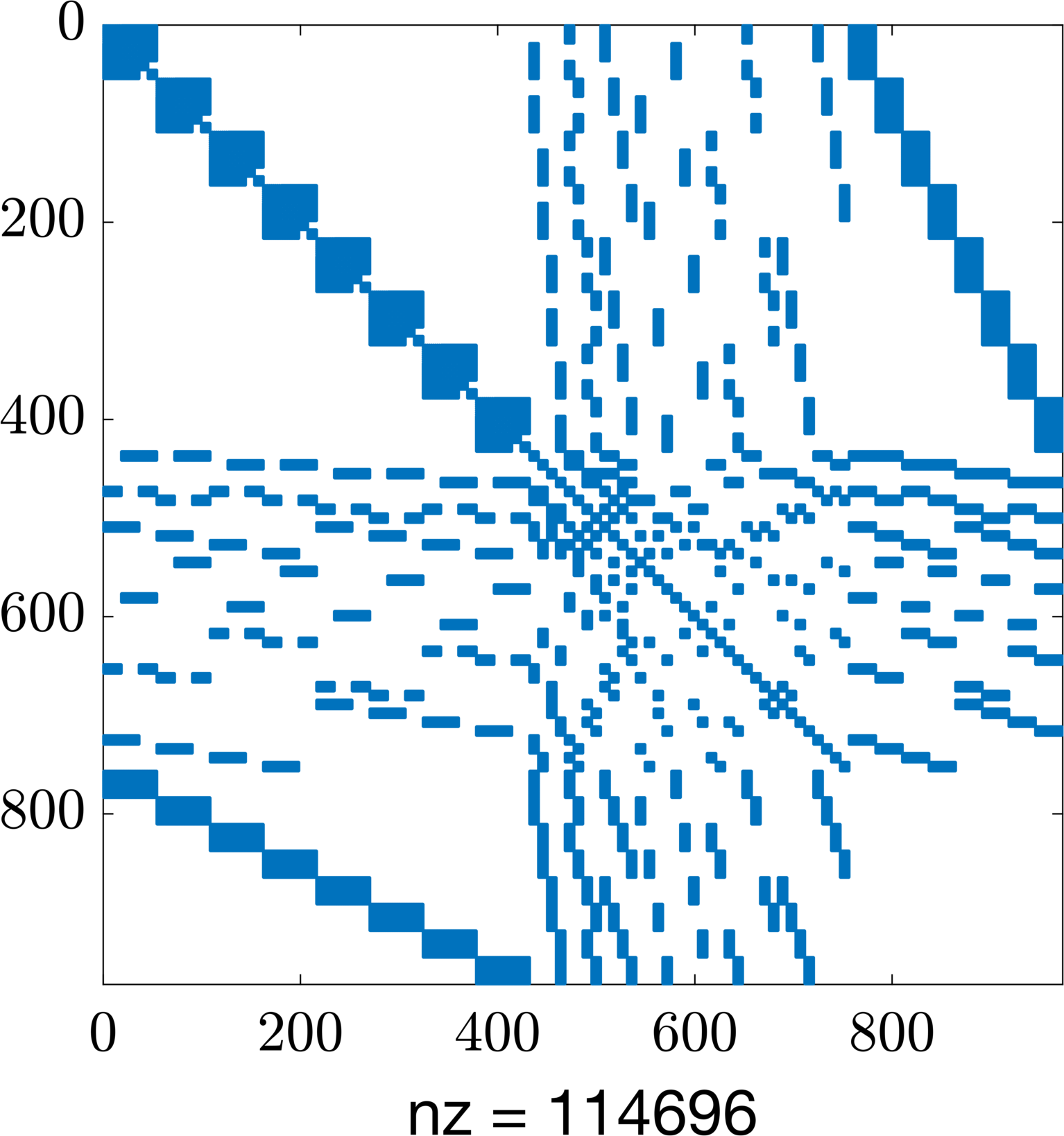}
	\end{subfigure}%
	\quad
	~ 
	\qquad
	\begin{subfigure}[b]{0.4\textwidth}
		\includegraphics[width=0.92\textwidth]{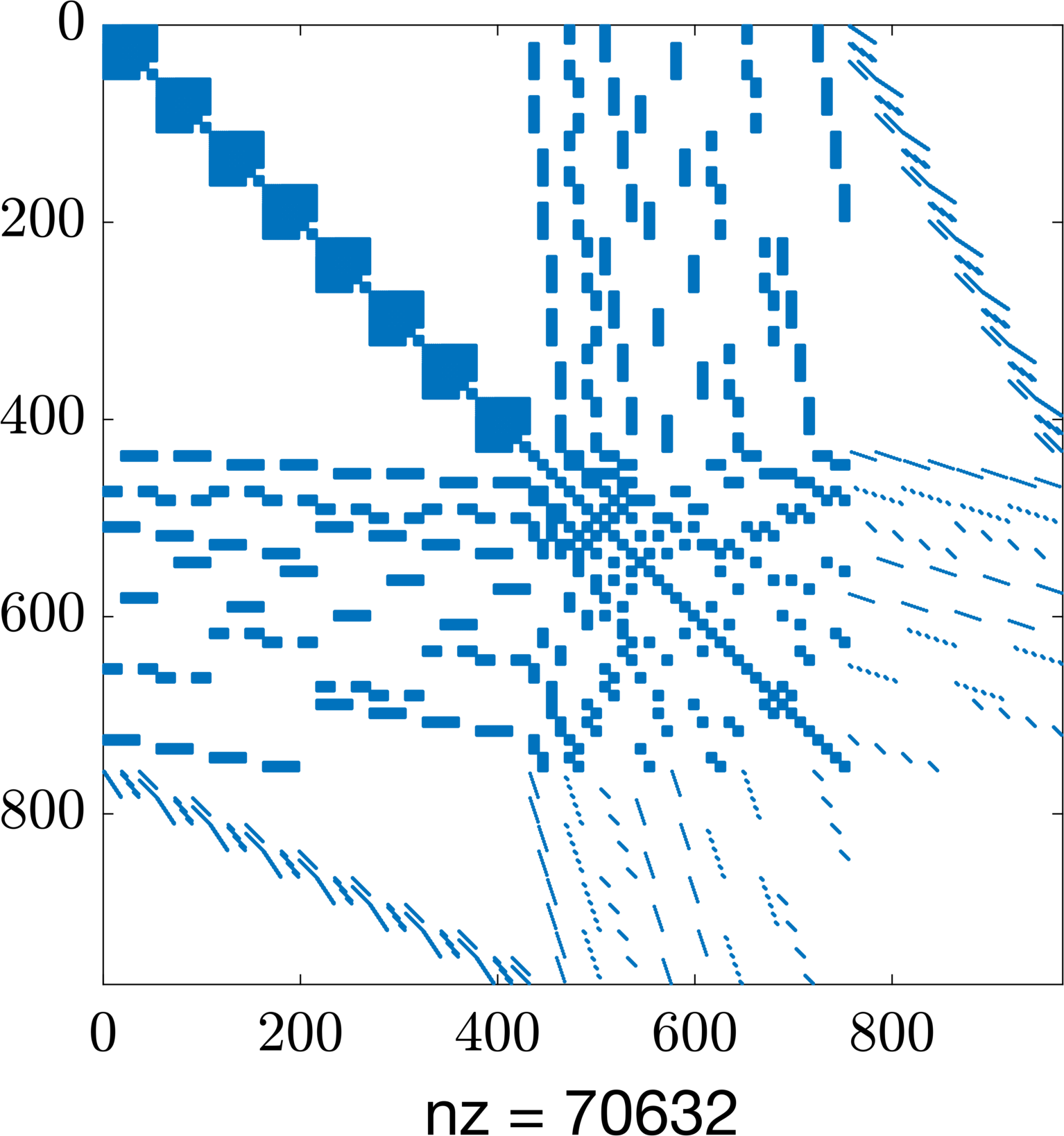}
	\end{subfigure}%
	\caption{Sparsity plots for \VJ{$N=3$,} i) \emph{Top-left}: $K_{el} = 1 \times 1 \times 1$, primal-primal formulation, ii) \emph{Top-right}: $K_{el} = 1 \times 1 \times 1$, primal-dual formulation, iii) \emph{Bottom-left}: $K_{el} = 2 \times 2 \times 2$, primal-primal formulation, iv) \emph{Bottom-right}: $K_{el} = 2 \times 2 \times 2$, primal-dual formulation.} \label{fig:spy_c_00}
\end{figure}

In Figure~\ref{fig:spy_c_00} we compare the sparsity pattern of primal-primal and primal-dual formulation for the algebraic formulations \eqref{eq:primal_primal_system} and \eqref{eq:primal_dual_system} for the given test case.
In the top-left and the top-right figure we show the sparsity plot for a single element with $N=3$.
The $8586$ non-zero entries of the primal-dual formulation are much less than the $14094$ non-zero entries of the primal-primal formulation.
In the bottom-left and bottom-right we show similar comparison but now for the multiple element case with {total} number of elements $K_{el}= 2 \times 2 \times 2$ and \VJ{$N=3$}.
The non-zero entities in primal-dual formulation $70632$ are much less than those in the primal-primal formulation ${114696}$.
%
\begin{table}[!htb]
	\caption{Condition number for the primal-primal formulation and the primal-dual formulation for {$K_{el}=1 \times 1 \times 1$ polynomial degree $N=2,4,8$}.}
	\centering
	\begin{tabular}{ccc}
		\hline
		$N$   & Primal-Primal & Primal-Dual \\
		\hline
		{2} & {362.2070} & 33.7474    \\
		{4} & {7.5959e+3} & 218.9917  \\
		{8} & {3.1730e+5} & 6.0411e+3 \\
		\hline
	\end{tabular}
	\label{tab:cond}
\end{table}

In Table~\ref{tab:cond} we list the condition number of the algebraic system of the two formulations for {$K_{el} = 1 \times 1 \times 1$ with polynomial degree $N=2,4,8$}.
We observe that the condition number for the primal-primal formulation is significantly higher than that of the primal-dual formulation.
{Also the rate of growth of condition number \varun{for increasing $N$} is higher for the primal-primal formulation}.
In this sense the use of dual polynomials can also be interpreted as a form of inverse type mixed preconditioning~\cite{2005Chen}
\begin{equation*}
\left(
\begin{array}{ll}
{\mathbb{M}}^{\bb{2}} &  {{\mathbb{E}}^{3,2}}^T {\mathbb{M}}^{\bb{3}} \nl
{\mathbb{M}}^{\bb{3}}{\mathbb{E}}^{3,2} & 0
\end{array}
\right)
=
\bb{
	\begin{array}{ll}
	\mathbb{I} & 0 \nl
	0 & {\mathbb{M}}^{\bb{3}}
	\end{array}
}
\bb{
	\begin{array}{ll}
	{\mathbb{M}}^{\bb{2}} & {{\mathbb{E}}^{3,2}}^T \nl
	{\mathbb{E}}^{3,2} & 0
	\end{array}
}
\bb{
	\begin{array}{ll}
	\mathbb{I} & 0 \nl
	0 & {\mathbb{M}}^{\bb{3}}
	\end{array}
} \;.
\end{equation*}

\begin{figure}[!htb]
	\centering
	\begin{subfigure}[b]{0.45\textwidth}
		\includegraphics[width=\textwidth]{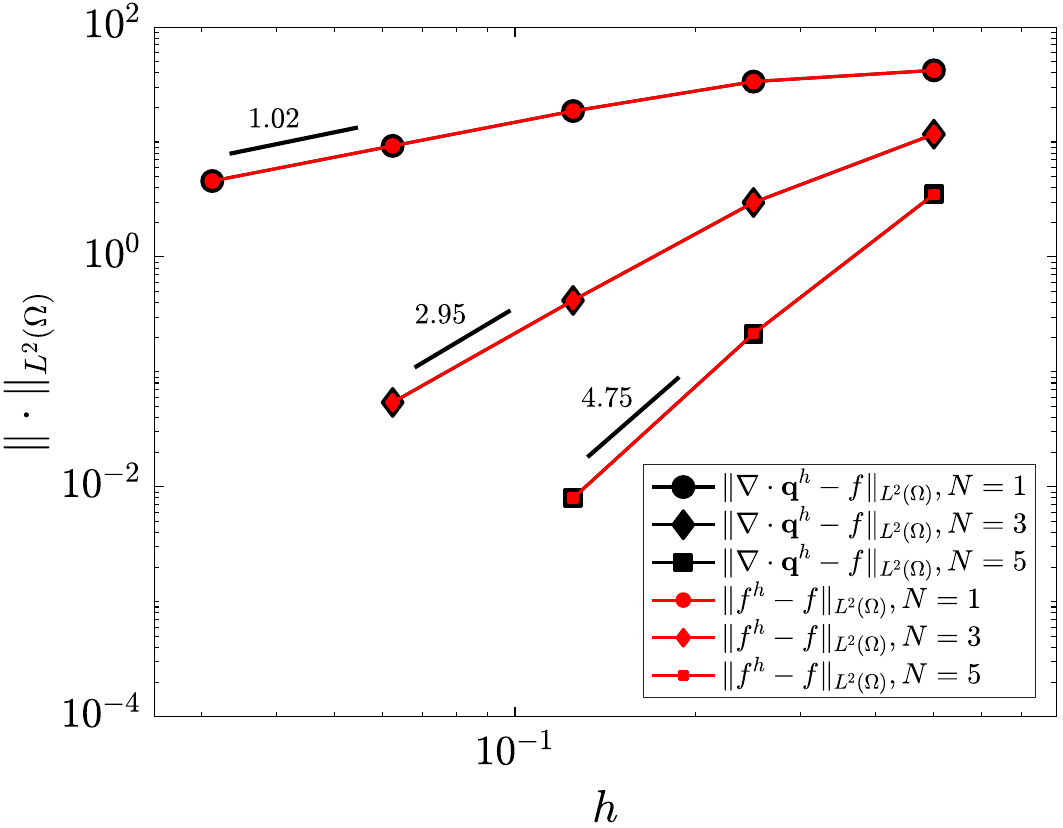}
	\end{subfigure}%
	\quad
	~ 
	~ 
	\begin{subfigure}[b]{0.45\textwidth}
		\includegraphics[width=\textwidth]{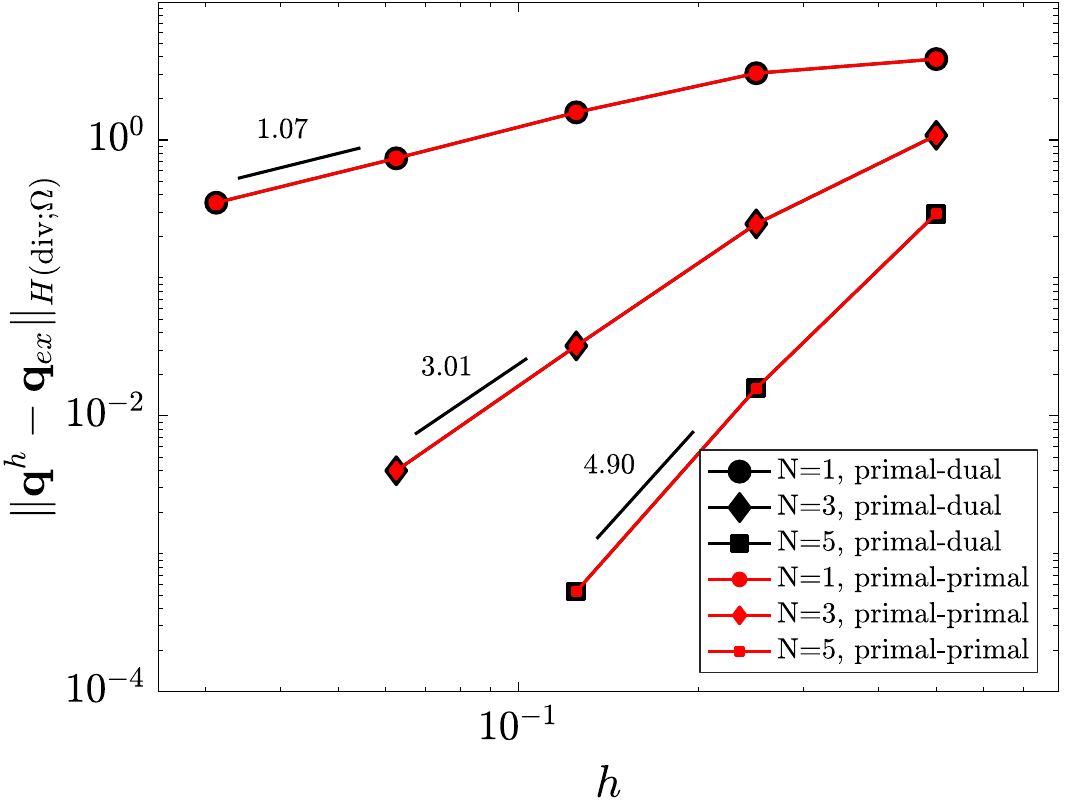}
	\end{subfigure}
	\quad
	\begin{subfigure}[b]{0.45\textwidth}
		\includegraphics[width=\textwidth]{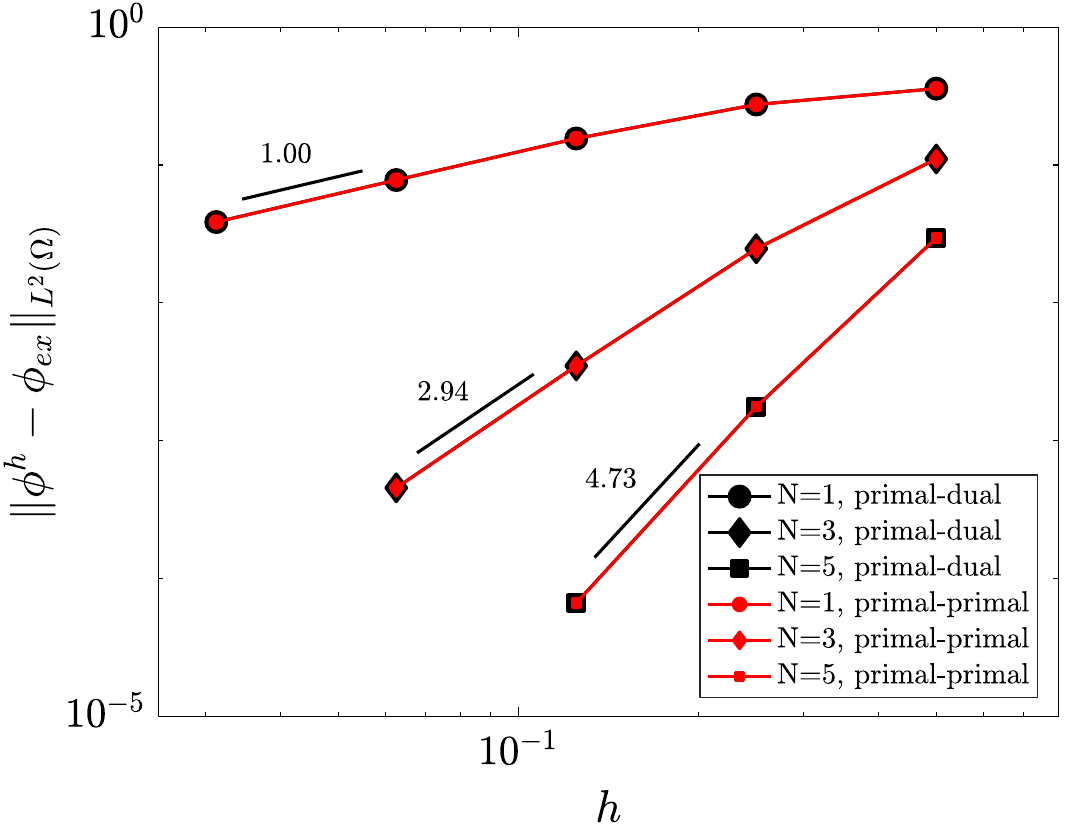}
	\end{subfigure}
	\caption{Top left: $L^2$-error in constraint $(\mbox{div}\, \bm{q}^h -f)$, and the interpolation error in RHS term $\bb{f^h -f}$. Top right: $H\bb{\mathrm{div}}$-error in flux $\bm{q}^h$. Bottom: $L^2$-error in $\phi ^h$.}
	\label{fig:results}
\end{figure}

In the top-left plot of \figref{fig:results} we show the $L^2$-error in the constraint $(\mbox{div}\, \bm{q}^h - f)$ and the interpolation error in the RHS term $\bb{f^h -f}$, for primal-dual formulation.
On the $x$-axis we have the element length \VJ{(assuming no deformation, i.e. $c=0$)} $h = 2/\sqrt[3]{K_{el}}$ and on the $y$-axis we have the $L^2$-error.
We see that the error in \varun{the} constraint is equal to the interpolation error of the \MIGnew{right hand side} term.
The error converges optimally with $\mathcal{O}\bb{{N}}$ because $\phi ^h$ is a piecewise polynomial of degree $N-1$.


In the bottom-left plot of \figref{fig:results} we see the {error convergence} of the fluxes in $H\bb{\mbox{div}}$-norm.
The results from the primal-primal formulation and the primal-dual formulation coincide with each other and the error converges optimally with $\mathcal{O}({N})$.

In the plot at the bottom right of \figref{fig:results} we see the convergence in $L^2$-error of $\phi ^h$.
The results from the primal-primal formulation and the primal-dual formulation overlap.
We see optimal \varun{rate of} convergence of errors, of $\mathcal{O}({N})$.

In terms of accuracy \figref{fig:results} shows that the results from primal-primal formulation and primal-dual formulation are {equal} up to machine precision, {and the rate of convergence are optimal for both the formulations.}

For the next test problem we choose a pair of div-grad problems that are dual to each other in a continuous setting.
We will show that this duality continues to hold true at \MIGnew{the} discrete level by using the primal-dual representations for these problems.

\section{The Dirichlet-Neumann problems}\label{sec:Dirichlet_Neumann}
{The mathematical theory of finite elements often makes use of the equivalence of dual problems. In general this equivalence no longer holds at the discrete level. In this section we want to show that this equivalence continues to hold at the finite dimensional level when dual representations are employed.}
In the proof of Lemma~2.2 {in} \cite{CarstensenDemkowiczGopalakrishnan}{, for instance,} use is made of an equivalence between a Dirichlet and a Neumann {problem}.
We start with the two problems: Given $\hat{\phi} \in H^{1/2}(\partial {\Omega})$
\begin{enumerate}
	\item The Dirichlet problem:
	Find $\phi \in H^1({\Omega})$ such that
	\begin{equation}
	\left \{ \begin{array}{ll}
	\phi = \hat{\phi} & \mbox{on } \partial {\Omega} \\[1ex]
	-\mbox{div}\, (\mbox{grad}\, \phi ) + \phi = 0 \quad \quad & \mbox{in } {\Omega}
	\end{array} \right . \;.
	\label{eq:Dirichlet}
	\end{equation}
	\item The Neumann problem:
	Find $\bm{q} \in H(\mathrm{div};{\Omega})$ such that
	\begin{equation}
	\left \{ \begin{array}{ll}
	\mbox{div}\,\bm{q} = \hat{\phi} & \mbox{on } \partial {\Omega} \\[1ex]
	-\mbox{grad}\, (\mbox{div}\, \bm{q} ) + \bm{q} = 0 \quad \quad & \mbox{in } {\Omega}
	\end{array} \right . \;.
	\label{eq:Neumann}
	\end{equation}
\end{enumerate}
If $\phi$ solves the Dirichlet problem \eqref{eq:Dirichlet}, then $\bm{q}$ solves the Neumann problem \eqref{eq:Neumann}, if and only if $\phi = \mbox{div}\, \bm{q}$. Furthermore, it follows that~\cite{CarstensenDemkowiczGopalakrishnan}
\[ \| \hat{\phi}\|_{H^{1/2}(\partial {\Omega})} = \| \phi\|_{H^1({\Omega})} = \| \bm{q} \|_{H(\mathrm{div};{\Omega})} \;.\]

The finite dimensional problem is to find suitable finite dimensional function spaces, $D\VJ{\bb{\Omega}} \subset H(\mathrm{div};{\Omega})$ and \VJ{ $\widetilde{S} \bb{\Omega} \subset L^2\bb{\Omega}$} for \VJ{$\hat{\phi}^h \in \widetilde{D}\bb{\partial \Omega}$}
\begin{enumerate}
	\item The Dirichlet problem:
	Find $\phi^h \in {\widetilde{S}} \bb{{\Omega}}$ such that 
	\begin{equation}
	\left \{ \begin{array}{ll}
	\phi^h = \hat{\phi}^{{h}}\ & \mbox{on } \partial {\Omega} \\[1ex]
	-\mbox{div}\, (\mbox{grad}\, \phi^h ) + \phi^h = 0 \quad \quad & \mbox{in } {\Omega}
	\end{array} \right . \;.
	\label{eq:discrete_Dirichlet}
	\end{equation}
	\item The Neumann problem:
	Find $\bm{q}^h \in D \bb{{\Omega}}$ such that
	\begin{equation}
	\left \{ \begin{array}{ll}
	\mbox{div}\,\bm{q}^h = \hat{\phi}^{{h}}\ & \mbox{on } \partial {\Omega} \\[1ex]
	-\mbox{grad}\, (\mbox{div}\, \bm{q}^h ) + \bm{q}^h = 0 \quad \quad & \mbox{in } {\Omega}
	\end{array} \right . \;,
	\label{eq:discrete_Neumann}
	\end{equation}
\end{enumerate}
such that the solutions $\phi^h$ and $\bm{q}^h$ will satisfy $\phi^h = \mbox{div}\,\bm{q}^h$ identically in the \VJ{domain} ${\Omega}$. Furthermore, we wish to prove that in this case
\begin{equation} \label{eq:norm_equality}
\VJ{\| \hat{\phi}^h \|_{\widetilde{D}\bb{\partial \Omega}}^2 =\| \widetilde{\phi}^h\|_{\VJ{H\bb{\widetilde{\mathrm{grad}}; \Omega}}} = \| \bm{q}^h \|_{H\bb{\mathrm{div};\Omega}}} \;.
\end{equation}
\subsection{The Neumann problem}
Consider ${\Omega} \subset \mathbb{R}^d$, with $d=2$.
Then the variational formulation of the Neumann problem, \eqref{eq:discrete_Neumann}, is given by: For $\hat{\phi}^h \in \VJ{\widetilde{D}}(\partial \Omega)$ find $\bm{q}^h \in D \bb{{\Omega}}$ such that
\begin{equation}\label{eq:variational_formulation_Neumann}
\left ( \mbox{div}\, \bm{p}^h, \mbox{div}\,\bm{q}^{{h}} \right )_{{\Omega}} + \left ( \bm{p}^h,\bm{q}^h \right )_{{\Omega}} = \int_{\partial {\Omega}} ( {\bm{p}}^h \cdot \bm{n}) \hat{\phi}^{h} \, \mathrm{d}\Gamma \;, \quad \quad \forall {\bm{p}}^h \in D \bb{{\Omega}} \;.
\end{equation}
\noindent
We \VJ{expand} $\bm{q}^h$ as in \eqref{eq:q_expansion_2D} \VJ{with transformation rules from \eqref{eq:u_with_mapping} as}
\[ \bm{q}^h({x,y}) = \Psi ^{{1}} ({x,y}) \mathcal{N}^{1}(\bm{q}^h) \;.
\]
Then using \eqref{eq:incidence_21} for divergence in the variational formulation \eqref{eq:variational_formulation_Neumann} we obtain
\begin{equation}
\left ( \mbox{div}\, {\bm{p}}^h, \mbox{div}\,\bm{q}^{{h}} \right )_{{\Omega}} + \left ( {\bm{p}}^h,\bm{q}^h \right )_{{\Omega}} =
\mathcal{N}^{1}({\bm{p}}^h)^\intercal {\mathbb{E}^{2,1}}^\intercal  \mathbb{M}^{(2)} \mathbb{E}^{2,1} \mathcal{N}^{1}(\bm{q}^h) + \mathcal{N}^{1}({\bm{p}}^h)^\intercal \mathbb{M}^{(1)} \mathcal{N}^{1}(\bm{q}^h) \;, \quad \quad {\forall \bm{p}^h \in D \bb{{\Omega}}} \;.
\label{eq:hdiv_discrete}
\end{equation}
The boundary terms on the right hand side of \eqref{eq:variational_formulation_Neumann} are evaluated \VJ{using Definition~\ref{def:nodal_sammpling_L2} as}
\VJ{\[ \widetilde{\mathcal{B}}^0 \bb{\hat{\phi}^h} = \int _{\partial \Omega} {\Psi^1_b\bb{x}}^\intercal \hat{\phi}\bb{x} dx \;. \]}
Using the fact that \eqref{eq:variational_formulation_Neumann} should hold for all $\mathcal{N}^1\bb{\VJ{\bm{p}}^h}$ gives
\begin{equation}
{\mathbb{E}^{2,1}}^T  \mathbb{M}^{(2)} \mathbb{E}^{2,1} \mathcal{N}^{1}(\bm{q}^h) + \mathbb{M}^{(1)} \mathcal{N}^{1}(\bm{q}^h) = \bdof \;.
\label{eq:Neumann_problem_primal}
\end{equation}
\subsection{The Dirichlet problem}
Consider now the Dirichlet problem given by \eqref{eq:discrete_Dirichlet} on the domain ${\Omega} \subset \mathbb{R}^d$, with $d=2$. The variational formulation for this problem is given by: For $\hat{\phi}^h \in \VJ{\widetilde{D}}(\partial \Omega)$ find $\phi^h \in {\VJ{\widetilde{S}}} \bb{\Omega}$, such that
\begin{equation}
\left ( \VJ{\widetilde{\mbox{grad}}}\,\VJ{\bb{\varphi^h,0}},{\VJ{\widetilde{\mbox{grad}}}}\,\VJ{\bb{\phi ^h , \hat{\phi}^h}} \right )_{{\Omega}} + \left ( {\varphi}^h,\phi ^h \right )_{{\Omega}} = \int_{\partial {\Omega}} \hat{\phi}^h\ \frac{\partial {\varphi}^h}{\partial \bm{n}} \, \mathrm{d}\Gamma \qquad \qquad \forall \quad \varphi ^h \in \widetilde{S}\bb{\Omega} \;.
\label{eq:variational_formulation_Dirichlet}
\end{equation}
\VJ{We} know that the gradient of $\phi^h$ is given by
\begin{equation}
{\widetilde{\mbox{grad}}}\, \VJ{\bb{\phi ^h , \hat{\phi}^h}}(x,y) = \widetilde{\Psi}^1({x,y}) \left ( - {\mathbb{E}^{2,1}}^\intercal \widetilde{\VJ{\mathcal{N}}}^0(\phi^h) + \bdof   \right ) \;.
\label{eq:gradient_of_phi}
\end{equation}
The gradient of the test functions ${\varphi}^h$ is discretized similarly, but then the variations on the boundary are set to zero, therefore
\[{ \widetilde{\mbox{grad}}}\, \VJ{\bb{\varphi ^h , 0}}({x,y}) = \widetilde{\Psi}^1({x,y}) \VJ{\bb{{- \mathbb{E}^{2,1}}^\intercal \widetilde{\VJ{\mathcal{N}}}^0(\varphi^h)}}  \;.\]
If we use this in the variational formulation \eqref{eq:variational_formulation_Dirichlet} we have \VJ{that the variations with the second argument on the RHS of \eqref{eq:gradient_of_phi} becomes zero, and from the remaining terms we have that}
\[
\left ({ \widetilde{\mbox{grad}}}\,\VJ{\bb{\varphi ^h , 0}},{ \widetilde{\mbox{grad}}}\,\VJ{\bb{\phi ^h , \hat{\phi}^h}} \right )_{{\Omega}}+ \left ( \varphi^h,\phi^h \right )_{{\Omega}} =
\widetilde{\mathcal{N}}^0(\varphi ^h)^\intercal \mathbb{E}^{2,1} \widetilde{\mathbb{M}}^{(1)} \left [  {\mathbb{E}^{2,1}}^\intercal \widetilde{\VJ{\mathcal{N}}}^0(\phi^h) - \bdof \right ] + \widetilde{\mathcal{N}}^0(\varphi^h)^\intercal \widetilde{\mathbb{M}}^{(0)} \widetilde{\VJ{\mathcal{N}}}^0(\phi^h) =0\;.
\]
Once again using the fact that equality should hold for all $\widetilde{\mathcal{N}}^0(\varphi^h)$ the discrete formulation is given by
\begin{equation}
\mathbb{E}^{2,1} \widetilde{\mathbb{M}}^{(1)} {\mathbb{E}^{2,1}}^\intercal \widetilde{\VJ{\mathcal{N}}}^0(\phi^h) + \widetilde{\mathbb{M}}^{(0)} \widetilde{\VJ{\mathcal{N}}}^0(\phi^h) = \mathbb{E}^{2,1} \widetilde{\mathbb{M}}^{(1)} \bdof \;.
\label{eq:Dirichlet_problem_dual}
\end{equation}

\subsection{Relation between Dirichlet and Neumann problem} \label{sec:equality_proof}
What we need to check now is that the solutions of \eqref{eq:Neumann_problem_primal} and \eqref{eq:Dirichlet_problem_dual} are related by $\phi^h = \mbox{div}\, \bm{q}^h$. This discrete relation translates into
\begin{equation}
\widetilde{\mathcal{N}}^0(\phi^h) = \mathbb{M}^{(2)} \mathbb{E}^{2,1} \mathcal{N}^{1}(\bm{q}^h) \;.
\label{eq:phi_q_relation}
\end{equation}
In order to establish this relation, we fill \eqref{eq:phi_q_relation} in \eqref{eq:Dirichlet_problem_dual} to obtain
\begin{equation}
\mathbb{E}^{2,1} \widetilde{\mathbb{M}}^{(1)} {\mathbb{E}^{2,1}}^\intercal \mathbb{M}^{(2)} \mathbb{E}^{2,1} \mathcal{N}^{1}(\bm{q}^h) + \widetilde{\mathbb{M}}^{(0)} \mathbb{M}^{(2)} \mathbb{E}^{2,1} \mathcal{N}^{1}(\bm{q}^h) = \mathbb{E}^{2,1} \widetilde{\mathbb{M}}^{(1)} \bdof \;.
\label{eq:proof_step_1}
\end{equation}
We substitute \eqref{eq:Neumann_problem_primal} in \eqref{eq:proof_step_1}  to get
\begin{equation*}
- \mathbb{E}^{2,1} \widetilde{\mathbb{M}}^{(1)} \mathbb{M}^{(1)} \mathcal{N}^{1}(\bm{q}^h) + \mathbb{E}^{2,1} \widetilde{\mathbb{M}}^{(1)} \bdof + \widetilde{\mathbb{M}}^{(0)} \mathbb{M}^{(2)} \mathbb{E}^{2,1} \mathcal{N}^{1}(\bm{q}^h) = \mathbb{E}^{2,1} \widetilde{\mathbb{M}}^{(1)} \bdof \;.
\label{eq:proof_step_2}
\end{equation*}
Then we use the fact that $\widetilde{\mathbb{M}}^{(1)} \mathbb{M}^{(1)} = \mathbb{I}$ and $\widetilde{\mathbb{M}}^{(0)} \mathbb{M}^{(2)} = \mathbb{I}$ to get
\begin{equation*}
-\mathbb{E}^{2,1} \mathcal{N}^{1}(\bm{q}^h) + \mathbb{E}^{2,1} \widetilde{\mathbb{M}}^{(1)} \bdof + \mathbb{E}^{2,1} \mathcal{N}^{1}(\bm{q}^h) = \mathbb{E}^{2,1} \widetilde{\mathbb{M}}^{(1)} \bdof \quad \;,
\end{equation*}
which proves the relation between the Dirichlet and the Neumann problem.

It remains to show that $\VJ{\| \widetilde{\phi}^h\|_{H\bb{\widetilde{\mathrm{grad}};\Omega}}}= \| \bm{q}^h \|_{\varun{H\bb{\mathrm{div};\Omega}}}$
Using \eqref{eq:gradient_of_phi} we have that
\begin{equation}
\VJ{\| \widetilde{\phi}^h\|^2_{H\bb{\widetilde{\mathrm{grad}};\Omega}}} = \widetilde{\mathcal{N}}^0(\phi^h)^\intercal \widetilde{\mathbb{M}}^{(0)}  \widetilde{\mathcal{N}}^0(\phi^h) + \left [  \VJ{\widetilde{\mathcal{B}}^0 \bb{\hat{\phi}^h}^\intercal \mathbb{N}_1^\intercal} - \widetilde{\mathcal{N}}^0(\phi^h)^T \mathbb{E}^{2,1} \right ] \widetilde{\mathbb{M}}^{(1)} \left [\bdof - {\mathbb{E}^{2,1}}^T \widetilde{\mathcal{N}}^0(\phi^h) \right ] \;.
\label{eq:H1-norm_phi_h}
\end{equation}
Since we have just established that, $\widetilde{\mathcal{N}}^0(\phi^h) = \mathbb{M}^{(2)}\mathbb{E}^{2,1} \mathcal{N}^1(\bm{q}^h)$, we can insert this in \eqref{eq:H1-norm_phi_h}
\begin{eqnarray}
\VJ{\| \widetilde{\phi}^h\|^2_{H\bb{\widetilde{\mathrm{grad}};\Omega}}} & = & \mathcal{N}^1(\bm{q}^h)^T {\mathbb{E}^{2,1}}^\intercal \mathbb{M}^{(2)} \widetilde{\mathbb{M}}^{(0)}  \mathbb{M}^{(2)}\mathbb{E}^{2,1} \mathcal{N}^1(\bm{q}^h) \nonumber \\
& + & \left [ \VJ{\widetilde{\mathcal{B}}^0 \bb{\hat{\phi}^h}^\intercal \mathbb{N}_1^\intercal} - \mathcal{N}^1(\bm{q}^h)^\intercal {\mathbb{E}^{2,1}}^\intercal \mathbb{M}^{(2)} \mathbb{E}^{2,1} \right ] \widetilde{\mathbb{M}}^{(1)} \left [ \bdof - {\mathbb{E}^{2,1}}^\intercal \mathbb{M}^{(2)}\mathbb{E}^{2,1} \mathcal{N}^1(\bm{q}^h) \right ] \nonumber \\
& \stackrel{\eqref{eq:Neumann_problem_primal}}{=} & \mathcal{N}^1(\bm{q}^h)^\intercal {\mathbb{E}^{2,1}}^\intercal \mathbb{M}^{(2)}\mathbb{E}^{2,1} \mathcal{N}^1(\bm{q}^h)
+  \mathcal{N}^1(\bm{q}^h)^\intercal \mathbb{M}^{(1)} \widetilde{\mathbb{M}}^{(1)} \mathbb{M}^{(1)}  \mathcal{N}^1(\bm{q}^h) \nonumber \\
& = & \mathcal{N}^1(\bm{q}^h)^\intercal  \mathbb{M}^{(1)}  \mathcal{N}^1(\bm{q}^h) + \mathcal{N}^1(\bm{q}^h)^\intercal {\mathbb{E}^{2,1}}^\intercal \mathbb{M}^{(2)}\mathbb{E}^{2,1} \mathcal{N}^1(\bm{q}^h) \nonumber \\
& = & \| \bm{q}^h \|_{\varun{H\bb{\mathrm{div};\Omega}}}^2 \;,
\label{eq:equal_norms}
\end{eqnarray}
where we used again that $\widetilde{\mathbb{M}}^{(0)}  \mathbb{M}^{(2)} = \mathbb{I}$ and $\widetilde{\mathbb{M}}^{(1)} \mathbb{M}^{(1)} = \mathbb{I}$ and the fact that the degrees of freedom of $\bm{q}^h$ satisfy~\eqref{eq:Neumann_problem_primal}.

\VJ{Finally, using \eqref{eq:grad_norm_2d}, for $\widetilde{\phi}^h = \bb{\phi^h \;, \hat{\phi}^h} \in \dualSb$ we have that
	\begin{equation}
	\| \hat{\phi}^h \|_{\widetilde{D}(\partial \Omega)} := \inf_{\varun{\phi^h \in \widetilde{S}(\Omega)}} \| \widetilde{\phi}^h \|_{H\bb{\widetilde{\mathrm{grad}}; \Omega}} \;.
	\end{equation}
	By taking variations, \varun{for fixed $\hat{\phi}^h$}, we see that the infimum is attained for the function $\phi^h$ which satisfies
	\begin{equation}
	\left ( \widetilde{\mbox{grad}}(\phi^h,\hat{\phi}^h), \widetilde{\mbox{grad}}(\varphi^h,0) \right )_{\Omega} + \left ( \phi^h, \varphi^h \right )_{\Omega} = 0 \;\quad \quad \forall \varphi^h \in \widetilde{S}(\Omega) \;,
	\end{equation}
	which is just the Dirichlet problem (\ref{eq:Dirichlet_problem_dual}).
Therefore, for the Dirichlet problem we have by definition that $\| \varun{\hat{\phi}^h \|_{\widetilde{D}(\partial \Omega)} = \| \widetilde{\phi}^h \|_{H\bb{\widetilde{\mathrm{grad}};\Omega}} }$.}
\subsection{Test case}
In this section we solve the Dirichlet \eqref{eq:variational_formulation_Dirichlet} and the Neumann \eqref{eq:variational_formulation_Neumann} problems on ${\Omega} \in [0,1]^2$ with one spectral element
for a non-trivial boundary condition $\hat{\phi}$ given by
\begin{equation*}
\hat{\phi} = \left \{ \begin{array}{ll}
0 \quad \quad & \mbox{for } x=0 \mbox{ and } y=0 \\[1ex]
- \sin(\pi y) \quad \quad & \mbox{for } x=1 \\[1ex]
-\ln \left (1-3x(1-x)\right ) \quad \quad & \mbox{for } y=1
\end{array} \right . \;.
\label{eq:boundary_condition}
\end{equation*}

\begin{figure}[!htb]
	\centering
	\begin{subfigure}[b]{0.4\textwidth} \label{fig:dom_a}
		\includegraphics[width=\textwidth]{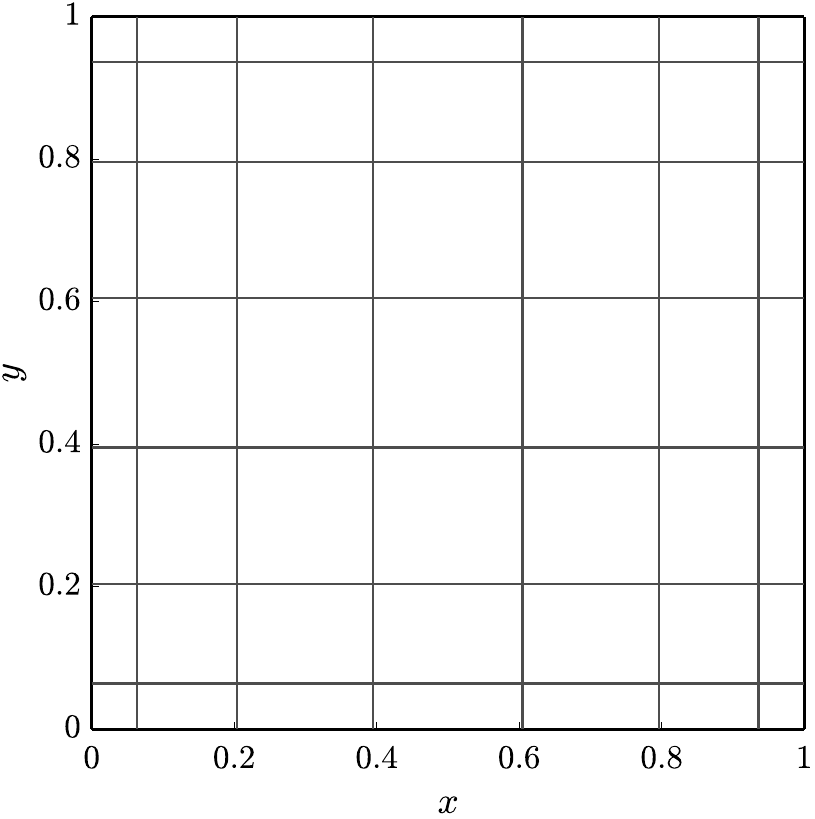}
		\caption{Orthogonal mesh, $c=0.0$}
	\end{subfigure}%
	\quad
	~ 
	\qquad
	\begin{subfigure}[b]{0.4\textwidth}
		\includegraphics[width=\textwidth]{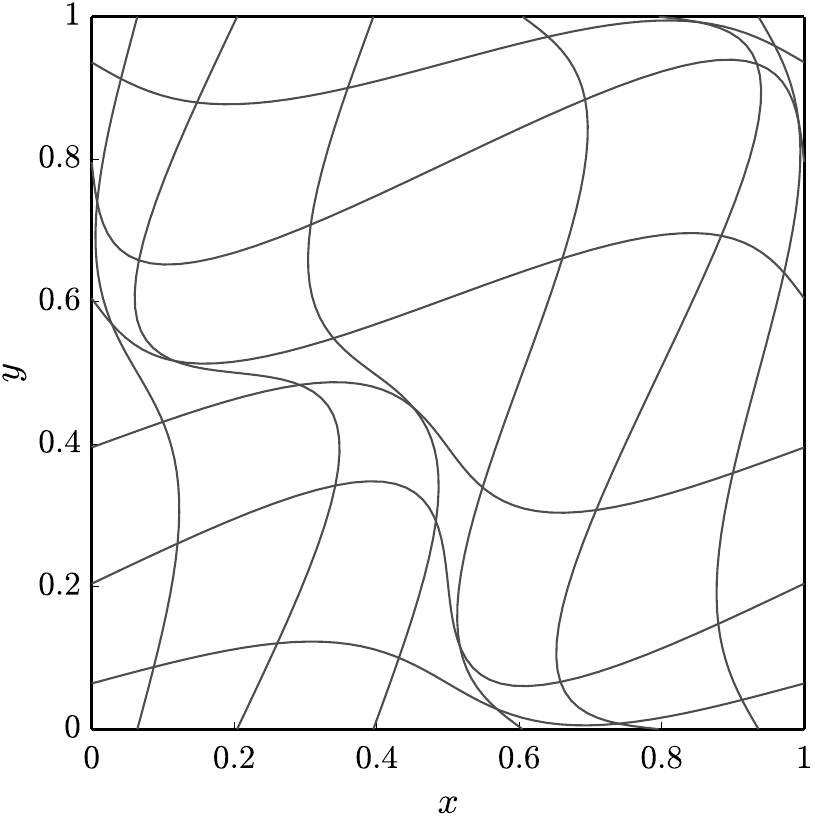}
		\caption{Curved mesh, $c=0.3$}
	\end{subfigure}
	\caption{Meshes generated by the transformation \VJ{of quadrature points according to} \eqref{eq:mapping_curved} for, $N={7}$, $K_{el} = 1$, $c=0.0$ and $c=0.3$.} \label{fig:domain}
\end{figure}

For this test case we use a `standard' orthogonal spectral element shown on the left, and a deformed spectral element shown on the right, of \figref{fig:domain}.

The deformed mesh coordinates $(x,y)$ are obtained by transforming the {reference} coordinates $(\xi,\eta) \VJ{\in [-1,1]^2}$ with the mapping
\begin{equation} \label{eq:mapping_curved}
\left\lbrace
\begin{array}{l}
x = \frac{1}{2}+ \frac{1}{2}(\xi + c \sin(\pi \xi)\sin(\pi \eta)) \\[1.5ex]
y = \frac{1}{2}+ \frac{1}{2}(\eta + c \sin(\pi \xi)\sin(\pi \eta))
\end{array} \right. \;,
\end{equation}
where $c$ is the deformation coefficient.
\begin{figure}[!htb]
	\centering
	\begin{subfigure}[b]{0.45\textwidth}
		\includegraphics[width=\textwidth]{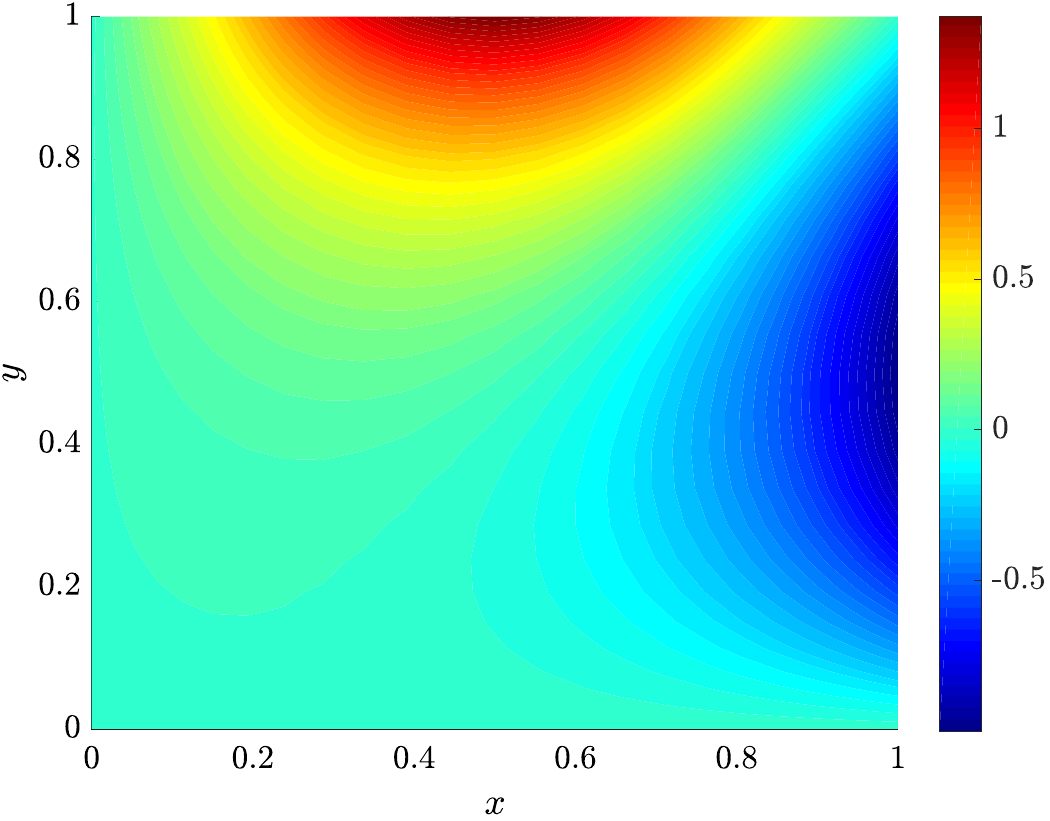}
		\caption{Solution $\mbox{div}\,\bm{q}^h$ of \eqref{eq:discrete_Neumann}}
		\label{fig:divq_c0}
	\end{subfigure}%
	\quad
	~ 
	\begin{subfigure}[b]{0.45\textwidth}
		\includegraphics[width=\textwidth]{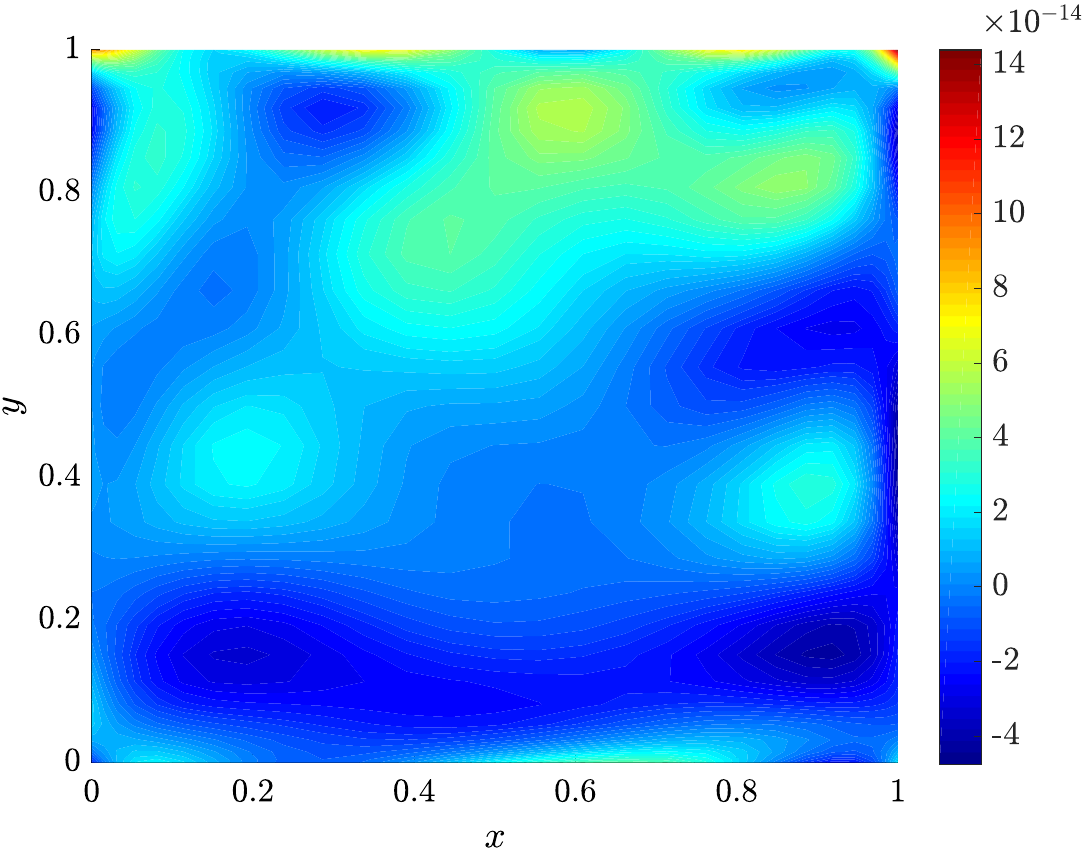}
		\caption{Difference between $\mbox{div}\,\bm{q}^h$ and the solution $\phi^h$ from \eqref{eq:discrete_Dirichlet}}
		\label{fig:divq_phi_c0}
	\end{subfigure}
	\caption{Comparison between $\phi^h$ obtained from \eqref{eq:discrete_Dirichlet} and $\mbox{div}\,\bm{q}^h$ calculated using \eqref{eq:discrete_Neumann} for $N=8$, $K_{el} = 1$, on an orthogonal mesh with $c=0$.} \label{fig:Dual_problems_on_orthog_mesh}
\end{figure}

In \figref{fig:divq_c0} the numerical solution $\mbox{div}\,\bm{q}^h$ is shown on the orthogonal mesh, $c=0$ for $N=8$.
The solution $\phi^h$ on the same mesh is visually indistinguishable from $\mathrm{div} \bm{q}^h$, therefore in \figref{fig:divq_phi_c0} the difference between $\mbox{div}\,\bm{q}^h$ and $\phi^h$ is shown.
The difference between both solutions is of {order $10^{-14}$}.

\begin{figure}[!htb]
	\centering
	\begin{subfigure}[b]{0.45\textwidth}
		\includegraphics[width=\textwidth]{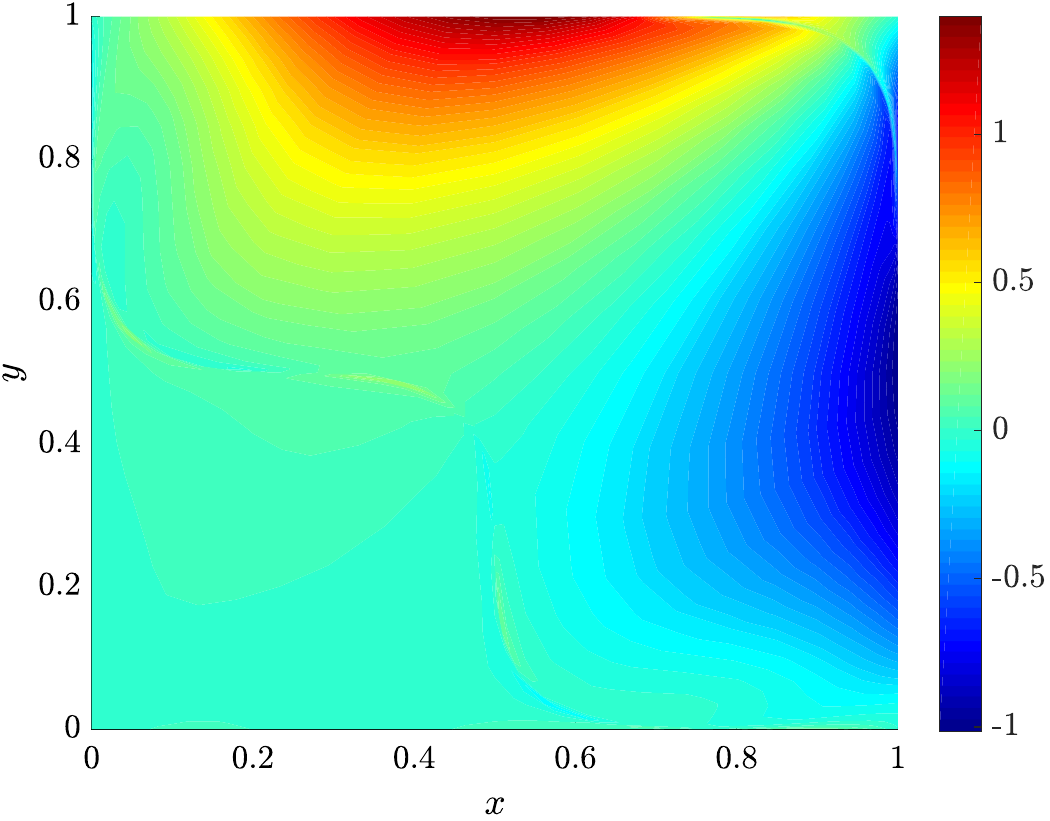}
		\caption{Solution $\mbox{div}\,\bm{q}^h$ of \eqref{eq:discrete_Neumann}}
		\label{fig:divq_c03}
	\end{subfigure}%
	\quad
	~ 
	\begin{subfigure}[b]{0.45\textwidth}
		\includegraphics[width=\textwidth]{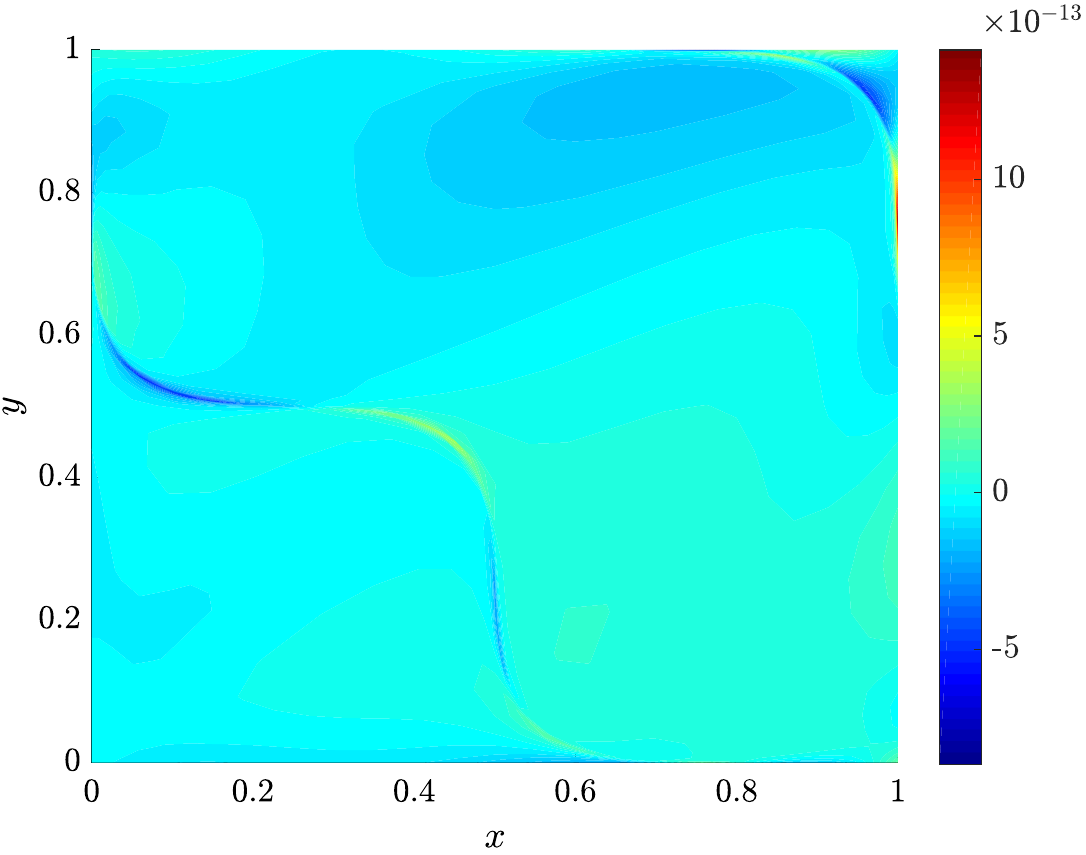}
		\caption{Difference between $\mbox{div}\,\bm{q}^h$ and the solution $\phi^h$ from \eqref{eq:discrete_Dirichlet}}
		\label{fig:divq_phi_c03}
	\end{subfigure}
	\caption{Comparison between $\phi^h$ obtained from \eqref{eq:discrete_Dirichlet} and $\mbox{div}\,\bm{q}^h$ calculated using \eqref{eq:discrete_Neumann} for $N=8$, $K_{el} = 1$, on curvilinear mesh with $c=0$.3.} \label{fig:dual_problems_on_curved_mesh}
\end{figure}

In \figref{fig:divq_c03} $\mbox{div}\,\bm{q}^h$ is plotted for the deformed grid with $c=0.3$, for $N=8$.
On the deformed mesh we expect the solution to be less accurate than on the orthogonal mesh, but $\phi^h$ computed on the same mesh is graphically still identical to \figref{fig:divq_c03}.
The difference between $\mbox{div}\,\bm{q}^h$ and $\phi^h$ is shown in \figref{fig:divq_phi_c03}.
This confirms that for this test case the discrete equivalence \eqref{eq:phi_q_relation} holds.

\begin{table}[htbp]
	\centering
	\caption{Norms $\VJ{\| \widetilde{\phi}^h\|_{H\bb{\widetilde{\mathrm{grad}};\Omega}}}$ and $\|\bm{q}^h\|_{\varun{H\bb{\mathrm{div};\Omega}}}$ on three different meshes as a function of the polynomial degree $N$.}
	\begin{tabular}{lcccccc}
		\hline
		& \multicolumn{2}{c}{$c=0$} & \multicolumn{2}{c}{$c=0.15$} & \multicolumn{2}{c}{$c=0.3$} \\
		$N$   & $\VJ{\| \widetilde{\phi}^h\|_{H\bb{\widetilde{\mathrm{grad}};\Omega}}}$ & $\left\|\boldsymbol{q}^{h}\right\|_{\varun{H\bb{\mathrm{div};\Omega}}}$ & $\VJ{\| \widetilde{\phi}^h\|_{H\bb{\widetilde{\mathrm{grad}};\Omega}}}$ & $\left\|\boldsymbol{q}^{h}\right\|_{\varun{H\bb{\mathrm{div};\Omega}}}$ & $\VJ{\| \widetilde{\phi}^h\|_{H\bb{\widetilde{\mathrm{grad}};\Omega}}}$ & $\left\|\boldsymbol{q}^{h}\right\|_{\varun{H\bb{\mathrm{div};\Omega}}}$ \\
		\hline
		2  & 2.45180494&	2.45180494&	2.45180494&	2.45180494&	 2.45180494&	2.45180494 \\
		4  & 2.37137238&	2.37137238&	2.35503380&	2.35503380&	 2.13797018&	2.13797018 \\
		6  & 2.35794814&	2.35794814&	2.35666554&	2.35666554&	 2.34310363&	2.34310363 \\
		8  & 2.35588158&	2.35588158&	2.35547353&	2.35547353&	 2.35133906&	2.35133906 \\
		10 & 2.35564418&	2.35564418&	2.35556015&	2.35556015&	 2.35443148&	2.35443148 \\
		12 & 2.35561580&	2.35561580&	2.35560124&	2.35560124&	 2.35534845&	2.35534845 \\
		14 & 2.35561268&	2.35561268&	2.35561045&	2.35561045&	 2.35555229&	2.35555229 \\
		16 & 2.35561231&	2.35561231&	2.35561199&	2.35561199&	 2.35559831&	2.35559831 \\
		18 & 2.35561227&	2.35561227&	2.35561223&	2.35561223&	 2.35560913&	2.35560913 \\
		\hline
	\end{tabular}%
	\label{tab:addlabel}%
\end{table}%
In order to corroborate that the norms $\VJ{\| \widetilde{\phi}^h\|_{H\bb{\widetilde{\mathrm{grad}};\Omega}}}$ and $\|\bm{q}^h\|_{\varun{H\bb{\mathrm{div};\Omega}}}$ are identical {for this problem}, Table~\ref{tab:addlabel} lists these norms on three different meshes, the orthogonal mesh, $c=0.0$, the slightly deformed mesh, $c=0.15$ and the highly deformed mesh, $c=0.3$.
This table shows that on all mesh configurations and for all polynomial degrees we have $\VJ{\| \widetilde{\phi}^h\|_{H\bb{\widetilde{\mathrm{grad}};\Omega}}} = \|\bm{q}^h\|_{\varun{H\bb{\mathrm{div};\Omega}}}$.
All the three mesh configurations show convergence to a limiting value of $2.35561$.

\section{Eigenvalue problems for the vector Laplace operator} \label{sec:eigenvalues2}
In this section we consider the eigenvalue problem of the vector Laplace operator given by: Let $\Omega \subset \mathbb{R}^2$ and find $\lambda \in \mathbb{R}$ and a non-vanishing $\bm{u}$ such that
\begin{equation}
\begin{array}{llll}
- \mbox{grad}\ \mbox{div}\, \bm{u} &= \lambda \bm{u} \quad &\mbox{in} & \Omega \\
\mbox{div}\, \bm{u} &=0& \mbox{on} & \partial \Omega
\end{array} \;.
\label{eq:Laplace_eigenvalue_primal}
\end{equation}
This eigenvalue problem is essentially the same as the eigenvalue problem discussed in \cite{2006Boffi} apart from the fact that here Neumann boundary conditions are applied. For a more thorough treatment of eigenvalue problems, we refer to \cite{2010Boffi2}.

The corresponding weak formulation reads: Find $\lambda \in \mathbb{R}$ and $\bm{u} \in H(\mbox{div};\Omega)$ such that
\begin{equation}
( \mbox{div}\ \bm{v},\mbox{div} \ \bm{u}) = \lambda (\bm{v},\bm{u}) \;,\quad \quad \forall \bm{v} \in H(\mbox{div};\Omega)\;.
\end{equation}
For the finite dimensional approximation of this problem, we choose $\bm{u}^h \in D(\Omega)$ which then leads to the generalized matrix eigenvalue problem
\begin{equation} \label{eq:alg_egig_primal}
{\mathbb{E}^{2,1}}^T \mathbb{M}^{(2)} \mathbb{E}^{2,1} \mathcal{N}^1(\bm{u}^h) = \lambda \mathbb{M}^{(1)} \mathcal{N}^1(\bm{u}^h) \;.
\end{equation}

Alternatively, we could rephrase this problem in mixed formulation as: Find $\lambda \in \mathbb{R}$ and $(\bm{u},p)\in H(\mbox{div};\Omega)\times L^2(\Omega)$ such that
\begin{equation}
\left \{ \begin{array}{lll} 
(\bm{v},\bm{u}) + (\mbox{div}\ \bm{v},p) & = 0 \quad \quad & \forall \bm{v} \in H(\mbox{div};\Omega) \\ [1ex]
(q,\mbox{div}\ \bm{u}) & = \lambda (q,p) & \forall q \in L^2(\Omega)
\end{array} \right. \;.
\label{eq:Laplace_eigenvalue_mixed}
\end{equation}
The main difference between (\ref{eq:Laplace_eigenvalue_primal}) and (\ref{eq:Laplace_eigenvalue_mixed}) is that (\ref{eq:Laplace_eigenvalue_primal}) gives a large number of eigenvalues $\lambda = 0$, because for any $\psi \in H(\mbox{curl};\Omega)$, $\bm{u} = \mbox{curl}\ \psi$ is an eigensolution associated to $\lambda=0$.
The mixed formulation (\ref{eq:Laplace_eigenvalue_mixed}) has no eigenvalue $\lambda = 0$.

Note that the mixed formulation is essentially an eigenvalue problem for $\mbox{div}\ \bm{u}$ and not $\bm{u}$.
In the finite dimensional setting (\ref{eq:Laplace_eigenvalue_mixed}) we take $\bm{u}^h \in \VJ{D}(\Omega)$ and $p^h \in \VJ{S}(\Omega)$ and this translates into the eigenvalue problem for the degrees of freedom
\begin{equation} \label{eq:mixed_eigen}
\left ( \begin{array}{cc}
\mathbb{M}^{(1)}  &  {\mathbb{E}^{2,1}}^T \mathbb{M}^{(2)} \\
\mathbb{M}^{(2)} \mathbb{E}^{2,1} & 0
\end{array} \right )
\left ( \begin{array}{c}
\mathcal{N}^1(\bm{u}^h) \\
\mathcal{N}^2(p^h)
\end{array} \right ) = \lambda \left ( \begin{array}{cc}
0 & 0 \\
0 & \mathbb{M}^{(2)}
\end{array} \right ) \left ( \begin{array}{c}
\mathcal{N}^1(\bm{u}^h) \\
\mathcal{N}^2(p^h)
\end{array} \right ) \;.
\end{equation}
\VJ{The eigenvalues and eigenmodes for this problem were obtained directly as an output of the standard MATLAB function.}

If we eliminate the degrees of freedom $\mathcal{N}^1(\bm{u}^h)$ \VJ{(using Schur complement method)} from \eqref{eq:mixed_eigen} and rephrase the resulting eigenvalue problem in terms of the dual degrees of freedom for $p^h$ by using $\widetilde{\mathcal{N}}^0(p^h) = \mathbb{M}^{(2)} \mathcal{N}^2(p^h)$, we obtain the dual formulation: Find $\lambda \in \mathbb{R}$ and $p^h \in \widetilde{S}^h(\Omega)$ such that
\begin{equation} \label{eq:discrete_eigen}
\mathbb{E}^{2,1} \widetilde{\mathbb{M}}^{1} {\mathbb{E}^{2,1}}^T \widetilde{\mathcal{N}}^0(p^h) = \lambda \widetilde{\mathbb{M}}^{(0)} \widetilde{\mathcal{N}}^0(p^h) \;.
\end{equation}
\begin{figure}[!htb]
	\centering
	\begin{subfigure}[b]{0.3\textwidth}
		\includegraphics[width=\textwidth]{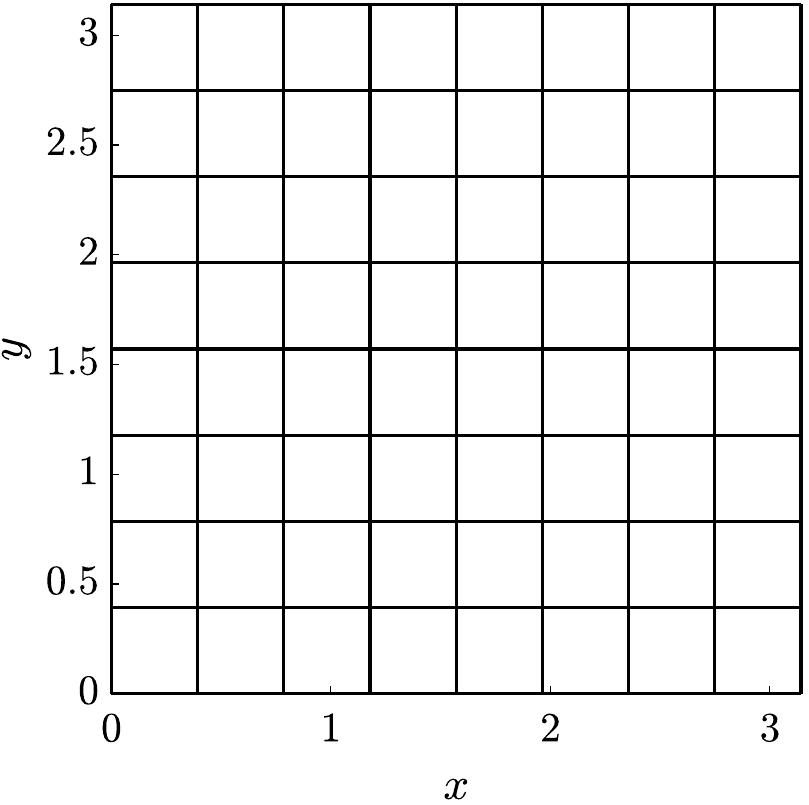}
	\end{subfigure}%
	~ 
	\begin{subfigure}[b]{0.3\textwidth}
		\includegraphics[width=\textwidth]{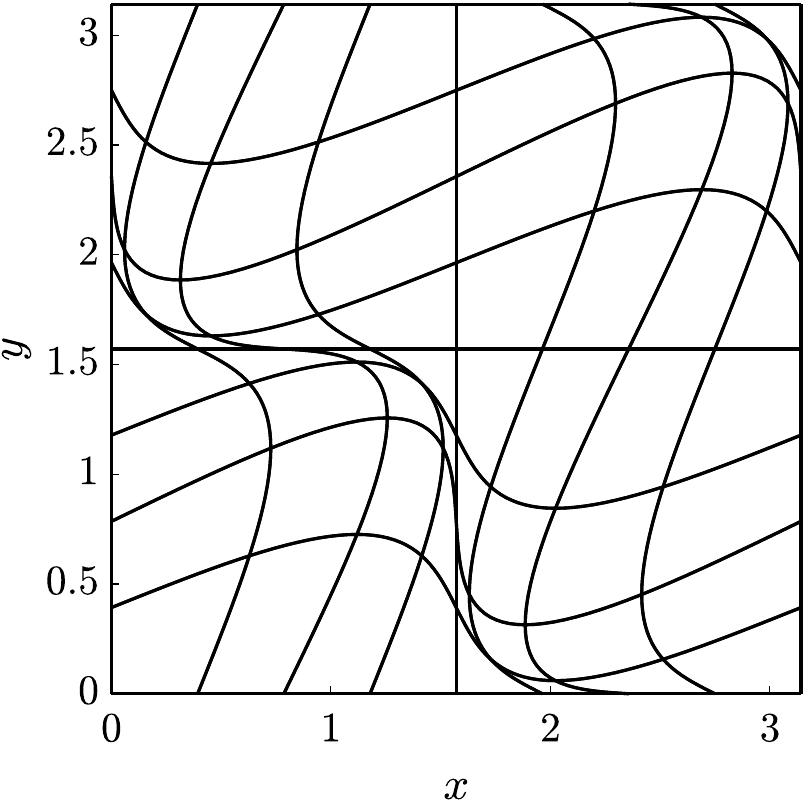}
	\end{subfigure}
	~
	\begin{subfigure}[b]{0.3\textwidth}
		\includegraphics[width=\textwidth]{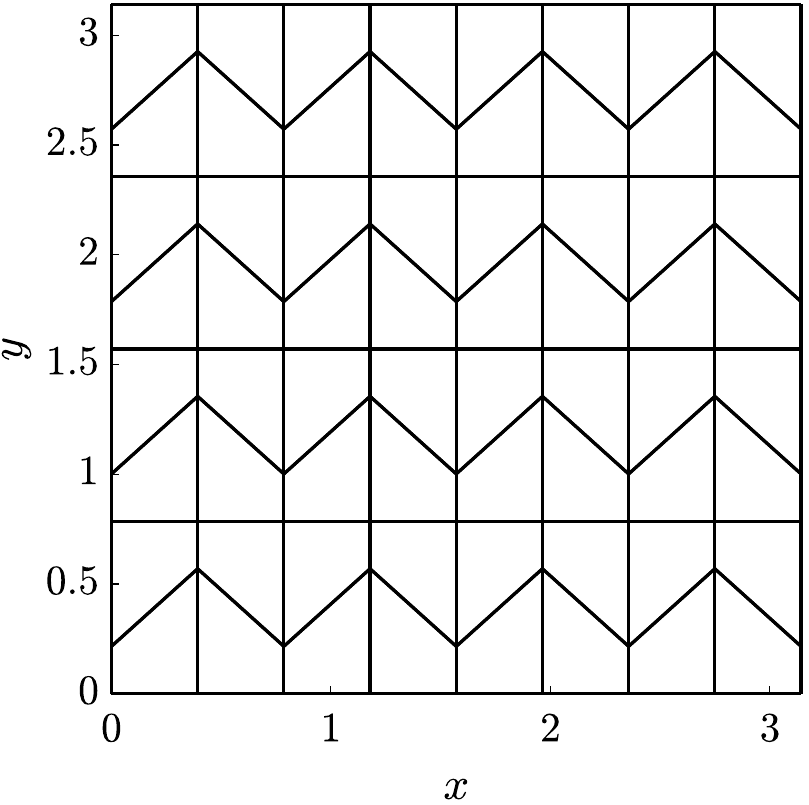}
	\end{subfigure}
	\caption{The three meshes on domain $[0,\pi]$ with $8 \times 8$ elements. Left : orthogonal mesh; Middle : curved mesh; Right : non-affine mesh} \label{fig:three_meshes}
\end{figure}
\begin{figure}[!htb]
	\centering
	\begin{subfigure}[b]{0.3\textwidth}
		\includegraphics[width=\textwidth]{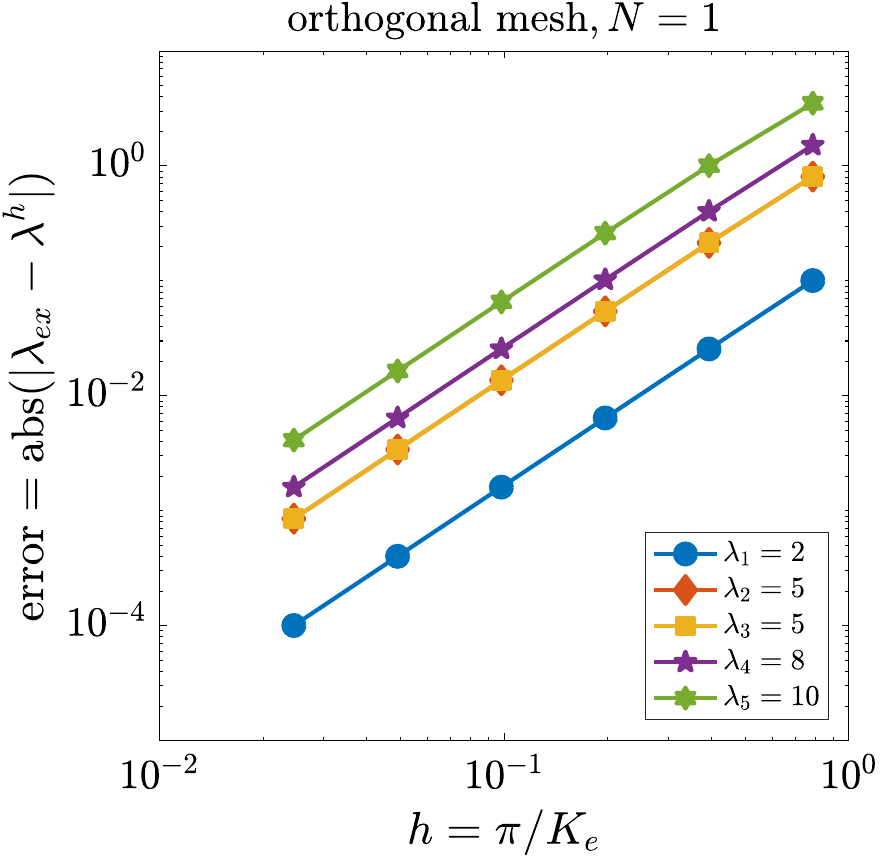}
	\end{subfigure}%
	~ 
	\begin{subfigure}[b]{0.3\textwidth}
		\includegraphics[width=\textwidth]{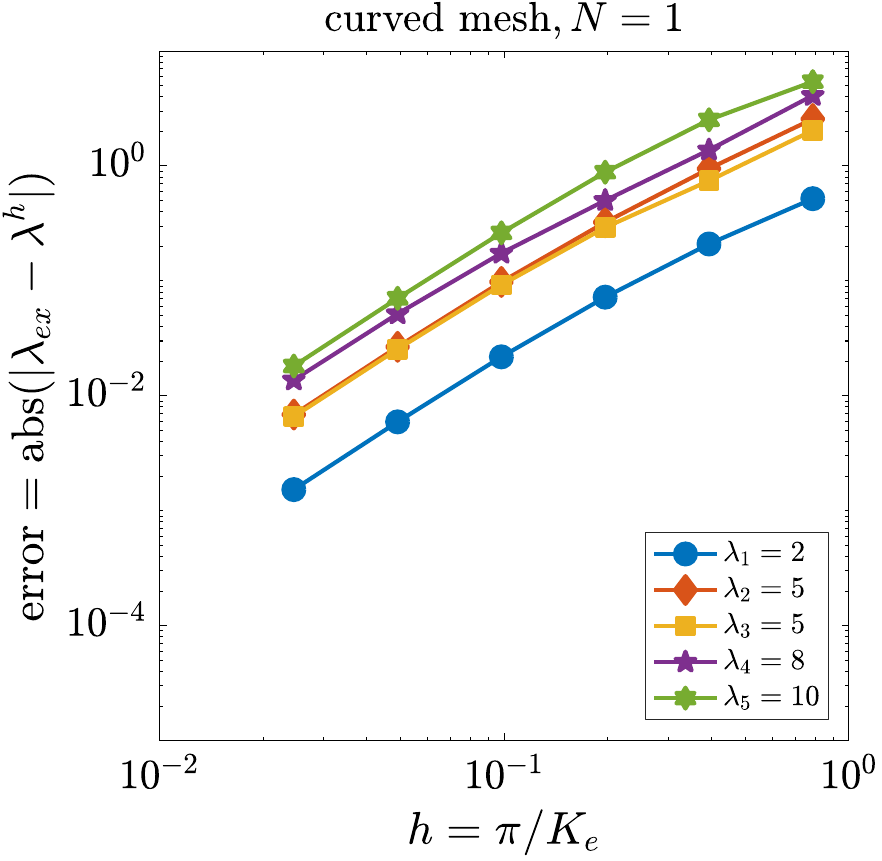}
	\end{subfigure}
	~
	\begin{subfigure}[b]{0.3\textwidth}
		\includegraphics[width=\textwidth]{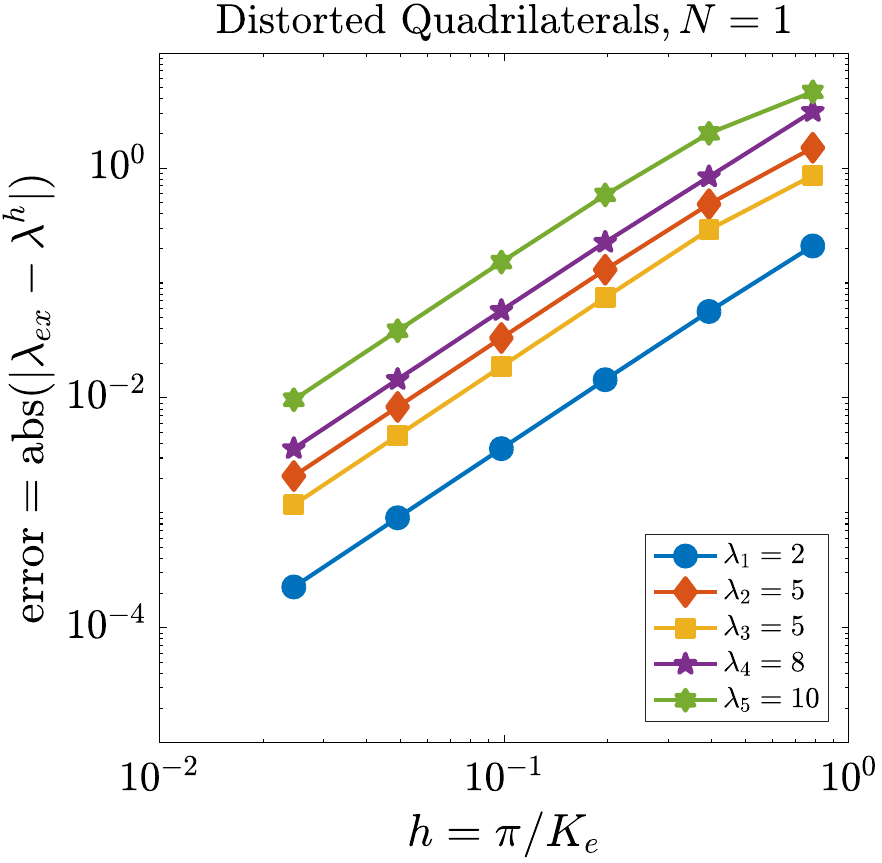}
	\end{subfigure} 
	
	\begin{subfigure}[b]{0.3\textwidth}
		\includegraphics[width=\textwidth]{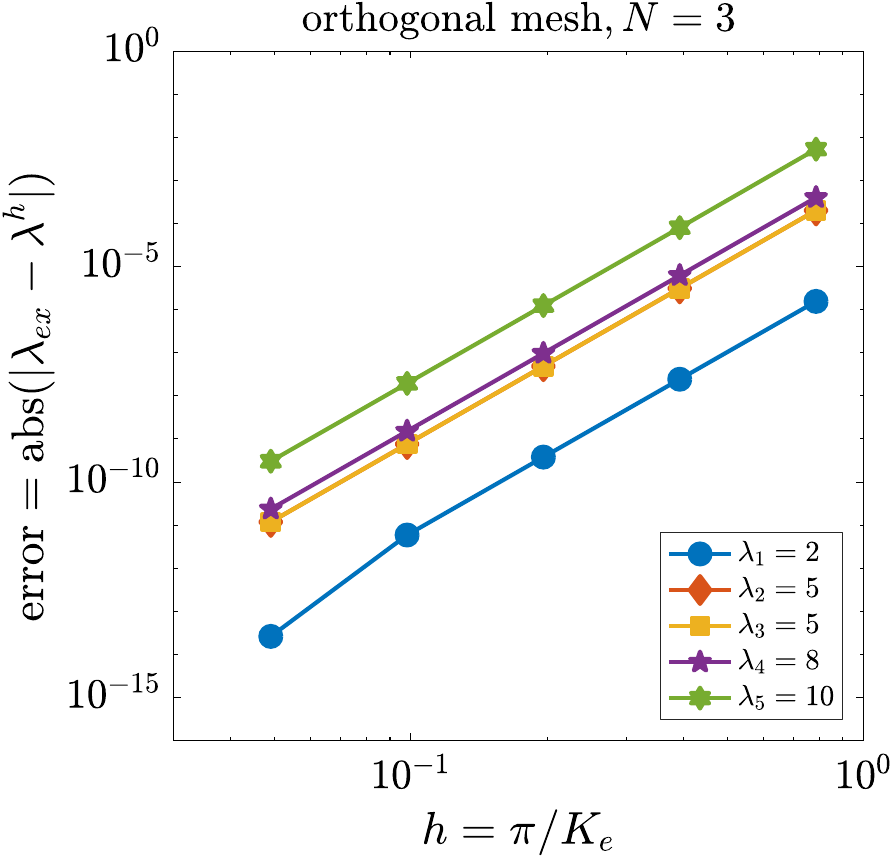}
	\end{subfigure}%
	~ 
	\begin{subfigure}[b]{0.3\textwidth}
		\includegraphics[width=\textwidth]{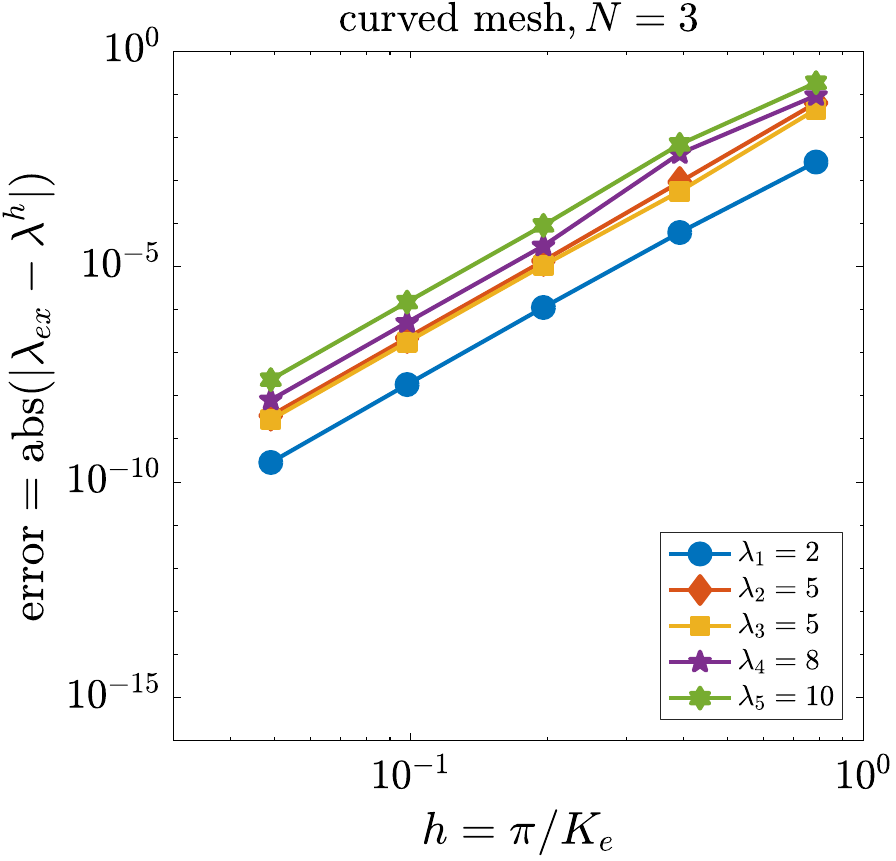}
	\end{subfigure}
	~
	\begin{subfigure}[b]{0.3\textwidth}
		\includegraphics[width=\textwidth]{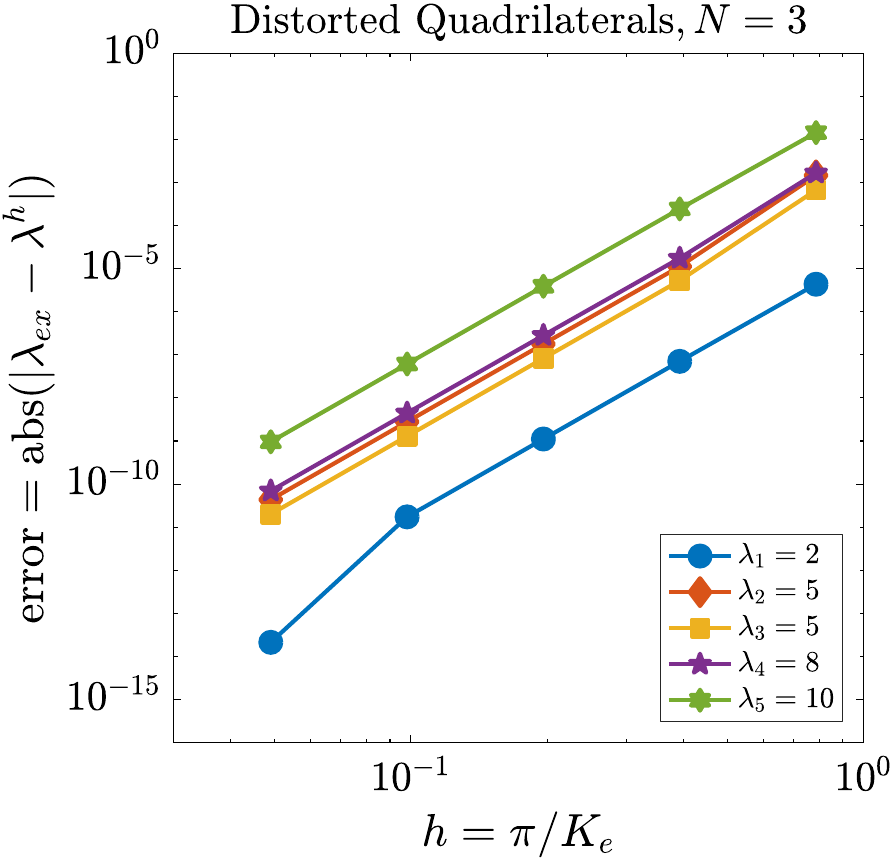}
	\end{subfigure} 
	
	\begin{subfigure}[b]{0.3\textwidth}
		\includegraphics[width=\textwidth]{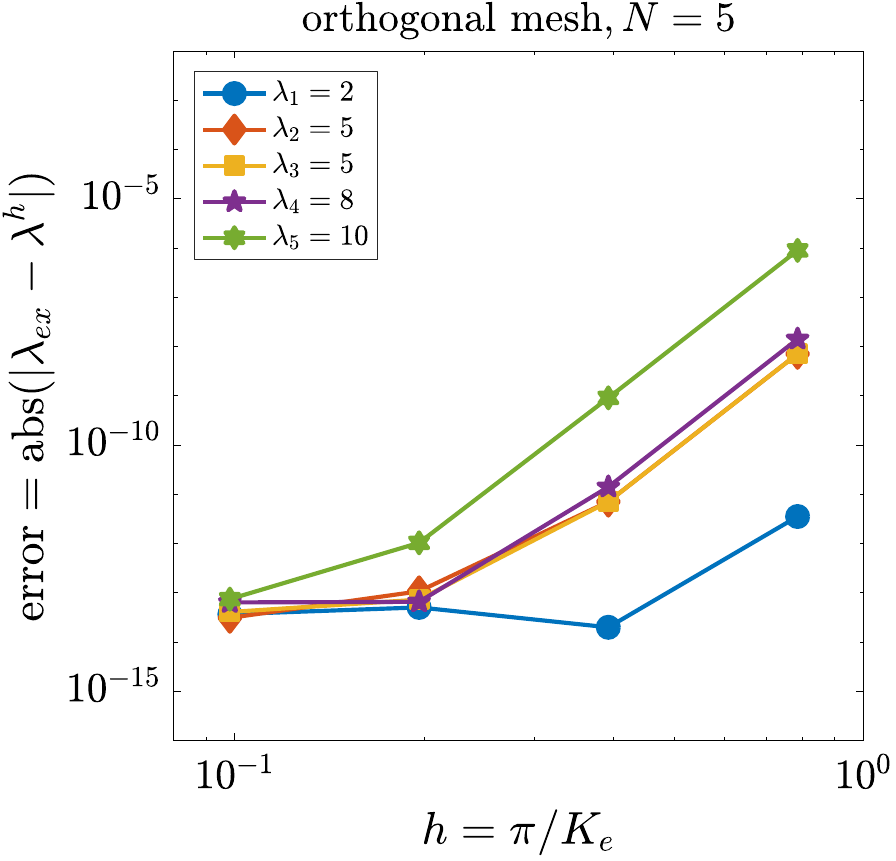}
	\end{subfigure}%
	~ 
	\begin{subfigure}[b]{0.3\textwidth}
		\includegraphics[width=\textwidth]{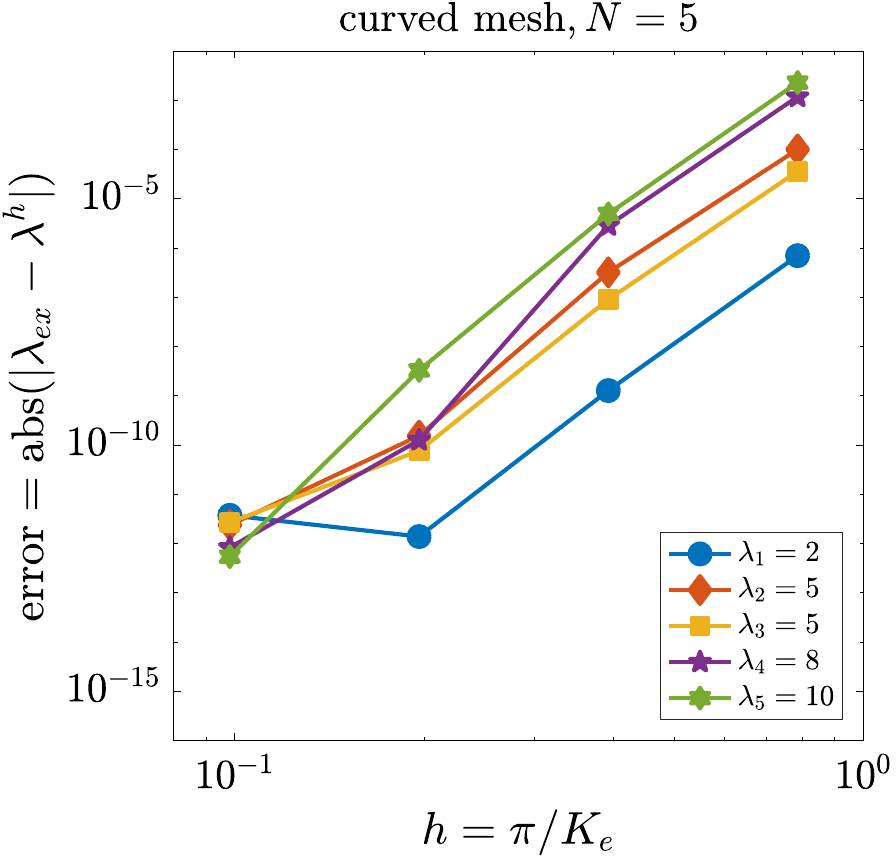}
	\end{subfigure}
	~
	\begin{subfigure}[b]{0.3\textwidth}
		\includegraphics[width=\textwidth]{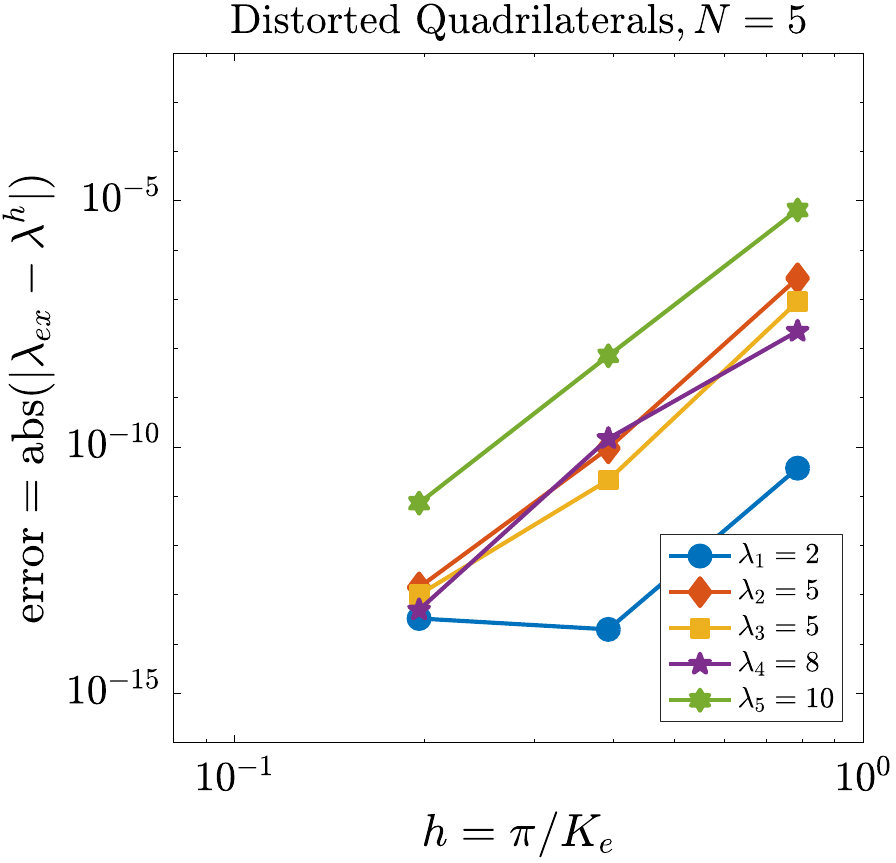}
	\end{subfigure}
	\caption{Convergence results for the first five eigenvalues. In the first column we have the results for the orthogonal mesh, in the second column for the curved mesh, and in the third column for the non-affine mesh.  In the first row we have the results for $N=1$, in the second row for $N=3$, and in the third row for $N=5$.}
	\label{fig:convergence_eigenvalue}
\end{figure}

{\begin{table}[!hbt]
		\centering
		\caption{Rate of convergence for the first five eigenvalues on 1) orthogonal mesh; 2) curved mesh; and 3) non-affine mesh for $N=1,3,5$.}
		\begin{tabular}{cccccccccc}
			\hline
			& \multicolumn{3}{c}{Orthogonal}              & \multicolumn{3}{c}{Curved} & \multicolumn{3}{c}{Non-affine}            \\
			\hline
			& \multicolumn{1}{c}{N=1} & N=3    & N=5     & N=1     & N=3     & N=5    & \multicolumn{1}{c}{N=1} & N=3    & N=5     \\
			\hline
			optimal rate  	& 2 & 6 & 10 & 2 & 6 & 10 & 2 & 6 & 10 \\
			\hline
			eigenvalue 1 & 1.9999                  & 6.0072 & 7.4718 & 1.9577  & 6.0060 & 9.4727 & 1.9996                  & 6.0019 & 10.8722 \\
			eigenvalue 2 & 1.9997                  & 6.0046 & 9.9669 & 1.9571  & 5.9851 & 9.6523  & 1.9983                  & 6.0052 & 10.4249 \\
			eigenvalue 3 & 1.9997                  & 6.0070 & 9.9692 & 1.9503  & 5.9807 & 9.4252  & 1.9994                  & 6.0102 & 9.8826  \\
			eigenvalue 4 & 1.9997                  & 6.0021 & 9.9717 & 1.9117  & 5.9741 & 10.1329  & 1.9982                  & 5.9998 & 9.3962 \\
			eigenvalue 5 & 1.9992                  & 6.0002 & 9.9458 & 1.9582  & 5.9777 & 9.7031 & 1.9961                  & 5.9978 & 9.8859  \\
			\hline
		\end{tabular}
		\label{tab:rate_of _convergence}
\end{table}}

The eigenvalues calculated for \eqref{eq:alg_egig_primal} and \eqref{eq:discrete_eigen} are identical up to 4 decimal places, therefore we only present results for \eqref{eq:discrete_eigen}.
We use three different mesh configurations 1) orthogonal mesh; 2) curved mesh; 3) non-affine mesh, as shown in \figref{fig:three_meshes}.

In \figref{fig:convergence_eigenvalue} we show the convergence plots for the first five eigenvalues of \eqref{eq:discrete_eigen} given by $\lambda = 2,5,5,8,10$.
On the x-axis we have the element length $h=\pi / K_e$, where $K_e$ is the number of elements in one direction.
On the y-axis we have the absolute error of the eigenvalues.
In \figref{fig:convergence_eigenvalue} along the horizontal rows we have the convergence plots with varying mesh: orthogonal, curved and non-affine, and in that order.
Along the vertical columns we have the convergence plots with varying $N=1,3,5$ and in that order.
The numerical eigenvalues are given in ~\ref{sec:eigenvalues}.
The rate of convergence for the eigenvalues with varying meshes and polynomial degrees $N$ are given in \tabref{tab:rate_of _convergence}.

For an order $N$ element the optimal convergence rate is given by $2N$, ~\cite{2010Boffi2}.
For all the three mesh configurations we observe optimal order convergence rates as can be seen in \tabref{tab:rate_of _convergence}.
\section{Conclusions}\label{sec:Conlusions}
In this paper a dual polynomial \VJ{space} is constructed.
The duality pairing between variables from a primal and a dual representation reduces to the vector product between the primal and dual degrees of freedom.
The grad, curl and div \varun{operators} applied to the dual \VJ{(and primal)} representation are topological relations between the dual degrees \VJ{(and primal degrees)} of freedom. 
These topological relations do not depend on the metric.
Furthermore, the addition of differential operations in the boundary ensures that the dual function spaces form a de Rham sequence.
The first example where the use of a dual representation is beneficial concerns the mixed formulation of the Poisson problem.
We derive the \emph{inf-sup} stability condition in terms of degrees of freedom {only}.
We show optimal convergence for multi-element test case on a curved 3D domain.
When a primal-dual formulation is used, two sub-matrices in the mixed formulation become very sparse, even though very high order methods are used and these two sub-matrices do not change when the mesh is deformed.
The second example shows the equivalence of a Dirichlet-Neumann pair of equations at the discrete level.
This equivalence is proven and illustrated with a test case.
And lastly, we \VJ{solve the grad-div} eigenvalue problem in terms of primal and dual representations, and show optimal convergence on both affine and non-affine meshes.

We have also seen that the use of dual representations allows us to work directly with the degrees of freedom, without explicitly referring to the basis functions.
It suffices to make use of its properties.
This allows for a direct comparison with staggered finite volume methods.
For \varun{example}, in (\ref{eq:primal_dual_system}), $\mathbb{E}^{3,2}$ acts directly on the degrees of freedom of $\bm{q}^h$ and ${\mathbb{E}^{3,2}}^T$ acts on the dual degrees of freedom for $\phi^h$.

\VJ{To deal with computationally demanding problems, in future, this work will be extended to the framework of domain decomposition methods. The compatible construction of the discrete trace space presented in this work allows one to set up a fully compatible hybrid finite element formulation. The advantage of a hybrid formulation is that the dual basis functions and degrees of freedom can be defined using element mass matrices only, instead of the global mass matrices presented in this paper. For preliminary work in this direction, see \cite{2018Jain}.}

	\VJ{
	Furthermore, in this paper the construction of dual polynomial spaces is based on multiplication with mass matrices (or inverse of mass matrices).
	These matrices change with change in shape and size of the element.
	As an improvement, we will present a construction of dual spaces where the mass matrix is also independent of the shape and the size of the element by using wedge product instead of the inner product, see \cite[\S6.1]{Palha2014}.}

\section*{Acknowledgements}
Varun Jain is supported by the Shell-NWO PhD75 Grant P/fom-d-69/3.
Yi Zhang is supported by Chinese Scholarship Council No. 201607720010. \MIGnew{The authors would like to thank the reviewers for their valuable suggestions which improved the paper in many ways.}

\bibliographystyle{elsarticle-num}
\bibliography{references}

\appendix
\section{Eigenvalues} \label{sec:eigenvalues}
\begin{table}[!hbt] 
	\caption{eigenvalue 1}
	\begin{tabular}{cccccccc}
		\hline
		& h    & $\pi / 4$        & $\pi / 8$        & $\pi / 16$        & $\pi / 32$        & $\pi / 64$        & $\pi / 128$        \\
		\hline
		Cartesian   & N=1 & 1.8993 & 1.9744 & 1.9936 & 1.9984 & 1.9996 & 1.9999 \\
		& N=3 & 2.0000 & 2.0000 & 2.0000 & 2.0000 & 2.0000 &        \\
		& N=5 & 2.0000 & 2.0000 & 2.0000 & 2.0000 &        &        \\
		\hline
		Curved      & N=1 & 1.4808 & 1.7905 & 1.9278 & 1.9782 & 1.9941 & 1.9985 \\
		& N=3 & 1.9973 & 1.9999 & 2.0000 & 2.0000 & 2.0000 &        \\
		& N=5 & 2.0000 & 2.0000 & 2.0000 & 2.0000 &        &        \\
		\hline
		Non-affine & N=1 & 1.7901 & 1.9435 & 1.9856 & 1.9964 & 1.9991 & 1.9998 \\
		& N=3 & 2.0000 & 2.0000 & 2.0000 & 2.0000 & 2.0000 &        \\
		& N=5 & 2.0000 & 2.0000 & 2.0000 & 2.0000 &        &       \\
		\hline
	\end{tabular}
	\label{tab:eigenvalue1}
\end{table}
\begin{table}
	\caption{eigenvalue2}
	\begin{tabular}{cccccccc}
		\hline
		& h    & $\pi / 4$        & $\pi / 8$        & $\pi / 16$        & $\pi / 32$        & $\pi / 64$        & $\pi / 128$        \\
		\hline
		Cartesian   & N=1 & 4.1919 & 4.7858 & 4.9457 & 4.9864 & 4.9966 & 4.9991 \\
		& N=3 & 4.9998 & 5.0000 & 5.0000 & 5.0000 & 5.0000 &        \\
		& N=5 & 5.0000 & 5.0000 & 5.0000 & 5.0000 &        &        \\
		\hline
		Curved      & N=1 & 2.4366 & 4.0535 & 4.6764 & 4.9020 & 4.9733 & 4.9931 \\
		& N=3 & 4.9374 & 4.9991 & 5.0000 & 5.0000 & 5.0000 &        \\
		& N=5 & 4.9999 & 5.0000 & 5.0000 & 5.0000 &        &        \\
		\hline
		Non-affine & N=1 & 3.4977 & 4.5154 & 4.8696 & 4.9668 & 4.9917 & 4.9979 \\
		& N=3 & 4.9985 & 5.0000 & 5.0000 & 5.0000 & 5.0000 &        \\
		& N=5 & 5.0000 & 5.0000 & 5.0000 &        &        &                \\
		\hline
	\end{tabular}
\end{table}
\begin{table}
	\caption{eigenvalue3}
	\begin{tabular}{cccccccc}
		\hline
		& h    & $\pi / 4$        & $\pi / 8$        & $\pi / 16$        & $\pi / 32$        & $\pi / 64$        & $\pi / 128$        \\
		\hline
		Cartesian   & N=1 & 4.1919 & 4.7858 & 4.9457 & 4.9864 & 4.9966 & 4.9991 \\
		& N=3 & 4.9998 & 5.0000 & 5.0000 & 5.0000 & 5.0000 &        \\
		& N=5 & 5.0000 & 5.0000 & 5.0000 & 5.0000 &        &        \\
		\hline
		Curved      & N=1 & 2.9601 & 4.2569 & 4.7085 & 4.9078 & 4.9745 & 4.9934 \\
		& N=3 & 4.9558 & 4.9994 & 5.0000 & 5.0000 & 5.0000 &        \\
		& N=5 & 5.0000 & 5.0000 & 5.0000 & 5.0000 &        &        \\
		\hline
		Non-affine & N=1 & 4.1387 & 4.7093 & 4.9255 & 4.9812 & 4.9953 & 4.9988 \\
		& N=3 & 4.9993 & 5.0000 & 5.0000 & 5.0000 & 5.0000 &        \\
		& N=5 & 5.0000 & 5.0000 & 5.0000 & 5.0000 &        &         \\
		\hline
	\end{tabular}
\end{table}
\begin{table}[!hbt]
	\caption{eigenvalue4}
	\begin{tabular}{cccccccc}
		\hline
		& h    & $\pi / 4$        & $\pi / 8$        & $\pi / 16$        & $\pi / 32$        & $\pi / 64$        & $\pi / 128$        \\
		\hline
		Cartesian   & N=1 & 6.4846 & 7.5971 & 7.8977 & 7.9743 & 7.9936 & 7.9984 \\
		& N=3 & 7.9996 & 8.0000 & 8.0000 & 8.0000 & 8.0000 &        \\
		& N=5 & 8.0000 & 8.0000 & 8.0000 & 8.0000 &        &        \\
		\hline
		Curved      & N=1 & 3.8951 & 6.6238 & 7.4993 & 7.8256 & 7.9482 & 7.9862 \\
		& N=3 & 7.9074 & 7.9958 & 8.0000 & 8.0000 & 8.0000 &        \\
		& N=5 & 8.0012 & 8.0000 & 8.0000 & 8.0000 &        &        \\
		\hline
		Non-affine & N=1 & 4.8904 & 7.1603 & 7.7740 & 7.9424 & 7.9855 & 7.9964 \\
		& N=3 & 7.9983 & 8.0000 & 8.0000 & 8.0000 & 8.0000 &        \\
		& N=5 & 8.0000 & 8.0000 & 8.0000 &        &        &           \\
		\hline
	\end{tabular}
\end{table}
\begin{table}[!hbt]
	\caption{eigenvalue5}
	\begin{tabular}{cccccccc}
		\hline
		& h    & $\pi / 4$        & $\pi / 8$        & $\pi / 16$        & $\pi / 32$        & $\pi / 64$        & $\pi / 128$        \\
		\hline
		Cartesian   & N=1 & 6.4846  & 8.9933  & 9.7395  & 9.9343  & 9.9835  & 9.9959 \\
		& N=3 & 9.9947  & 9.9999  & 10.0000 & 10.0000 & 10.0000 &        \\
		& N=5 & 10.0000 & 10.0000 & 10.0000 & 10.0000 &         &        \\
		\hline
		Curved      & N=1 & 4.5854  & 7.4787  & 9.1173  & 9.7381  & 9.9294  & 9.9818 \\
		& N=3 & 9.8156  & 9.9932  & 9.9999  & 10.0000 & 10.0000 &        \\
		& N=5 & 9.9977  & 10.0000 & 10.0000 & 10.0000 &         &        \\
		\hline
		Non-affine & N=1 & 5.4005  & 7.9946  & 9.4163  & 9.8479  & 9.9616  & 9.9904 \\
		& N=3 & 9.9857  & 9.9998  & 10.0000 & 10.0000 & 10.0000 &        \\
		& N=5 & 10.0000 & 10.0000 & 10.0000 &         &         &           \\
		\hline
	\end{tabular}
\end{table}
\clearpage

\varun{\section{Identity at domain boundary}\label{section:motivation}
		For an $n$-dimensional domain we wish to prove the following identity
		\begin{equation}
			\mathbb{N}_{k-1}\mathbb{N}_{k-1}^{\intercal}\left(\mathbb{E}^{k, k-1}\right)^{\intercal}\mathbb{N}_{k} = \left(\mathbb{E}^{k, k-1}\right)^{\intercal}\mathbb{N}_{k}\,, \qquad \text{with } 0 < k < n\,.
\end{equation}
To prove the identity we first need to identify some properties of the different matrices. The first point to note is that by changing the numbering of the degrees of freedom we simply permute the rows and columns of the matrices. Naturally, this change in the numbering of the degrees of freedom does not impact the identity. For this reason we are free to choose any numbering that best suites the proof, and the result will hold for any other numbering. Therefore, in this proof we choose the following numbering. For each set of degrees of freedom (associated to geometrical objects of dimension $k$, with $0 \leq k \leq n$) we first number the degrees of freedom on the boundary and then the ones on the interior of the domain. Therefore, for each set of degrees of freedom associated to a $k$-dimensional geometric object, the boundary degrees of freedom are numbered from 1 to $d^{b}_{k}$ and the ones in the interior from $d^{b}_{k} + 1$ to $d_{k}$, where $d^{b}_{k}$ is the number of degrees of freedom on the boundary and $d_{k}$ is the total number of degrees of freedom.
}

\varun{If we make this choice of numbering, then $\mathbb{N}_{k}$ is a $d_{k}\times d_{k}^{b}$ matrix with the following form
		\begin{equation}
			\mathbb{N}_{k} = 
			\left[
				\begin{array}{c}
					D_{k} \\
					\hline
					0
				\end{array}
			\right]\,, \label{eq:n_matrix_property}
		\end{equation} 
		where $D_{k}$ is a $d_{k}^{b}\times d_{k}^{b}$ diagonal matrix filled with only -1 and 1. In the same way the incidence matrix $\mathbb{E}^{k, k-1}$ is a $d_{k}\times d_{k-1}$ matrix given as
		\begin{equation}
			\mathbb{E}^{k,k-1} = 
			\left[
				\begin{array}{c|c}
					\mathbb{E}_{bb}^{k, k-1} & 0 \\
					\hline
					\mathbb{E}_{bi}^{k, k-1} & \mathbb{E}_{ii}^{k, k-1}
				\end{array}
			\right]\,, \label{eq:incidence_matrix}
		\end{equation} 
		where $\mathbb{E}_{bb}^{k, k-1}$ contains the terms of the incidence matrix that only relate boundary degrees of freedom, $\mathbb{E}_{bi}^{k, k-1}$ contains the contribution of the boundary degrees of freedom to the interior, and finally $\mathbb{E}_{ii}^{k, k-1}$ contains the contributions that only relate to the interior degrees of freedom. What is important to note here is the zero block that states that the interior degrees of freedom do not contribute to the boundary degrees of freedom, this is essential for the proof.}
		
\varun{If we now use \eqref{eq:n_matrix_property} we can easily see that
		\begin{equation}
			\mathbb{N}_{k}\mathbb{N}_{k}^{\intercal} = 
			\left[
				\begin{array}{c}
					D_{k} \\
					\hline
					0
				\end{array}
			\right]
				\left[
				\begin{array}{c|c}
					D_{k} & 0
				\end{array}
				\right]
				= 
				\left[
				\begin{array}{c|c}
					\mathbb{I}_{k}^{b} & 0 \\
					\hline
					0 & 0
				\end{array}
				\right]\,. \label{eq:n_matrix_identity}
		\end{equation}
		where $\mathbb{I}_{k}^{b}$ is the $d^{b}_{k}\times d^{b}_{k}$ identity matrix.}
		
\varun{Combining \eqref{eq:n_matrix_property} and \eqref{eq:incidence_matrix} we get
		\begin{equation}
			\left(\mathbb{E}^{k,k-1}\right)^{\intercal} \mathbb{N}_{k} = 
			\left[
				\begin{array}{c|c}
					\left(\mathbb{E}_{bb}^{k, k-1}\right)^{\intercal} & \left(\mathbb{E}_{bi}^{k, k-1}\right)^{\intercal} \\
					\hline
					0 & \left(\mathbb{E}_{ii}^{k, k-1}\right)^{\intercal}
				\end{array}
			\right]
			\left[
				\begin{array}{c}
					D_{k} \\
					\hline
					0
				\end{array}
			\right]
			=
			\left[
				\begin{array}{c}
					\left(\mathbb{E}_{bb}^{k, k-1}\right)^{\intercal} D_{k} \\
					\hline
					0
				\end{array}
			\right]\,. \label{eq:incidence_n}
		\end{equation}}
		
\varun{Using \eqref{eq:n_matrix_identity} together with \eqref{eq:incidence_n}, it is straight forward to see that
		\begin{equation}
			\mathbb{N}_{k-1}\mathbb{N}_{k-1}^{\intercal} \left(\mathbb{E}^{k,k-1}\right)^{\intercal} \mathbb{N}_{k}  = 
			\left[
				\begin{array}{c|c}
					\mathbb{I}_{k-1}^{b} & 0 \\
					\hline
					0 & 0
				\end{array}
				\right]
				\left[
				\begin{array}{c}
					\left(\mathbb{E}_{bb}^{k, k-1}\right)^{\intercal} D_{k} \\
					\hline
					0
				\end{array}
			\right]
			=
			\left[
				\begin{array}{c}
					\left(\mathbb{E}_{bb}^{k, k-1}\right)^{\intercal} D_{k} \\
					\hline
					0
				\end{array}
			\right]
			= 
			\left(\mathbb{E}^{k,k-1}\right)^{\intercal} \mathbb{N}_{k}\,.
		\end{equation}}
\end{document}